\PassOptionsToPackage{unicode}{hyperref}
\PassOptionsToPackage{hyphens}{url}
\PassOptionsToPackage{dvipsnames,svgnames,x11names}{xcolor}
\documentclass[
]{article}
\usepackage{xcolor}
\usepackage{amsmath,amssymb}
\usepackage{iftex}
\ifPDFTeX
  \usepackage[T1]{fontenc}
  \usepackage[utf8]{inputenc}
  \usepackage{textcomp} % provide euro and other symbols
\else % if luatex or xetex
  \usepackage{unicode-math} % this also loads fontspec
  \defaultfontfeatures{Scale=MatchLowercase}
  \defaultfontfeatures[\rmfamily]{Ligatures=TeX,Scale=1}
\fi
\IfFileExists{upquote.sty}{\usepackage{upquote}}{}
\IfFileExists{microtype.sty}{% use microtype if available
  \usepackage[]{microtype}
  \UseMicrotypeSet[protrusion]{basicmath} % disable protrusion for tt fonts
}{}
\makeatletter
\@ifundefined{KOMAClassName}{% if non-KOMA class
  \IfFileExists{parskip.sty}{%
    \usepackage{parskip}
  }{% else
    \setlength{\parindent}{0pt}
    \setlength{\parskip}{6pt plus 2pt minus 1pt}}
}{% if KOMA class
  \KOMAoptions{parskip=half}}
\makeatother
\makeatletter
\ifx\paragraph\undefined\else
  \let\oldparagraph\paragraph
  \renewcommand{\paragraph}{
    \@ifstar
      \xxxParagraphStar
      \xxxParagraphNoStar
  }
  \newcommand{\xxxParagraphStar}[1]{\oldparagraph*{#1}\mbox{}}
  \newcommand{\xxxParagraphNoStar}[1]{\oldparagraph{#1}\mbox{}}
\fi
\ifx\subparagraph\undefined\else
  \let\oldsubparagraph\subparagraph
  \renewcommand{\subparagraph}{
    \@ifstar
      \xxxSubParagraphStar
      \xxxSubParagraphNoStar
  }
  \newcommand{\xxxSubParagraphStar}[1]{\oldsubparagraph*{#1}\mbox{}}
  \newcommand{\xxxSubParagraphNoStar}[1]{\oldsubparagraph{#1}\mbox{}}
\fi
\makeatother

\usepackage{color}
\usepackage{fancyvrb}

\DefineVerbatimEnvironment{Highlighting}{Verbatim}{commandchars=\\\{\}}
\usepackage{framed}
\definecolor{shadecolor}{RGB}{241,243,245}
\newenvironment{Shaded}{\begin{snugshade}}{\end{snugshade}}

\newcommand{\BuiltInTok}[1]{\textcolor[rgb]{0.00,0.23,0.31}{#1}}
\newcommand{\CharTok}[1]{\textcolor[rgb]{0.13,0.47,0.30}{#1}}
\newcommand{\CommentTok}[1]{\textcolor[rgb]{0.37,0.37,0.37}{#1}}

\newcommand{\ControlFlowTok}[1]{\textcolor[rgb]{0.00,0.23,0.31}{\textbf{#1}}}

\newcommand{\DecValTok}[1]{\textcolor[rgb]{0.68,0.00,0.00}{#1}}

\newcommand{\FloatTok}[1]{\textcolor[rgb]{0.68,0.00,0.00}{#1}}

\newcommand{\ImportTok}[1]{\textcolor[rgb]{0.00,0.46,0.62}{#1}}

\newcommand{\KeywordTok}[1]{\textcolor[rgb]{0.00,0.23,0.31}{\textbf{#1}}}
\newcommand{\NormalTok}[1]{\textcolor[rgb]{0.00,0.23,0.31}{#1}}
\newcommand{\OperatorTok}[1]{\textcolor[rgb]{0.37,0.37,0.37}{#1}}

\newcommand{\PreprocessorTok}[1]{\textcolor[rgb]{0.68,0.00,0.00}{#1}}

\newcommand{\SpecialCharTok}[1]{\textcolor[rgb]{0.37,0.37,0.37}{#1}}
\newcommand{\SpecialStringTok}[1]{\textcolor[rgb]{0.13,0.47,0.30}{#1}}
\newcommand{\StringTok}[1]{\textcolor[rgb]{0.13,0.47,0.30}{#1}}
\newcommand{\VariableTok}[1]{\textcolor[rgb]{0.07,0.07,0.07}{#1}}

\usepackage{longtable,booktabs,array}
\usepackage{calc} % for calculating minipage widths
\usepackage{etoolbox}
\makeatletter
\patchcmd\longtable{\par}{\if@noskipsec\mbox{}\fi\par}{}{}
\makeatother
\IfFileExists{footnotehyper.sty}{\usepackage{footnotehyper}}{\usepackage{footnote}}
\makesavenoteenv{longtable}
\usepackage{graphicx}
\makeatletter
\newsavebox\pandoc@box
\newcommand*\pandocbounded[1]{% scales image to fit in text height/width
  \sbox\pandoc@box{#1}%
  \Gscale@div\@tempa{\textheight}{\dimexpr\ht\pandoc@box+\dp\pandoc@box\relax}%
  \Gscale@div\@tempb{\linewidth}{\wd\pandoc@box}%
  \ifdim\@tempb\p@<\@tempa\p@\let\@tempa\@tempb\fi% select the smaller of both
  \ifdim\@tempa\p@<\p@\scalebox{\@tempa}{\usebox\pandoc@box}%
  \else\usebox{\pandoc@box}%
  \fi%
}
\def\fps@figure{htbp}
\makeatother

\NewDocumentCommand\citeproctext{}{}

\makeatletter
 \let\@cite@ofmt\@firstofone
 \def\@biblabel#1{}
 \def\@cite#1#2{{#1\if@tempswa , #2\fi}}
\makeatother
\newlength{\cslhangindent}
\newlength{\csllabelwidth}
\newenvironment{CSLReferences}[2] % #1 hanging-indent, #2 entry-spacing
 {\begin{list}{}{%
  \setlength{\itemindent}{0pt}
  \setlength{\leftmargin}{0pt}
  \setlength{\parsep}{0pt}
  \ifodd #1
   \setlength{\leftmargin}{\cslhangindent}
   \setlength{\itemindent}{-1\cslhangindent}
  \fi
  \setlength{\itemsep}{#2\baselineskip}}}
 {\end{list}}
\usepackage{calc}

\providecommand{\tightlist}{%
  \setlength{\itemsep}{0pt}\setlength{\parskip}{0pt}}

\usepackage{url}
\usepackage{glossaries}
\newacronym{ai}{AI}{Artificial Intelligence}
\newacronym{bo}{BO}{Bayesian Optimization}
\newacronym{cart}{CART}{Classification And Regression Tree}
\newacronym{ccd}{CCD}{Central Composite Design}
\newacronym{cnn}{CNN}{Convolutional Neural Network}
\newacronym{cpm}{CPM}{Compressor Performance Map}
\newacronym{cv}{CV}{Cross Validation}
\newacronym{cvfdt}{CVFDT}{Concept-adapting Very Fast Decision Tree}
\newacronym{dace}{DACE}{Design and Analysis of Computer Experiments}
\newacronym{ddm}{DDM}{Drift Detection Method}
\newacronym{dl}{DL}{Deep Learning}
\newacronym{doe}{DOE}{Design of Experiments}
\newacronym{efdt}{EFDT}{Extremely Fast Decision Tree}
\newacronym{gbrt}{gbrt}{Gradient Boosting Regression Tree}
\newacronym{gcd}{GCD}{Greatest Common Divisor}
\newacronym{gra}{GRA}{Global Recurring Abrupt}
\newacronym{hat}{HAT}{Hoeffding Adaptive Tree}
\newacronym{hatc}{HATC}{Hoeffding Adaptive Tree Classifier}
\newacronym{hatr}{HATR}{Hoeffding Adaptive Tree Regressor}
\newacronym{hcf}{HCF}{High-Cycle Fatigue}
\newacronym{hpt}{HPT}{Hyperparameter Tuning}
\newacronym{ht}{HT}{Hoeffding Tree}
\newacronym{htc}{HTC}{Hoeffding Tree Classifier}
\newacronym{htr}{HTR}{Hoeffding Tree Regressor}
\newacronym{ig}{IG}{Integrated Gradient}
\newacronym{ki}{KI}{Künstliche Intelligenz}
\newacronym{kpi}{KPI}{Key Performance Indicator}
\newacronym{lasso}{Lasso}{Least Absolute Shrinkage and Selection operator}
\newacronym{mae}{MAE}{Mean Absolute Error}
\newacronym{ml}{ML}{Machine Learning}
\newacronym{moa}{MOA}{Massive Online Analysis}
\newacronym{mse}{MSE}{Mean Squared Error}
\newacronym{nn}{NN}{Neural Network}
\newacronym{ocba}{OCBA}{Optimal Computational Budget Allocation}
\newacronym{pa}{PA}{Passive-Aggressive}
\newacronym{pca}{PCA}{Principal Component Analysis}
\newacronym{rf}{RF}{Random Forest}
\newacronym{rsm}{RSM}{Response Surface Methodology}
\newacronym{river}{river}{River: Online machine learning in Python}
\newacronym{rmoa}{RMOA}{Massive Online Analysis in R}
\newacronym{rocauc}{ROC AUC}{AUC (Area Under The Curve) ROC (Receiver Operating Characteristics)}
\newacronym{sea}{SEA}{SEA synthetic dataset}
\newacronym{shap}{SHAP}{SHapley Additive exPlanations}
\newacronym{sklearn}{sklearn}{scikit-learn: Machine Learning in Python}
\newacronym{smbo}{SMBO}{Surrogate Model Based Optimization}
\newacronym{smote}{SMOTE}{Synthetic Minority Oversampling Technique}
\newacronym{spo}{SPO}{Sequential Parameter Optimization}
\newacronym{spot}{SPOT}{Sequential Parameter Optimization Toolbox}
\newacronym{spotpython}{spotPython}{Sequential Parameter Optimization Toolbox for Python}
\newacronym{spotriver}{spotRiver}{Sequential Parameter Optimization Toolbox for River}
\newacronym{sgd}{SGD}{Stochastic Gradient Descent}
\newacronym{svm}{SVM}{Support Vector Machine}
\newacronym{vfdt}{VFDT}{Very Fast Decision Tree}
\newacronym{xai}{XAI}{eXplainable Artificial Intelligence}
\newacronym{xgb}{XGBoost}{eXtreme Gradient Boosting}
\newacronym{blue}{BLUE}{ResponsiBle models, Legal issues, trUst in predictions, Ethical issues}
\newacronym{red}{RED}{Research on data, Explore models, Debug models}
\usepackage{arxiv}
\usepackage{orcidlink}
\usepackage{amsmath}
\usepackage[T1]{fontenc}
\makeatletter
\@ifpackageloaded{tcolorbox}{}{\usepackage[skins,breakable]{tcolorbox}}
\@ifpackageloaded{fontawesome5}{}{\usepackage{fontawesome5}}
\definecolor{quarto-callout-color}{HTML}{909090}
\definecolor{quarto-callout-note-color}{HTML}{0758E5}
\definecolor{quarto-callout-important-color}{HTML}{CC1914}
\definecolor{quarto-callout-warning-color}{HTML}{EB9113}
\definecolor{quarto-callout-tip-color}{HTML}{00A047}
\definecolor{quarto-callout-caution-color}{HTML}{FC5300}
\definecolor{quarto-callout-color-frame}{HTML}{acacac}
\definecolor{quarto-callout-note-color-frame}{HTML}{4582ec}
\definecolor{quarto-callout-important-color-frame}{HTML}{d9534f}
\definecolor{quarto-callout-warning-color-frame}{HTML}{f0ad4e}
\definecolor{quarto-callout-tip-color-frame}{HTML}{02b875}
\definecolor{quarto-callout-caution-color-frame}{HTML}{fd7e14}
\makeatother
\makeatletter
\@ifpackageloaded{caption}{}{\usepackage{caption}}
\AtBeginDocument{%
\ifdefined\contentsname
  \renewcommand*\contentsname{Table of contents}
\else
  \newcommand\contentsname{Table of contents}
\fi
\ifdefined\listfigurename
  \renewcommand*\listfigurename{List of Figures}
\else
  \newcommand\listfigurename{List of Figures}
\fi
\ifdefined\listtablename
  \renewcommand*\listtablename{List of Tables}
\else
  \newcommand\listtablename{List of Tables}
\fi
\ifdefined\figurename
  \renewcommand*\figurename{Figure}
\else
  \newcommand\figurename{Figure}
\fi
\ifdefined\tablename
  \renewcommand*\tablename{Table}
\else
  \newcommand\tablename{Table}
\fi
}
\@ifpackageloaded{float}{}{\usepackage{float}}
\floatstyle{ruled}
\@ifundefined{c@chapter}{\newfloat{codelisting}{h}{lop}}{\newfloat{codelisting}{h}{lop}[chapter]}
\floatname{codelisting}{Listing}

\usepackage{amsthm}
\theoremstyle{definition}
\newtheorem{example}{Example}[section]
\theoremstyle{remark}
\AtBeginDocument{}

\makeatother
\makeatletter
\@ifpackageloaded{caption}{}{\usepackage{caption}}
\@ifpackageloaded{subcaption}{}{\usepackage{subcaption}}
\makeatother
\usepackage{bookmark}
\IfFileExists{xurl.sty}{\usepackage{xurl}}{} % add URL line breaks if available
\hypersetup{
  pdftitle={Multi-Objective Optimization and Hyperparameter Tuning With Desirability Functions},
  pdfauthor={Thomas Bartz-Beielstein},
  pdfkeywords={desirability function, multi-objective
optimization, surrogate modeling, hyperparameter tuning, Morris-Mitchell
criterion, maximin criterion, sequential parameter optimization},
  colorlinks=true,
  linkcolor={blue},
  filecolor={Maroon},
  citecolor={Blue},
  urlcolor={Blue},
  pdfcreator={LaTeX via pandoc}}

\newcommand{\runninghead}{A Preprint }
\renewcommand{\runninghead}{MO Optimization and HPT With Desirability
Functions }
\title{Multi-Objective Optimization and Hyperparameter Tuning With
Desirability Functions}
\usepackage{etoolbox}
\makeatletter
\providecommand{\subtitle}[1]{% add subtitle to \maketitle
  \apptocmd{\@title}{\par {\large #1 \par}}{}{}
}
\makeatother
\subtitle{An Approach Using the Python Package spotdesirability}
\author{\textbf{Thomas
Bartz-Beielstein}~\orcidlink{0000-0002-5938-5158}\\\\Bartz \& Bartz
GmbH\\\\\\\\51643 Gummersbach,
Germany\\\href{mailto:bartzbeielstein@gmail.com}{bartzbeielstein@gmail.com}}
\date{}
\begin{document}
\maketitle
\begin{abstract}
The desirability-function approach is a widely adopted method for
optimizing multiple-response processes. Kuhn (2016) implemented the
packages \texttt{desirability} and \texttt{desirability2} in the
statistical programming language \texttt{R}, but no comparable package
exists for \texttt{Python}. The goal of this article is to provide an
introduction to the desirability function approach using the Python
package \texttt{spotdesirability}, which is available as part of the
\texttt{sequential\ parameter\ optimization} framework. After a brief
introduction to the desirability function approach, three examples are
given that demonstrate how to use the desirability functions for (i)
classical optimization, (ii) surrogate-model based optimization, and
(iii) hyperparameter tuning. An extended Morris-Mitchell criterion,
which allows the computation of the search-space coverage, is proposed
and used in a fourth example to handle the exploration-exploitation
trade-off in optimization. Finally, infill-diagnostic plots are
introduced as a tool to visualize the locations of the infill points
with respect to already existing points.
\end{abstract}
{\bfseries \emph Keywords}
\def\sep{\textbullet\ }
desirability function \sep multi-objective optimization \sep surrogate
modeling \sep hyperparameter tuning \sep Morris-Mitchell
criterion \sep maximin criterion \sep 
sequential parameter optimization

\section{Introduction}\label{introduction}

The desirability-function approach is a widely adopted method in
industry for optimizing multiple-response processes (National Institute
of Standards and Technology 2021). It operates on the principle that the
overall ``quality'' of a product or process with multiple quality
characteristics is deemed unacceptable if any characteristic falls
outside the ``desired'' limits. This approach identifies operating
conditions that yield the most ``desirable'' response values,
effectively balancing multiple objectives to achieve optimal outcomes.
Often, different scales are used for various objectives. When combining
these objectives into a single new one, the challenge arises of how to
compare the scales with each other. The fundamental idea of the
desirability index is to transform the deviations of the objective value
from its target value into comparable desirabilities, i.e., onto a
common scale. For this, a target value as well as a lower and/or upper
specification limit must be known for each objective involved. A result
outside the specification limits is assigned a desirability of 0, while
a result at the target value is assigned a desirability of 1. Linear or
nonlinear transformation, such as a power transformation, can be chosen
as the transformation between the specification limits. The desirability
index according to Derringer and Suich (1980) is then the geometric mean
of the desirabilities of the various objectives (Weihs and Jessenberger
1999).

The \texttt{desirability} package (Kuhn 2016), which is written in the
statistical programming language \texttt{R}, contains \texttt{S3}
classes for multivariate optimization using the desirability function
approach of Harington (1965) with functional forms described by
Derringer and Suich (1980). It is available on CRAN, see
\url{https://cran.r-project.org/package=desirability}. A newer version,
the \texttt{desirability2} package, improves on the original
desirability package by enabling in-line computations that can be used
with dplyr pipelines (Kuhn 2025). It is also available on CRAN, see
\url{https://cran.r-project.org/web/packages/desirability2/index.html}.

Hyperparameter Tuning (or Hyperparameter Optimization) is crucial for
configuring machine learning algorithms, as hyperparameters
significantly impact performance (Bartz et al. 2022; Bischl et al. 2023)
To avoid manual, time-consuming, and irreproducible trial-and-error
processes, these tuning methods can be used. They include simple
techniques like grid and random search, as well as advanced approaches
such as evolution strategies, surrogate optimization, Hyperband, and
racing. The tuning process has to consider several objectives, such as
maximizing the model's performance while minimizing the training time or
model complexity. The desirability function approach is a suitable
method for multi-objective optimization, as it allows for the
simultaneous optimization of multiple objectives by combining them into
a single desirability score.

This paper is structured as follows: After presenting the desirability
function approach in Section~\ref{sec-desirability}, we introduce the
\texttt{Python} package \texttt{spotdesirability}, which is a
\texttt{Python} implementation of the \texttt{R} package
\texttt{desirability}. It is available on PyPI, see
\url{https://pypi.org/project/spotdesirability/} and GitHub, see
\url{https://github.com/sequential-parameter-optimization/spotdesirability}.
The introduction is based on several ``hands-on'' examples.
Section~\ref{sec-related-work} provides an overview of related work in
the field of multi-objective optimization and hyperparameter tuning.

Section~\ref{sec-example-chemical-reaction} presents an example of a
chemical reaction with two objectives: conversion and activity. The
example is based on a response surface experiment described by Myers,
Montgomery, and Anderson-Cook (2016) and also used by Kuhn (2016). It
allows a direct comparison of the results obtained with the \texttt{R}
package \texttt{desirability} and the \texttt{Python} package
\texttt{spotdesirability}.

Section~\ref{sec-maximizing-desirability} describes how to maximize the
desirability function using the Nelder-Mead algorithm from the
\texttt{scipy.optimize.minimize} function. This approach is common in
\gls{rsm} (Box and Wilson 1951; Myers, Montgomery, and Anderson-Cook
2016). The optimization process is illustrated using the chemical
reaction example from Section~\ref{sec-example-chemical-reaction}. This
example is based on the example presented in Kuhn (2016), so that,
similar to the comparison in
Section~\ref{sec-example-chemical-reaction}, a comparison of the results
obtained with the \texttt{R} and \texttt{Python} packages is possible.

Section~\ref{sec-surrogate} presents an example of surrogate model-based
optimization (Gramacy 2020; Forrester, Sóbester, and Keane 2008) using
the \texttt{spotdesirability} package. Results from the RSM optimization
can be compared with the results from surrogate model-based
optimization. The surrogate model is based on the \texttt{spotoptim}
package, which is the successor of the \texttt{spotpython} package
(Bartz-Beielstein 2023).

Section~\ref{sec-hyperparameter-tuning} presents an example of
hyperparameter tuning of a neural network implemented in
\texttt{PyTorch} using the \texttt{spotdesirability} package. The goal
of this example is to demonstrate how to use the desirability function
approach for hyperparameter tuning in a deep learning context.
Section~\ref{sec-space-fillingness} presents an example of
space-fillingness as an objective in multi-objective optimization using
the \texttt{spotdesirability} package. We also introduce
infill-diagnostic plots as a tool to visualize the locations of the
infill points with respect to already existing points. The article
concludes with a summary and outlook in Section~\ref{sec-conclusion}.
Supplementary material will be provided in the Sequential Parameter
Optimization Cookbook, which is available at:
\url{https://sequential-parameter-optimization.github.io/spotoptim-cookbook/}.

\section{Desirability}\label{sec-desirability}

\subsection{Basic Desirability
Functions}\label{basic-desirability-functions}

The desirability function approach to simultaneously optimizing multiple
equations was originally proposed by Harington (1965). The approach
translates the functions to a common scale (\([0, 1]\)), combines them
using the geometric mean, and optimizes the overall metric. The
equations can represent model predictions or other equations. Kuhn
(2016) notes that desirability functions are popular in \gls{rsm} (Box
and Wilson 1951; Myers, Montgomery, and Anderson-Cook 2016) to
simultaneously optimize a series of quadratic models. A response surface
experiment may use measurements on a set of outcomes, where instead of
optimizing each outcome separately, settings for the predictor variables
are sought to satisfy all outcomes at once.

Kuhn (2016) explains that originally, Harrington used exponential
functions to quantify desirability. In our \texttt{Python}
implementation, which is based on the \texttt{R} package
\texttt{desirablity} from Kuhn (2016), the simple discontinuous
functions of Derringer and Suich (1980) are adopted. For simultaneous
optimization of equations, individual ``desirability'' functions are
constructed for each function, and Derringer and Suich (1980) proposed
three forms of these functions corresponding to the optimization goal
type. Kuhn (2016) describes the \texttt{R} implementation as follows:

\begin{quote}
Suppose there are \(R\) equations or functions to simultaneously
optimize, denoted \(f_r(\vec{x})\) (\(r = 1 \ldots R\)). For each of the
\(R\) functions, an individual ``desirability'' function is constructed
that is high when \(f_r(\vec{x})\) is at the desirable level (such as a
maximum, minimum, or target) and low when \(f_r(\vec{x})\) is at an
undesirable value. Derringer and Suich (1980) proposed three forms of
these functions, corresponding to the type of optimization goal, namely
maximization, minimization, or target optimization. The associated
desirability functions are denoted \(d_r^{\text{max}}\),
\(d_r^{\text{min}}\), and \(d_r^{\text{target}}\).
\end{quote}

\subsubsection{Maximization}\label{maximization}

For maximization of \(f_r(\vec{x})\) (``larger-is-better''), the
following function is used:

\[
d_r^{\text{max}} =
\begin{cases}
    0 & \text{if } f_r(\vec{x}) < A \\
    \left(\frac{f_r(\vec{x}) - A}{B - A}\right)^s & \text{if } A \leq f_r(\vec{x}) \leq B \\
    1 & \text{if } f_r(\vec{x}) > B,
\end{cases}
\]

where \(A\), \(B\), and \(s\) are chosen by the investigator.

\subsubsection{Minimization}\label{minimization}

For minimization (``smaller-is-better''), the following function is
proposed:

\[
d_r^{\text{min}} =
\begin{cases}
    0 & \text{if } f_r(\vec{x}) > B \\
    \left(\frac{f_r(\vec{x}) - B}{A - B}\right)^s & \text{if } A \leq f_r(\vec{x}) \leq B \\
    1 & \text{if } f_r(\vec{x}) < A
\end{cases}
\]

\subsubsection{Target Optimization}\label{target-optimization}

In ``target-is-best'' situations, the following function is used:

\[
d_r^{\text{target}} =
\begin{cases}
    \left(\frac{f_r(\vec{x}) - A}{t_0 - A}\right)^{s_1} & \text{if } A \leq f_r(\vec{x}) \leq t_0 \\
    \left(\frac{f_r(\vec{x}) - B}{t_0 - B}\right)^{s_2} & \text{if } t_0 \leq f_r(\vec{x}) \leq B \\
    0 & \text{otherwise.}
\end{cases}
\]

\begin{figure}

\centering{

\pandocbounded{\includegraphics[keepaspectratio]{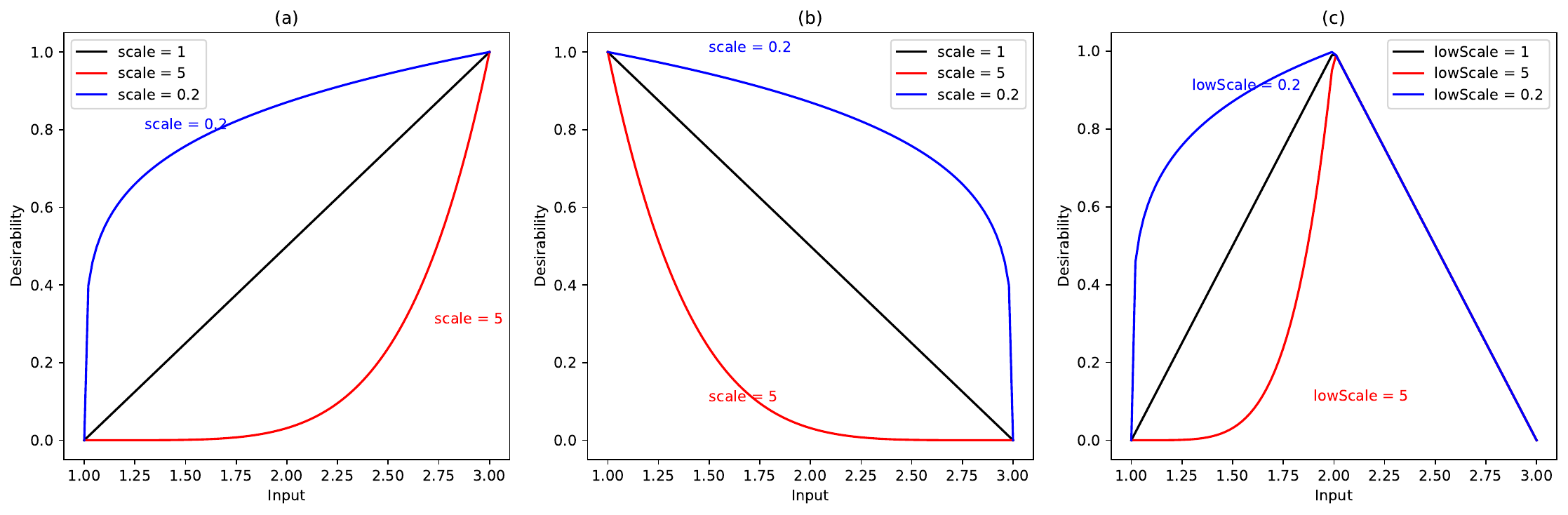}}

}

\caption{\label{fig-kuhn16a-1}Examples of the three primary desirability
functions. Panel (a) shows an example of a larger--is--better function,
panel (b) shows a smaller--is--better desirability function and panel
(c) shows a function where the optimal value corresponds to a target
value. Not that increasing the scale parameter makes it more difficult
to achieve higher desirability, while values smaller than 1 make it
easier to achieve good results.}

\end{figure}

Kuhn (2016) explains that these functions, which are shown in
Figure~\ref{fig-kuhn16a-1}, share the same scale and are discontinuous
at specific points \(A\), \(B\), and \(t_0\). The values of \(s\),
\(s_1\), or \(s_2\) can be chosen so that the desirability criterion is
easier or more difficult to satisfy. For example:

\begin{itemize}
\tightlist
\item
  If \(s\) is chosen to be less than 1 in \(d_r^{\text{min}}\),
  \(d_r^{\text{min}}\) is near 1 even if the model \(f_r(\vec{x})\) is
  not low.
\item
  As values of \(s\) move closer to 0, the desirability reflected by
  \(d_r^{\text{min}}\) becomes higher.
\item
  Values of \(s\) greater than 1 will make \(d_r^{\text{min}}\) harder
  to satisfy in terms of desirability.
\end{itemize}

These scaling factors are useful when one equation holds more importance
than others. Any function can reflect model desirability, i.e., Del
Castillo, Montgomery, and McCarville (1996) developed alternative
functions suitable for gradient-based optimizations.

For each of these three desirability functions (and the others discussed
in Section~\ref{sec-nonstandard-desirabilities}), there are
\texttt{print\_class\_attributes}, \texttt{plot}, and \texttt{predict}
methods similar to the \texttt{R} implementation (Kuhn 2016). The
\texttt{print\_attributes} method prints the class attributes, the
\texttt{plot} method plots the desirability function, and the
\texttt{predict} method predicts the desirability for a given input.

\subsection{Overall Desirability}\label{overall-desirability}

Given the \(R\) desirability functions \(d_1 \ldots d_r\) are on the
{[}0,1{]} scale, they can be combined to achieve an overall desirability
function, \(D\). One method of doing this is by the geometric mean:

\[
D = \left(\prod_{r=1}^R d_r\right)^{1/R}.
\]

The geometric mean has the property that if any one model is undesirable
(\(d_r = 0\)), the overall desirability is also unacceptable
(\(D = 0\)). Once \(D\) has been defined and the prediction equations
for each of the \(R\) equations have been computed, it can be used to
optimize or rank the predictors.

\subsection{Non-Standard Features}\label{sec-nonstandard}

The \texttt{R} package \texttt{desirability} (Kuhn 2016) offers a few
non-standard features. These non-standard features are also included in
the \texttt{Python} implementation and will be discussed in the
following. First, we will consider the non-informative desirability and
missing values, followed by zero-desirability tolerances, and finally
non-standard desirability functions.

\subsubsection{Non-Informative Desirability and Missing
Values}\label{sec-missing}

According to Kuhn (2016), if inputs to desirability functions are
uncomputable, the package estimates a non-informative value by computing
desirabilities over the possible range and taking the mean.\\
If an input to a desirability function is \texttt{NA}, by default, it is
replaced by this non-informative value. Setting \texttt{object\$missing}
to \texttt{NA} (in \texttt{R}) changes the calculation to return an
\texttt{NA} for the result, where \texttt{object} is the result of a
call to one of the desirability functions.A similar procedure is
implemented in the \texttt{Python} package. The non-informative value is
plotted as a broken line in default \texttt{plot} methods.

\subsubsection{Zero-Desirability
Tolerances}\label{zero-desirability-tolerances}

Kuhn (2016) highlights that in high-dimensional outcomes, finding
feasible solutions where every desirability value is acceptable can be
challenging. Each desirability R function has a \texttt{tol} argument,
which can be set between {[}0, 1{]} (default is \texttt{NULL}). If not
\texttt{NULL}, zero desirability values are replaced by \texttt{tol}.

\subsubsection{Non-Standard Desirability
Functions}\label{sec-nonstandard-desirabilities}

Kuhn (2016) mentions scenarios where the three discussed desirability
functions are inadequate for user requirements. In this case, the
\texttt{dArb} function (\texttt{Arb} stands for ``Arbitary'') can be
used to create a custom desirability function. \texttt{dArb} accepts
numeric vector inputs with matching desirabilities to approximate other
functional forms.

\begin{example}[Logistic Desirability
Function]\protect\hypertarget{exm-logistic}{}\label{exm-logistic}

A logistic function can be used as a desirability function. The logistic
function is defined as \[
d(\vec{x}) = \frac{1}{1+\exp(-\vec{x})}.
\] For inputs outside the range \(\pm5\), desirability values remain
near zero and one. The desirability function is defined using 20
computation points on this range.

\phantomsection\label{kuhn16a-logistic}
\begin{Shaded}
\begin{Highlighting}[]
\KeywordTok{def}\NormalTok{ foo(u):}
    \ControlFlowTok{return} \DecValTok{1} \OperatorTok{/}\NormalTok{ (}\DecValTok{1} \OperatorTok{+}\NormalTok{ np.exp(}\OperatorTok{{-}}\NormalTok{u))}
\NormalTok{x\_input }\OperatorTok{=}\NormalTok{ np.linspace(}\OperatorTok{{-}}\DecValTok{5}\NormalTok{, }\DecValTok{5}\NormalTok{, }\DecValTok{20}\NormalTok{)}
\NormalTok{logistic\_d }\OperatorTok{=}\NormalTok{ DArb(x\_input, foo(x\_input))}
\NormalTok{logistic\_d.print\_class\_attributes()}
\end{Highlighting}
\end{Shaded}

\begin{verbatim}

Class: DArb
x: [-5.         -4.47368421 -3.94736842 -3.42105263 -2.89473684 -2.36842105
 -1.84210526 -1.31578947 -0.78947368 -0.26315789  0.26315789  0.78947368
  1.31578947  1.84210526  2.36842105  2.89473684  3.42105263  3.94736842
  4.47368421  5.        ]
d: [0.00669285 0.01127661 0.0189398  0.03164396 0.05241435 0.08561266
 0.1368025  0.21151967 0.31228169 0.43458759 0.56541241 0.68771831
 0.78848033 0.8631975  0.91438734 0.94758565 0.96835604 0.9810602
 0.98872339 0.99330715]
tol: None
missing: 0.5
\end{verbatim}

Inputs in-between these grid points are linearly interpolated. Using
this method, extreme values are applied outside the input range.
Figure~\ref{fig-kuhn16a-7} displays a \texttt{plot} of the
\texttt{logisticD} object.

\begin{Shaded}
\begin{Highlighting}[]
\NormalTok{logistic\_d.plot(figsize}\OperatorTok{=}\NormalTok{(}\DecValTok{4}\NormalTok{, }\DecValTok{3}\NormalTok{))}
\end{Highlighting}
\end{Shaded}

\begin{figure}[H]

\centering{

\pandocbounded{\includegraphics[keepaspectratio]{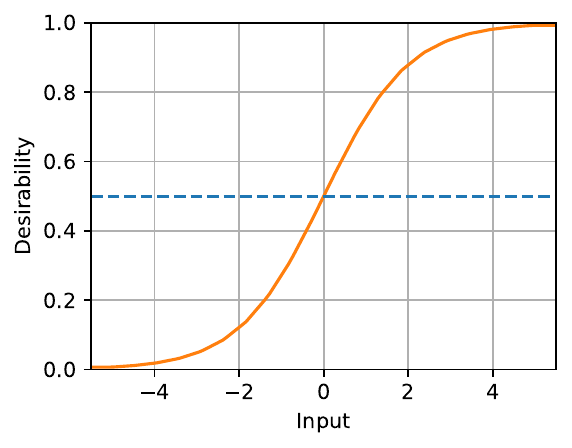}}

}

\caption{\label{fig-kuhn16a-7}An example of a desirability function
created using the \texttt{DArb} function. The desirability function is a
logistic curve that is defined by 20 points on the range {[}-5, 5{]}.}

\end{figure}%

\end{example}

Kuhn also proposes a desirability function for implementing box
constraints on an equation.

\begin{example}[Box-Constraint Desirability
Function]\protect\hypertarget{exm-box}{}\label{exm-box}

The \texttt{DBox} class can be used to create a desirability function
for box constraints. For example, assigning zero desirability to values
beyond \(\pm 1.682\) in the design region, instead of penalizing.
Figure~\ref{fig-kuhn16a-8} demonstrates an example function.

\begin{Shaded}
\begin{Highlighting}[]
\NormalTok{box\_desirability }\OperatorTok{=}\NormalTok{ DBox(low}\OperatorTok{={-}}\FloatTok{1.682}\NormalTok{, high}\OperatorTok{=}\FloatTok{1.682}\NormalTok{)}
\NormalTok{box\_desirability.plot(non\_inform}\OperatorTok{=}\VariableTok{False}\NormalTok{, figsize}\OperatorTok{=}\NormalTok{(}\DecValTok{4}\NormalTok{,}\DecValTok{3}\NormalTok{))}
\end{Highlighting}
\end{Shaded}

\begin{figure}[H]

\centering{

\pandocbounded{\includegraphics[keepaspectratio]{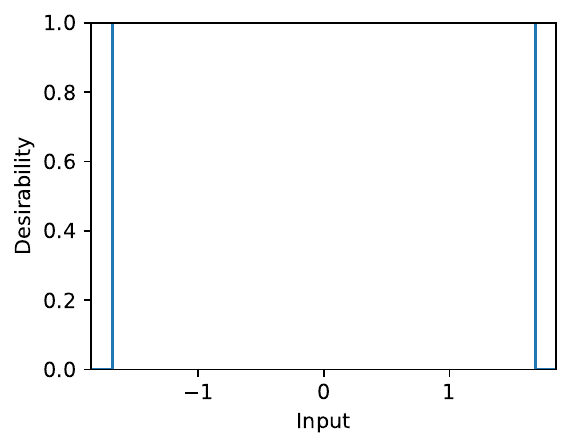}}

}

\caption{\label{fig-kuhn16a-8}An example of a box-like desirability
function that assigns zero desirability to values outside of the range
{[}-1.682, 1.682{]}.}

\end{figure}%

\end{example}

\begin{example}[Desirability Function for Categorical
Inputs]\protect\hypertarget{exm-categorical}{}\label{exm-categorical}

Kuhn concludes by mentioning another non-standard application involving
categorical inputs. Desirabilities are assigned to each value.
Figure~\ref{fig-kuhn16a-9} visualizes a plot of desirability profiles
for this setup. Using \texttt{spotdesirability}, the
\texttt{DCategorical} class can be used to create a desirability
function for categorical inputs as follows: first, the desirability
values for each category are defined, and then the \texttt{DCategorical}
class is used to create a desirability function for the categorical
inputs.

\phantomsection\label{kuhn16a-categorical}
\begin{Shaded}
\begin{Highlighting}[]
\NormalTok{values }\OperatorTok{=}\NormalTok{ \{}\StringTok{"value1"}\NormalTok{: }\FloatTok{0.1}\NormalTok{, }\StringTok{"value2"}\NormalTok{: }\FloatTok{0.9}\NormalTok{, }\StringTok{"value3"}\NormalTok{: }\FloatTok{0.2}\NormalTok{\}}
\NormalTok{grouped\_desirabilities }\OperatorTok{=}\NormalTok{ DCategorical(values)}
\BuiltInTok{print}\NormalTok{(}\StringTok{"Desirability values for categories:"}\NormalTok{)}
\ControlFlowTok{for}\NormalTok{ category, desirability }\KeywordTok{in}\NormalTok{ grouped\_desirabilities.values.items():}
    \BuiltInTok{print}\NormalTok{(}\SpecialStringTok{f"}\SpecialCharTok{\{}\NormalTok{category}\SpecialCharTok{\}}\SpecialStringTok{: }\SpecialCharTok{\{}\NormalTok{desirability}\SpecialCharTok{\}}\SpecialStringTok{"}\NormalTok{)}

\CommentTok{\# Example usage: Predict desirability for a specific category}
\NormalTok{category }\OperatorTok{=} \StringTok{"value2"}
\NormalTok{predicted\_desirability }\OperatorTok{=}\NormalTok{ grouped\_desirabilities.predict([category])}
\BuiltInTok{print}\NormalTok{(}\SpecialStringTok{f"}\CharTok{\textbackslash{}n}\SpecialStringTok{Predicted desirability for \textquotesingle{}}\SpecialCharTok{\{}\NormalTok{category}\SpecialCharTok{\}}\SpecialStringTok{\textquotesingle{}: }\SpecialCharTok{\{}\NormalTok{predicted\_desirability[}\DecValTok{0}\NormalTok{]}\SpecialCharTok{\}}\SpecialStringTok{"}\NormalTok{)}
\end{Highlighting}
\end{Shaded}

\begin{verbatim}
Desirability values for categories:
value1: 0.1
value2: 0.9
value3: 0.2

Predicted desirability for 'value2': 0.9
\end{verbatim}

\begin{Shaded}
\begin{Highlighting}[]
\NormalTok{grouped\_desirabilities.plot(non\_inform}\OperatorTok{=}\VariableTok{False}\NormalTok{, figsize}\OperatorTok{=}\NormalTok{(}\DecValTok{4}\NormalTok{,}\DecValTok{3}\NormalTok{))}
\end{Highlighting}
\end{Shaded}

\begin{figure}[H]

\centering{

\pandocbounded{\includegraphics[keepaspectratio]{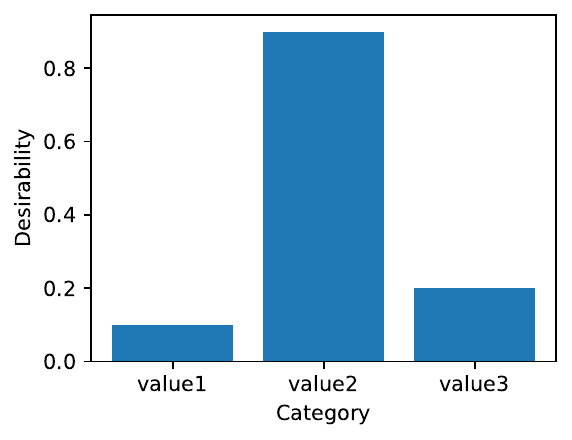}}

}

\caption{\label{fig-kuhn16a-9}Desirability function for categorical
values. The desirability values are assigned to three categories:
`value1', `value2', and `value3'.}

\end{figure}%

\end{example}

\section{Related Work}\label{sec-related-work}

Multiobjective approaches are established optimization tools (Emmerich
and Deutz 2018). The weighted-sum approach is a simple and widely used
method for multi-objective optimization, but probably only, because its
disadvantages are unknown. Compared to the weighted-sum approach, the
desirability-function approach is a better choice for multi-objective
optimization. The desirability function approach also allows for more
flexibility in defining the objectives and their trade-offs. Nino et al.
(2015) discuss the use of Experimental Designs and \gls{rsm} to optimize
conflicting responses in the development of a 3D printer prototype.
Specifically, they focus on an interlocking device designed to recycle
polyethylene terephthalate water bottles. The optimization involves two
conflicting goals: maximizing load capacity and minimizing mass. A
Box-Behnken Design (BBD) was used for the experimental setup, and
desirability functions were applied to identify the best trade-offs.
Karl et al. (2023) describe multi-objective optimization in maschine
learning. Coello et al. (2021) give an overview of Multi-Objective
Evolutionary Algorithms. Bartz-Beielstein (2025) provides an
introduction to the desirability function approach to multi-objective
optimization (direct and surrogate model-based), and multi-objective
hyperparameter tuning.

\section{An Example With Two Objectives: Chemical
Reaction}\label{sec-example-chemical-reaction}

Similar to the presentation in Kuhn (2016), we will use the example of a
chemical reaction to illustrate the desirability function approach. The
example is based on a response surface experiment described by Myers,
Montgomery, and Anderson-Cook (2016). The goal is to maximize the
percent conversion of a chemical reaction while keeping the thermal
activity within a specified range.

\begin{figure}

\centering{

\pandocbounded{\includegraphics[keepaspectratio]{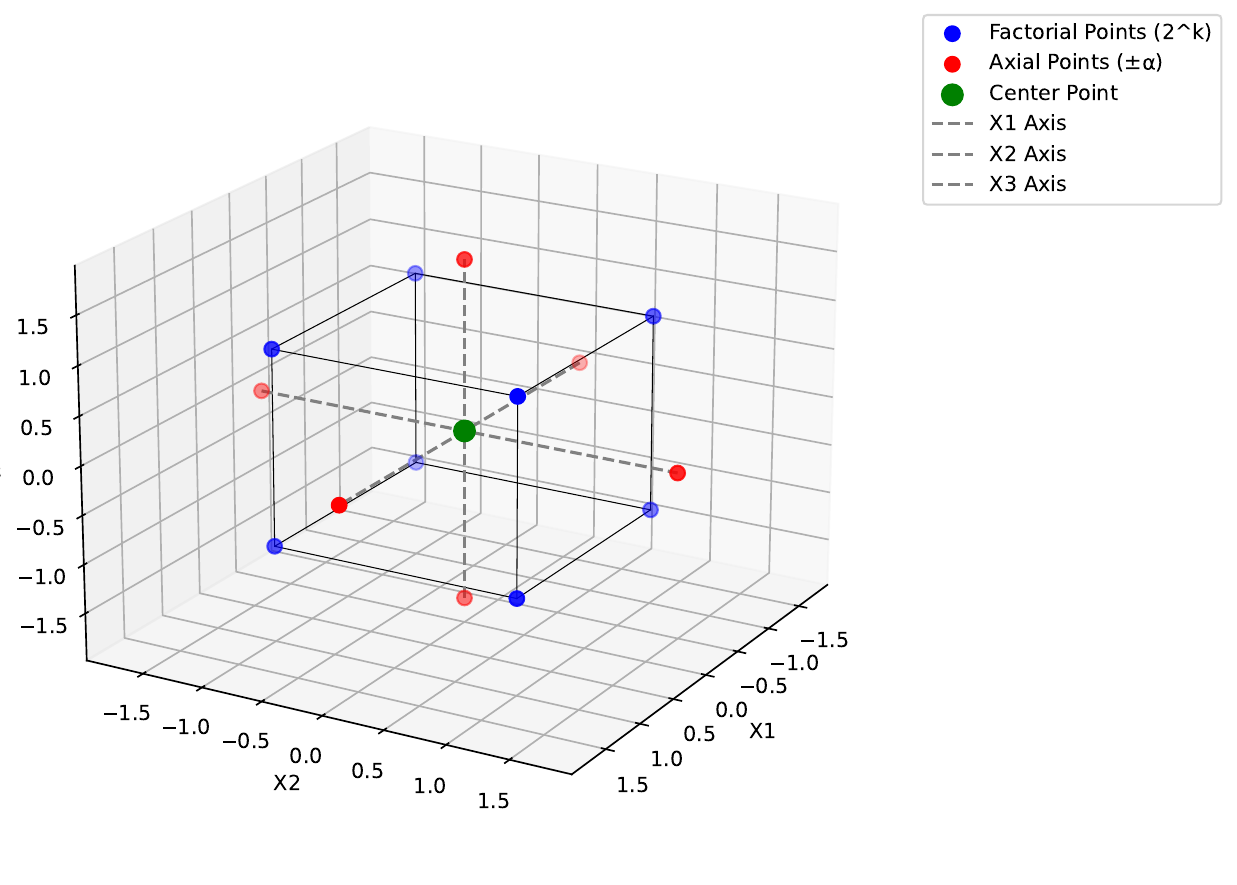}}

}

\caption{\label{fig-ccd}Central composite design (CCD) for the chemical
reaction example}

\end{figure}

The \gls{ccd} is the most popular class of designs used for fitting
second-order response surface models (Montgomery 2001). Since the
location of the optimum is unknown before the \gls{rsm} starts, Box and
Hunter (1957) suggested that the design should be rotatable (it provides
equal precision of estimation in all directions or stated differently,
the variance of the predicted response is the same at all points that
are the same distance from the center of the design space). A \gls{ccd}
is made rotable by using an axis distance value of
\(\alpha = (n_F)^{1/4}\), where \(n_F\) is the number of points (here
\(2^3 = 8\)) (Montgomery 2001). Figure~\ref{fig-ccd} shows the design
space for the chemical reaction example. The design space is defined by
three variables: reaction time, reaction temperature, and percent
catalyst. This rotable \gls{ccd} consists of a full factorial design
with three factors, each at two levels, plus a center point and six
(\(2\times k\)) axial points. The axial points are located at a distance
of \(\pm \alpha\) from the center point in each direction.

Montgomery (2001) notes that it is \emph{not} important to have exact
rotability. From a prediction variance point of view, the best choice is
to set \(\alpha = \sqrt{k}\), which results in a so-called spherical
\gls{ccd}.

\subsection{The Two Objective Functions: Conversion and
Activity}\label{the-two-objective-functions-conversion-and-activity}

Myers, Montgomery, and Anderson-Cook (2016) present two equations for
the fitted quadratic response surface models, the first one for the
conversion and the second one for the activity.

\begin{align*}
f_{\text{con}}(x) =
&
 81.09
+
1.0284 \cdot x_1
+
4.043 \cdot x_2
+
6.2037 \cdot x_3
+
1.8366 \cdot x_1^2
+
2.9382 \cdot x_2^2 \\
&
+
5.1915 \cdot x_3^2
+
2.2150 \cdot x_1 \cdot x_2
+
11.375 \cdot x_1 \cdot x_3
+
3.875 \cdot x_2 \cdot x_3
\end{align*} and \begin{align*}
f_{\text{act}}(x) = 
 & 
 59.85
+ 3.583 \cdot x_1
+ 0.2546 \cdot x_2
+ 2.2298 \cdot x_3
+ 0.83479 \cdot x_1^2
+ 0.07484 \cdot x_2^2
\\
&
+ 0.05716 \cdot x_3^2
+ 0.3875 \cdot x_1 \cdot x_2
+ 0.375 \cdot x_1 \cdot x_3
+ 0.3125 \cdot x_2 \cdot x_3. 
\end{align*}

They are implemented as Python functions that take a vector of three
parameters (\(x_1\): reaction time, \(x_2\): reaction temperature,
\(x_3\): percent catalyst) and return the predicted values for the
percent conversion and thermal activity. They are available in the
\texttt{spotdesirability} package as \texttt{conversion\_pred(x)} and
\texttt{activity\_pred(x)}.

The goal of the analysis in Myers, Montgomery, and Anderson-Cook (2016)
was to

\begin{itemize}
\tightlist
\item
  maximize conversion while
\item
  keeping the thermal activity between 55 and 60 units. An activity
  target of 57.5 was used in the analysis.
\end{itemize}

Plots of the response surface models are shown in
Figure~\ref{fig-kuhn16a-2} and Figure~\ref{fig-kuhn16a-3}, where
reaction time and percent catalyst are plotted while the reaction
temperature was varied at four different levels. Both quadratic models,
as pointed out by Kuhn (2016), are saddle surfaces, and the stationary
points are outside of the experimental region. To determine predictor
settings for these models, a constrained optimization can be used to
stay inside the experimental region. Kuhn notes:

\begin{quote}
In practice, we would just use the \texttt{predict} method for the
linear model objects to get the prediction equation. Our results are
slightly different from those given by Myers and Montgomery because they
used prediction equations with full floating-point precision.
\end{quote}

\subsection{Contour-Plot Generation}\label{sec-contour-plot-generation}

We will generate contour plots for the percent conversion and thermal
activity models. The contour-plot generation comprehends the following
steps:

\begin{itemize}
\tightlist
\item
  generating a grid of points in the design space,
\item
  evaluating the response surface models at these points,
\item
  and plotting the contour plots for the response surface models.
\end{itemize}

\begin{figure}

\centering{

\pandocbounded{\includegraphics[keepaspectratio]{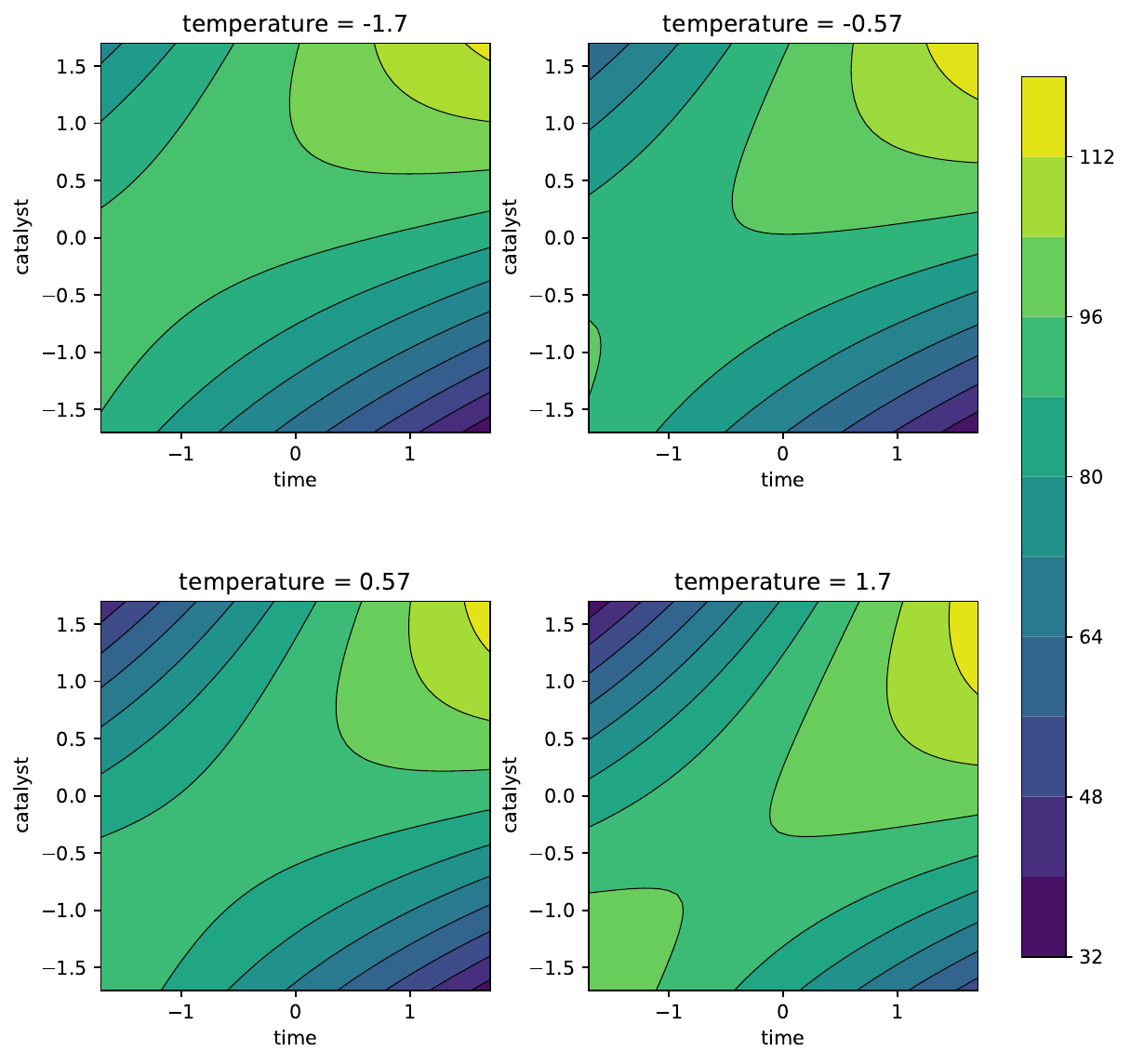}}

}

\caption{\label{fig-kuhn16a-2}The response surface for the percent
conversion model. To plot the model contours, the temperature variable
was fixed at four diverse levels.}

\end{figure}

\begin{figure}

\centering{

\pandocbounded{\includegraphics[keepaspectratio]{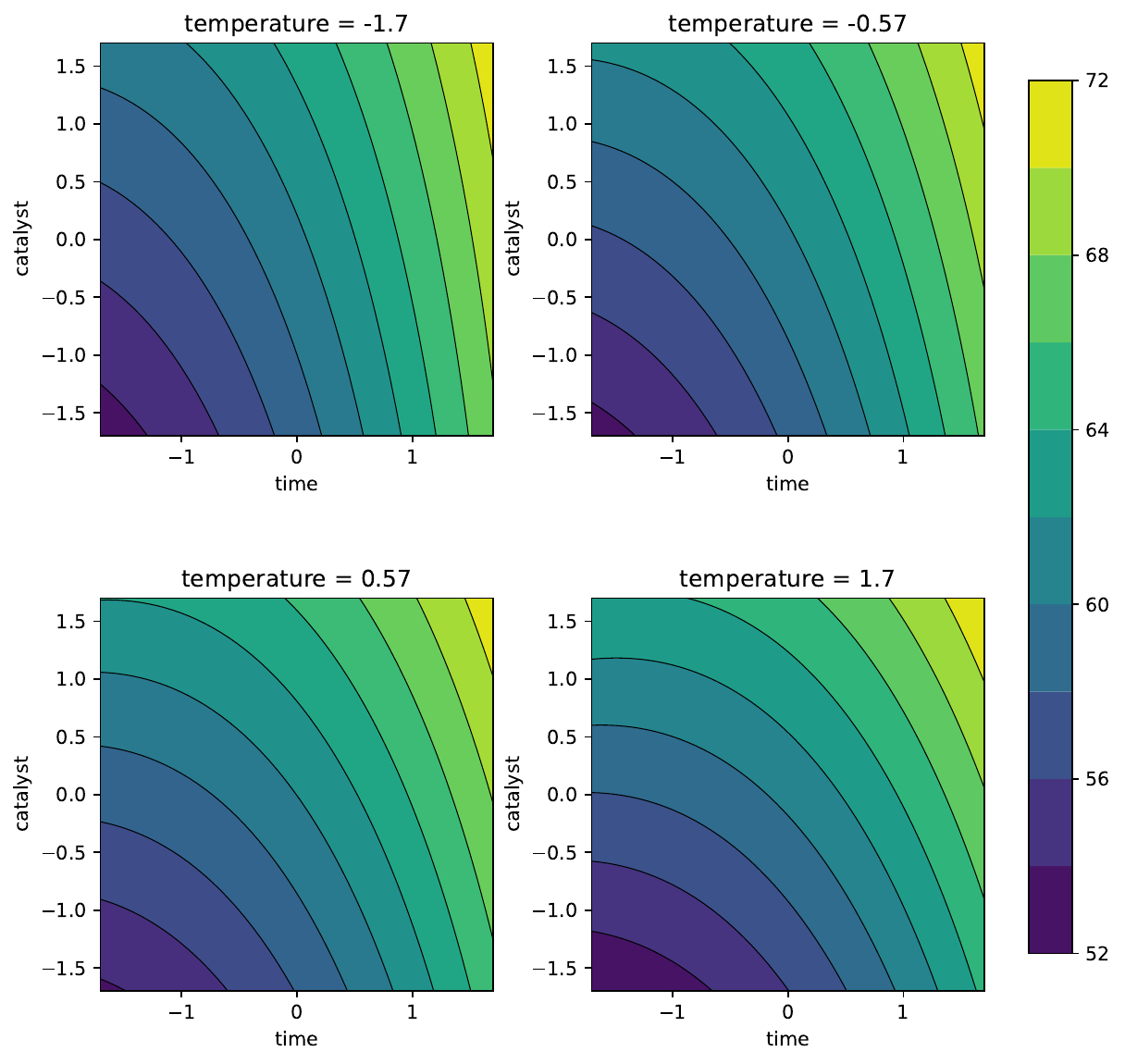}}

}

\caption{\label{fig-kuhn16a-3}The response surface for the thermal
activity model. To plot the model contours, the temperature variable was
fixed at four diverse levels.}

\end{figure}

We will use the function \texttt{mo\_generate\_plot\_grid} to generate
the grid and the function \texttt{mo\_contourf\_plots} for creating the
contour plots for the response surface models. Both functions are
available in the \texttt{spotoptim} package. First we define the
variables, their ranges, the resolutions for the grid, and the objective
functions using \texttt{Python}'s \texttt{dictionary} data structure.
The \texttt{variables} dictionary contains the variable names as keys
and their ranges as values. The \texttt{resolutions} dictionary contains
the variable names as keys and their resolutions as values. The
\texttt{functions} dictionary contains the function names as keys and
the corresponding functions as values. Next we can generate the Pandas
DataFrame \texttt{plot\_grid} using \texttt{spotoptim}'s
\texttt{mo\_generate\_plot\_grid} function. It has the columns
\texttt{time}, \texttt{temperature}, \texttt{catalyst},
\texttt{conversionPred}, and \texttt{activityPred}. The
\texttt{plot\_grid} DataFrame is used to generate the contour plots for
the response surface models using \texttt{spotoptim}'s
\texttt{contourf\_plot} function. Figure~\ref{fig-kuhn16a-2} shows the
response surface for the percent conversion model. To plot the model
contours, the temperature variable was fixed at four diverse levels. The
largest effects in the fitted model are due to the time \(\times\)
catalyst interaction and the linear and quadratic effects of catalyst.
Figure~\ref{fig-kuhn16a-3} shows the response surface for the thermal
activity model. To plot the model contours, the temperature variable was
fixed at four diverse levels. The main effects of time and catalyst have
the largest effect on the fitted model.

\subsection{Defining the Desirability
Functions}\label{sec-defining-desirability}

\begin{figure}[H]

\centering{

\pandocbounded{\includegraphics[keepaspectratio]{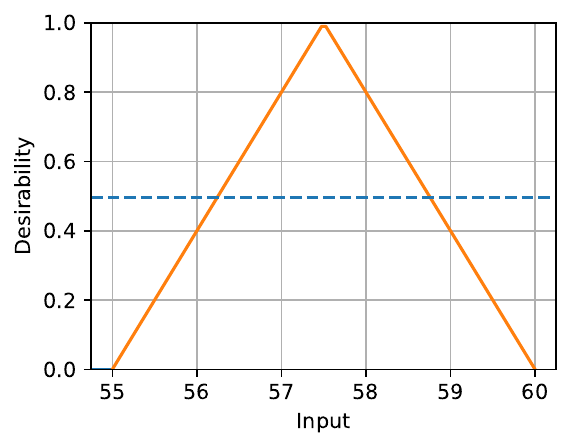}}

}

\caption{\label{fig-kuhn16a-desirability-plot-activity}The desirability
function for the activity outcome.}

\end{figure}

Following the steps described in Kuhn (2016), translating the
experimental goals to desirability functions, a larger-is-better
function (\(d_r^{\text{max}}\)) is used for percent conversion with
values \(A = 80\) and \(B = 97\). A target-oriented desirability
function (\(d_r^{\text{target}}\)) was used for thermal activity with
\(t_0 = 57.5\), \(A = 55\), and \(B = 60\). Kuhn emphasizes that to
construct the overall desirability functions, objects must be created
for the individual functions. In the following, we will use classes of
the \texttt{Python} package \texttt{spotdesirability} to create the
desirability objects. The \texttt{spotdesirability} package is part of
the \texttt{sequential\ parameter\ optimization} framework
(Bartz-Beielstein 2023). It is available on
\url{https://github.com/sequential-parameter-optimization/spotdesirability}
and on \url{https://pypi.org/project/spotdesirability} and can be
installed via \texttt{pip\ install\ spotdesirability}. The desirability
objects can be created as shown below.
Figure~\ref{fig-kuhn16a-desirability-plot-conversion} and
Figure~\ref{fig-kuhn16a-desirability-plot-activity} show the
desirability functions for the conversion and activity outcomes,
respectively.

\phantomsection\label{kuhn16a-desirability-obj}
\begin{Shaded}
\begin{Highlighting}[]
\NormalTok{conversionD }\OperatorTok{=}\NormalTok{ DMax(}\DecValTok{80}\NormalTok{, }\DecValTok{97}\NormalTok{)}
\NormalTok{activityD }\OperatorTok{=}\NormalTok{ DTarget(}\DecValTok{55}\NormalTok{, }\FloatTok{57.5}\NormalTok{, }\DecValTok{60}\NormalTok{)}
\end{Highlighting}
\end{Shaded}

\begin{Shaded}
\begin{Highlighting}[]
\NormalTok{conversionD.plot(figsize}\OperatorTok{=}\NormalTok{(}\DecValTok{4}\NormalTok{,}\DecValTok{3}\NormalTok{))}
\end{Highlighting}
\end{Shaded}

\begin{figure}[H]

\centering{

\pandocbounded{\includegraphics[keepaspectratio]{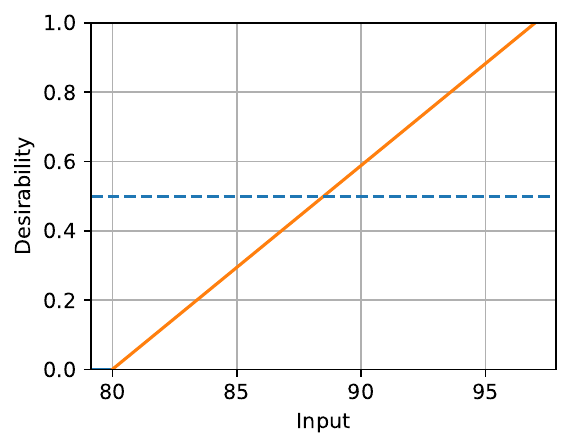}}

}

\caption{\label{fig-kuhn16a-desirability-plot-conversion}The
desirability function for the conversion outcome.}

\end{figure}%

\begin{Shaded}
\begin{Highlighting}[]
\NormalTok{activityD.plot(figsize}\OperatorTok{=}\NormalTok{(}\DecValTok{4}\NormalTok{,}\DecValTok{3}\NormalTok{))}
\end{Highlighting}
\end{Shaded}

Although the original analysis in Myers, Montgomery, and Anderson-Cook
(2016) used numerous combinations of scaling parameters, following the
presentation in Kuhn (2016), we will only show analyses with the default
scaling factor values. The use of the scaling parameter is illustrated
later in Section~\ref{sec-space-fillingness}.

\begin{example}[Computing Desirability at the Center
Point]\protect\hypertarget{exm-desirability}{}\label{exm-desirability}

Using these desirability objects \texttt{conversionD}and
\texttt{activityD}, the following code segment shows how to predict the
desirability for the center point of the experimental design. The center
point is defined as {[}0, 0, 0{]}. First, we compute the predicted
outcomes for the center point using the response surface models. Then
these predicted outcomes are used to compute the desirability for each
outcome.

\phantomsection\label{predict-outcomes}
\begin{Shaded}
\begin{Highlighting}[]
\NormalTok{pred\_outcomes }\OperatorTok{=}\NormalTok{ [}
\NormalTok{    conversion\_pred([}\DecValTok{0}\NormalTok{, }\DecValTok{0}\NormalTok{, }\DecValTok{0}\NormalTok{]),}
\NormalTok{    activity\_pred([}\DecValTok{0}\NormalTok{, }\DecValTok{0}\NormalTok{, }\DecValTok{0}\NormalTok{])}
\NormalTok{]}
\CommentTok{\# Predict desirability for each outcome}
\NormalTok{conversion\_desirability }\OperatorTok{=}\NormalTok{ conversionD.predict(pred\_outcomes[}\DecValTok{0}\NormalTok{])}
\NormalTok{activity\_desirability }\OperatorTok{=}\NormalTok{ activityD.predict(pred\_outcomes[}\DecValTok{1}\NormalTok{])}
\end{Highlighting}
\end{Shaded}

\phantomsection\label{print-desirability-outcomes}
\begin{verbatim}
Predicted Outcomes: [np.float64(81.09), np.float64(59.85)]
Conversion Desirability: [0.06411765]
Activity Desirability: [0.06]
\end{verbatim}

These two desirabilities are close to zero, because the center point of
the experimental design is not a good choice for this experiment.
Similar to the implementation in Kuhn (2016), to get the overall score
for these settings of the experimental factors, the \texttt{dOverall}
function is used to combine the objects and \texttt{predict} is used to
get the final score. The \texttt{print\_class\_attributes} method prints
the class attributes of the \texttt{DOverall} object.

\begin{Shaded}
\begin{Highlighting}[]
\NormalTok{overallD }\OperatorTok{=}\NormalTok{ DOverall(conversionD, activityD)}
\NormalTok{overallD.print\_class\_attributes()}
\end{Highlighting}
\end{Shaded}

\begin{verbatim}

Class: DOverall
d_objs: [

  Class: DMax
  low: 80
  high: 97
  scale: 1
  tol: None
  missing: 0.5

  Class: DTarget
  low: 55
  target: 57.5
  high: 60
  low_scale: 1
  high_scale: 1
  tol: None
  missing: 0.4949494949494951
]
\end{verbatim}

Note: The attribute \texttt{missing} denotes the value that is used for
missing values, see Section~\ref{sec-missing} for more details. Finally,
we can print the overall desirability for the center point of the
experimental design.

\phantomsection\label{kuhn16a-desirability-predict-tree}
\begin{Shaded}
\begin{Highlighting}[]
\NormalTok{overall\_desirability }\OperatorTok{=}\NormalTok{ overallD.predict(pred\_outcomes, }\BuiltInTok{all}\OperatorTok{=}\VariableTok{True}\NormalTok{)}
\BuiltInTok{print}\NormalTok{(}\StringTok{"Conversion Desirability:"}\NormalTok{, overall\_desirability[}\DecValTok{0}\NormalTok{][}\DecValTok{0}\NormalTok{])}
\BuiltInTok{print}\NormalTok{(}\StringTok{"Activity Desirability:"}\NormalTok{, overall\_desirability[}\DecValTok{0}\NormalTok{][}\DecValTok{1}\NormalTok{])}
\BuiltInTok{print}\NormalTok{(}\StringTok{"Overall Desirability:"}\NormalTok{, overall\_desirability[}\DecValTok{1}\NormalTok{])}
\end{Highlighting}
\end{Shaded}

\begin{verbatim}
Conversion Desirability: [0.06411765]
Activity Desirability: [0.06]
Overall Desirability: [0.06202466]
\end{verbatim}

Similar to the single desirabilities, the overall desirability is also
close to zero, because the center point of the experimental design is
not a good choice for this experiment.

\end{example}

\subsection{The Desirability DataFrame and Desirability Contour
Plots}\label{the-desirability-dataframe-and-desirability-contour-plots}

A DataFrame \texttt{d\_values\_df} is created to store the individual
desirability values for each outcome, and the overall desirability value
is added as a new column. First, we predict desirability values and
extract the individual and overall desirability values. Note: The
\texttt{all=True} argument indicates that both individual and overall
desirability values should be returned. We add the individual and
overall desirability values to the \texttt{plot\_grid} DataFrame, that
was created earlier in Section~\ref{sec-contour-plot-generation}.

We will use \texttt{spotoptim}'s \texttt{contourf\_plot} function to
create the contour plots for the individual desirability surfaces and
the overall desirability surface. The \texttt{plot\_grid} DataFrame
contains the predicted values for the conversion and activity models,
which are used to create the contour plots.

\begin{figure}

\centering{

\pandocbounded{\includegraphics[keepaspectratio]{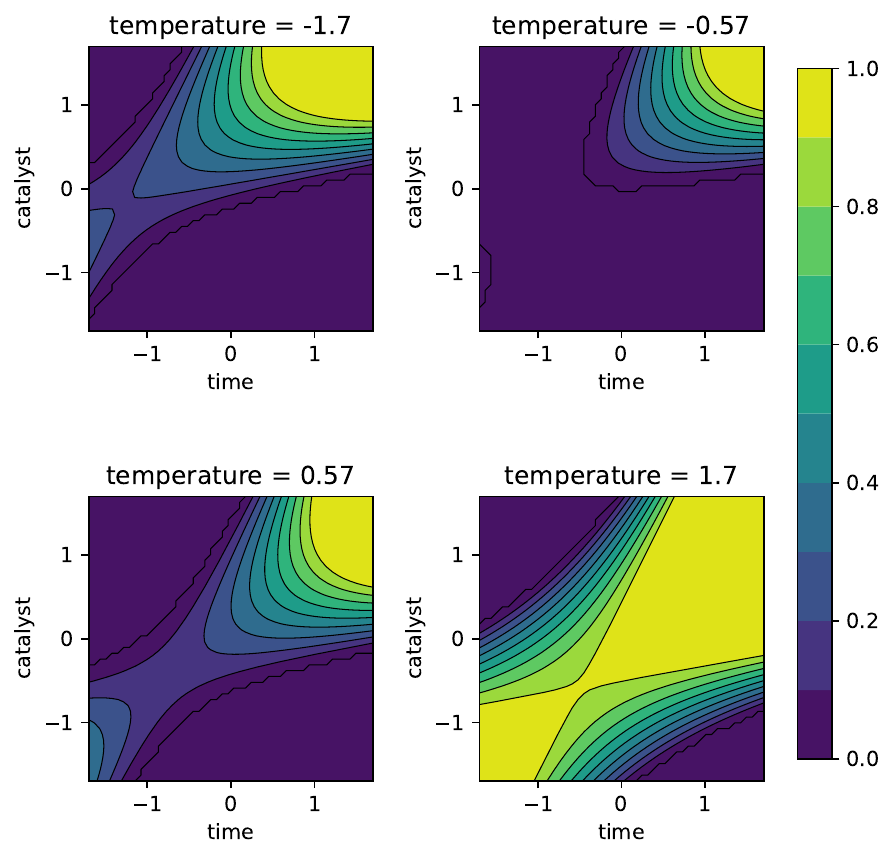}}

}

\caption{\label{fig-kuhn16a-4}The individual desirability surface for
the percent conversion outcome using \texttt{dMax(80,\ 97)}}

\end{figure}

\begin{figure}

\centering{

\pandocbounded{\includegraphics[keepaspectratio]{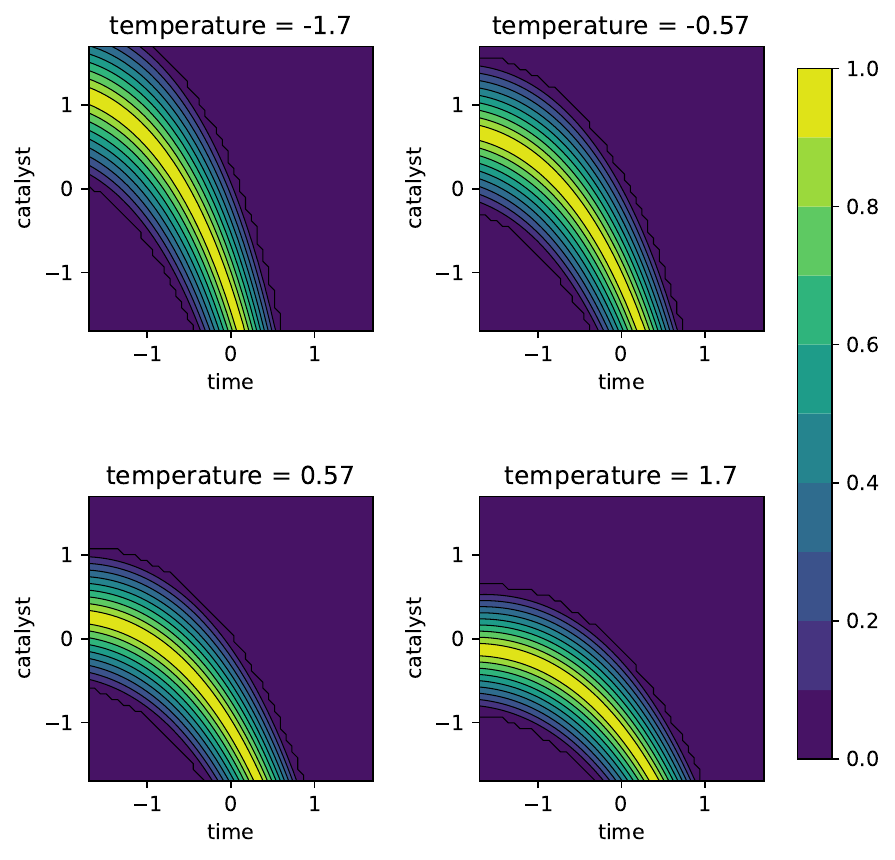}}

}

\caption{\label{fig-kuhn16a-5}The individual desirability surface for
the thermal activity outcome using \texttt{dTarget(55,\ 57.5,\ 60)}}

\end{figure}

\begin{figure}

\centering{

\pandocbounded{\includegraphics[keepaspectratio]{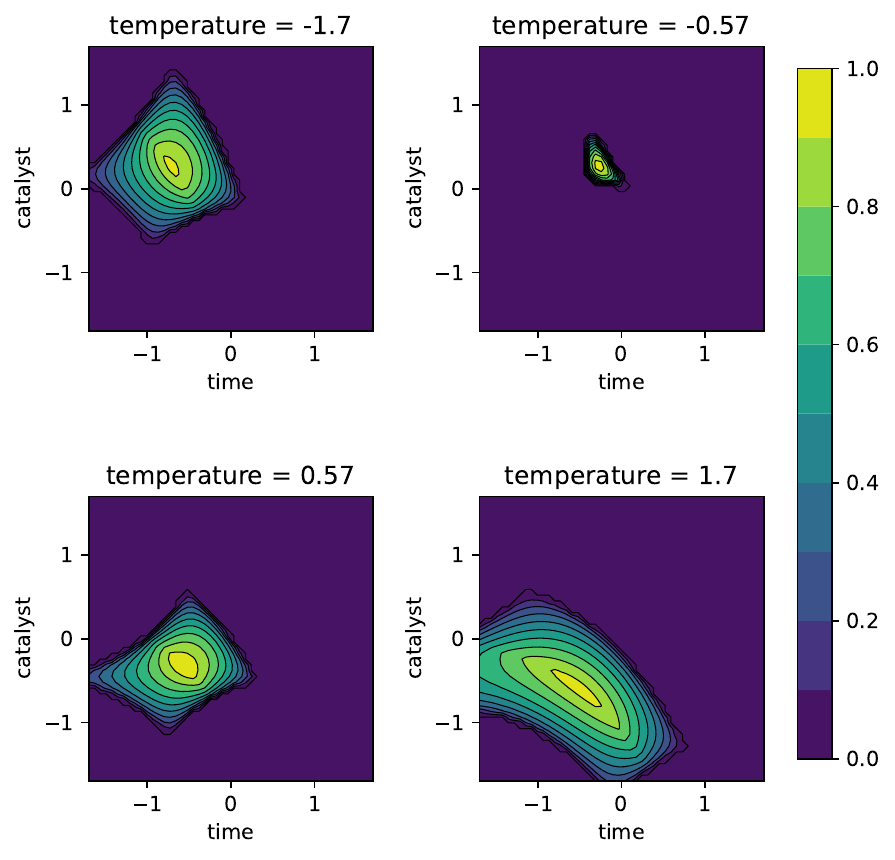}}

}

\caption{\label{fig-kuhn16a-6}The overall desirability surface for the
combined outcomes of percent conversion and thermal activity}

\end{figure}

Figure~\ref{fig-kuhn16a-4}, Figure~\ref{fig-kuhn16a-5}, and
Figure~\ref{fig-kuhn16a-6} show contour plots of the individual
desirability function surfaces and the overall surface. These plots are
in correspondence with the figures in Kuhn (2016), but the color schemes
are different. The \texttt{plot\_grid} DataFrame contains the predicted
values for the conversion and activity models, which are used to create
the contour plots. The individual desirability surface for the
\texttt{percent\ conversion} outcome is shown in
Figure~\ref{fig-kuhn16a-4} and the individual desirability surface for
the thermal activity outcome is shown in Figure~\ref{fig-kuhn16a-5}.
Finally, the overall desirability surface is shown in
Figure~\ref{fig-kuhn16a-6}.

\section{Multi-Objective Optimization and Maximizing
Desirability}\label{sec-maximizing-desirability}

Kuhn (2016) indicates that as described by Myers, Montgomery, and
Anderson-Cook (2016), desirability can be maximized within a cuboidal
region defined by the axial point values. The objective function, i.e.,
\texttt{rsm\_opt} from the \texttt{spotdesirability} package, utilizes a
penalty approach: if a candidate point extends beyond the cuboidal
design region, desirability is set to zero. These penalties are
implemented in the \texttt{rsm\_opt} function, which is used to optimize
the desirability function. An \(\alpha\) value of 1.682
(\(\approx (2^k)^{1/4}\) with \(k=3\) in our case), see Montgomery
(2001), is used as the limit for both circular and square spaces. After
checking the bounds, predictions for all provided functions are
calculated, and the overall desirability is predicted using the
\texttt{predict} method of the \texttt{DOverall} object. Since the
optimizer minimize the objective function, the negative desirability is
returned to maximize the desirability function, symbolically speaking:
\(\min f(-x) \Leftrightarrow \max f(x)\).

\phantomsection\label{kuhn16a-optimization-rsm}
\begin{Shaded}
\begin{Highlighting}[]
\KeywordTok{def}\NormalTok{ rsm\_opt(x, d\_object, prediction\_funcs, space}\OperatorTok{=}\StringTok{"square"}\NormalTok{, alpha}\OperatorTok{=}\FloatTok{1.682}\NormalTok{) }\OperatorTok{{-}\textgreater{}} \BuiltInTok{float}\NormalTok{:}
    \ControlFlowTok{if}\NormalTok{ space }\OperatorTok{==} \StringTok{"circular"}\NormalTok{:}
        \ControlFlowTok{if}\NormalTok{ np.sqrt(np.}\BuiltInTok{sum}\NormalTok{(np.array(x) }\OperatorTok{**} \DecValTok{2}\NormalTok{)) }\OperatorTok{\textgreater{}}\NormalTok{ alpha:}
            \ControlFlowTok{return} \FloatTok{0.0}
    \ControlFlowTok{elif}\NormalTok{ space }\OperatorTok{==} \StringTok{"square"}\NormalTok{:}
        \ControlFlowTok{if}\NormalTok{ np.}\BuiltInTok{any}\NormalTok{(np.}\BuiltInTok{abs}\NormalTok{(np.array(x)) }\OperatorTok{\textgreater{}}\NormalTok{ alpha):}
            \ControlFlowTok{return} \FloatTok{0.0}
    \ControlFlowTok{else}\NormalTok{:}
        \ControlFlowTok{raise} \PreprocessorTok{ValueError}\NormalTok{(}\StringTok{"space must be \textquotesingle{}square\textquotesingle{} or \textquotesingle{}circular\textquotesingle{}"}\NormalTok{)}
\NormalTok{    predictions }\OperatorTok{=}\NormalTok{ [func(x) }\ControlFlowTok{for}\NormalTok{ func }\KeywordTok{in}\NormalTok{ prediction\_funcs]}
\NormalTok{    desirability }\OperatorTok{=}\NormalTok{ d\_object.predict(np.array([predictions]))}
    \ControlFlowTok{return} \OperatorTok{{-}}\NormalTok{desirability}
\end{Highlighting}
\end{Shaded}

Note: Instead of using the penatlty approach, alternatively the
desirability function for box-constraints can be used. Furthermore,
\texttt{scipy.optimize} provides a \texttt{bounds} argument for some
optimizers to restrict the search space.

Kuhn (2016) used \texttt{R}'s \texttt{optim} function to implement the
Nelder-Mead simplex method (Nelder and Mead 1965; Olsson and Nelson
1975). This direct search method relies on function evaluations without
using gradient information. Although this method may converge to a local
optimum, it is fast with efficient functions, allowing for multiple
feasible region restarts to find the best result. We will use the
\texttt{scipy.optimize.minimize} function to implement the Nelder-Mead
simplex method in Python. Alternatively, methods like simulated
annealing (Bohachevsky 1986), also available in \texttt{R}'s
\texttt{optim} function, might better suit global optimum searches,
though they might need parameter tuning for effective performance.

Putting the pieces together, the following code segment shows how to
create the desirability objects and use them in the optimization
process. First, a \texttt{search\_grid} is created using
\texttt{numpy}'s \texttt{meshgrid} function to generate a grid of
restarts points in the design space. For each (restart) point in the
search grid, the \texttt{rsm\_opt} function is called to calculate the
desirability for that point. The \texttt{conversion\_pred} and
\texttt{activity\_pred} functions are used as prediction functions, and
the \texttt{DOverall} object is created using the individual
desirability objects for conversion and activity. The \texttt{overallD}
(overall desirability) is passed to tne \texttt{rsm\_opt} function. The
\texttt{minimize} function from \texttt{scipy.optimize} is used to find
the optimal parameters that minimize the negative desirability.

\phantomsection\label{kuhn16a-optimization}
\begin{Shaded}
\begin{Highlighting}[]
\CommentTok{\# Define the search grid}
\NormalTok{time }\OperatorTok{=}\NormalTok{ np.linspace(}\OperatorTok{{-}}\FloatTok{1.5}\NormalTok{, }\FloatTok{1.5}\NormalTok{, }\DecValTok{5}\NormalTok{)}
\NormalTok{temperature }\OperatorTok{=}\NormalTok{ np.linspace(}\OperatorTok{{-}}\FloatTok{1.5}\NormalTok{, }\FloatTok{1.5}\NormalTok{, }\DecValTok{5}\NormalTok{)}
\NormalTok{catalyst }\OperatorTok{=}\NormalTok{ np.linspace(}\OperatorTok{{-}}\FloatTok{1.5}\NormalTok{, }\FloatTok{1.5}\NormalTok{, }\DecValTok{5}\NormalTok{)}

\NormalTok{search\_grid }\OperatorTok{=}\NormalTok{ pd.DataFrame(}
\NormalTok{    np.array(np.meshgrid(time, temperature, catalyst)).T.reshape(}\OperatorTok{{-}}\DecValTok{1}\NormalTok{, }\DecValTok{3}\NormalTok{),}
\NormalTok{    columns}\OperatorTok{=}\NormalTok{[}\StringTok{"time"}\NormalTok{, }\StringTok{"temperature"}\NormalTok{, }\StringTok{"catalyst"}\NormalTok{]}
\NormalTok{)}

\CommentTok{\# List of prediction functions}
\NormalTok{prediction\_funcs }\OperatorTok{=}\NormalTok{ [conversion\_pred, activity\_pred]}

\CommentTok{\# Individual desirability objects}
\NormalTok{conversionD }\OperatorTok{=}\NormalTok{ DMax(}\DecValTok{80}\NormalTok{, }\DecValTok{97}\NormalTok{)}
\NormalTok{activityD }\OperatorTok{=}\NormalTok{ DTarget(}\DecValTok{55}\NormalTok{, }\FloatTok{57.5}\NormalTok{, }\DecValTok{60}\NormalTok{)}

\CommentTok{\# Desirability object (DOverall)}
\NormalTok{overallD }\OperatorTok{=}\NormalTok{ DOverall(conversionD, activityD)}

\CommentTok{\# Initialize the best result}
\NormalTok{best }\OperatorTok{=} \VariableTok{None}

\CommentTok{\# Perform optimization for each point in the search grid}
\ControlFlowTok{for}\NormalTok{ i, row }\KeywordTok{in}\NormalTok{ search\_grid.iterrows():}
\NormalTok{    initial\_guess }\OperatorTok{=}\NormalTok{ row.values  }\CommentTok{\# Initial guess for optimization}

    \CommentTok{\# Perform optimization using scipy\textquotesingle{}s minimize function}
    \CommentTok{\# rsm\_opt returns the negative desirability}
\NormalTok{    result }\OperatorTok{=}\NormalTok{ minimize(}
\NormalTok{        rsm\_opt,}
\NormalTok{        initial\_guess,}
\NormalTok{        args}\OperatorTok{=}\NormalTok{(overallD, prediction\_funcs, }\StringTok{"square"}\NormalTok{), }
\NormalTok{        method}\OperatorTok{=}\StringTok{"Nelder{-}Mead"}\NormalTok{,}
\NormalTok{        options}\OperatorTok{=}\NormalTok{\{}\StringTok{"maxiter"}\NormalTok{: }\DecValTok{1000}\NormalTok{, }\StringTok{"disp"}\NormalTok{: }\VariableTok{False}\NormalTok{\}}
\NormalTok{    )}

    \CommentTok{\# Update the best result if necessary}
    \CommentTok{\# Compare based on the negative desirability}
    \ControlFlowTok{if}\NormalTok{ best }\KeywordTok{is} \VariableTok{None} \KeywordTok{or}\NormalTok{ result.fun }\OperatorTok{\textless{}}\NormalTok{ best.fun:}
\NormalTok{        best }\OperatorTok{=}\NormalTok{ result}
\BuiltInTok{print}\NormalTok{(}\StringTok{"Best Parameters:"}\NormalTok{, best.x)}
\BuiltInTok{print}\NormalTok{(}\StringTok{"Best Desirability:"}\NormalTok{, }\OperatorTok{{-}}\NormalTok{best.fun)}
\end{Highlighting}
\end{Shaded}

\begin{verbatim}
Best Parameters: [-0.51207663  1.68199987 -0.58609664]
Best Desirability: 0.9425092694688632
\end{verbatim}

Because the optimizer ``sees'' only desirability and not the underlying
objective function, using the best input parameters found by the
optimizer, the predicted values for conversion and activity can be
calculated as follows:

\begin{Shaded}
\begin{Highlighting}[]
\BuiltInTok{print}\NormalTok{(}\SpecialStringTok{f"Conversion pred(x): }\SpecialCharTok{\{}\NormalTok{conversion\_pred(best.x)}\SpecialCharTok{\}}\SpecialStringTok{"}\NormalTok{)}
\BuiltInTok{print}\NormalTok{(}\SpecialStringTok{f"Activity pred(x): }\SpecialCharTok{\{}\NormalTok{activity\_pred(best.x)}\SpecialCharTok{\}}\SpecialStringTok{"}\NormalTok{)}
\end{Highlighting}
\end{Shaded}

\begin{verbatim}
Conversion pred(x): 95.10150374903237
Activity pred(x): 57.49999992427212
\end{verbatim}

Instead of generationg contour plots for four different
\texttt{temperature} values as in Figure~\ref{fig-kuhn16a-2} and
Figure~\ref{fig-kuhn16a-3}, we extract the best temperature from the
best parameters and remove it from the best parameters for plotting. The
\texttt{best.x} array contains the best parameters found by the
optimizer, where the second element corresponds to the temperature
variable. Then we set the values of \texttt{temperature} to the best
temperature in the \texttt{plot\_grid\_df} and recalculate the predicted
values for \texttt{conversion} and \texttt{activity} using the
\texttt{conversion\_pred} and \texttt{activity\_pred} functions. A copy
of the \texttt{plot\_grid} DataFrame is created, and the
\texttt{temperature} column is updated with the best temperature value.

\begin{figure}

\centering{

\pandocbounded{\includegraphics[keepaspectratio]{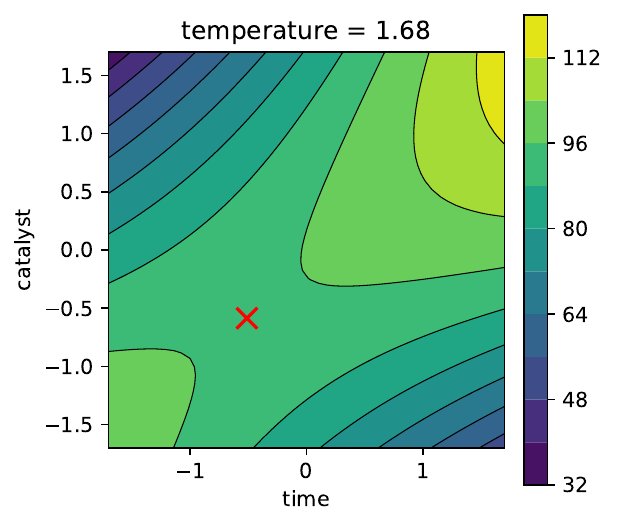}}

}

\caption{\label{fig-kuhn16a-best-conversion}The response surface for the
percent conversion model. To plot the model contours, the temperature
variable was fixed at the best value found by the optimizer.}

\end{figure}

\begin{figure}

\centering{

\pandocbounded{\includegraphics[keepaspectratio]{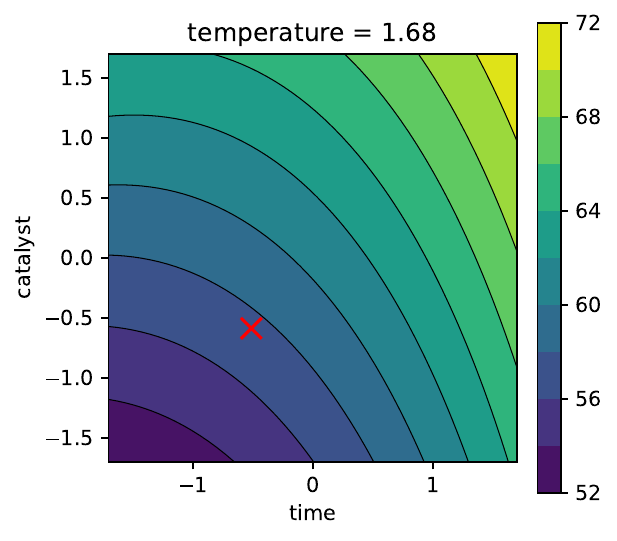}}

}

\caption{\label{fig-kuhn16a-best-activity}The response surface for the
thermal activity model. To plot the model contours, the temperature
variable was fixed at the best value found by the optimizer.}

\end{figure}

Now we are ready to plot the response surfaces for the best parameters
found by the optimizer. The \texttt{contourf\_plot} function is used to
create the contour plots for the response surface models. The
\texttt{highlight\_point} argument is used to highlight the best point
found by the optimizer in the contour plots. First, the response surface
for the \texttt{percent\ conversion} model is plotted. The temperature
variable is fixed at the best value found by the optimizer, see
Figure~\ref{fig-kuhn16a-best-conversion}. Second, the response surface
for the \texttt{thermal\ activity} model is plotted. The temperature
variable is fixed at the best value found by the optimizer, see
Figure~\ref{fig-kuhn16a-best-activity}.

\begin{tcolorbox}[enhanced jigsaw, colbacktitle=quarto-callout-important-color!10!white, opacityback=0, title=\textcolor{quarto-callout-important-color}{\faExclamation}\hspace{0.5em}{Analysing the Best Values From the Nelder-Mead Optimizer}, opacitybacktitle=0.6, arc=.35mm, left=2mm, bottomrule=.15mm, breakable, toprule=.15mm, coltitle=black, colframe=quarto-callout-important-color-frame, colback=white, rightrule=.15mm, bottomtitle=1mm, toptitle=1mm, titlerule=0mm, leftrule=.75mm]

\begin{itemize}
\tightlist
\item
  Objective function values for the best parameters found by the
  optimizer are:

  \begin{itemize}
  \tightlist
  \item
    \texttt{conversion} = 95.1
  \item
    \texttt{activity} = 57.5
  \end{itemize}
\item
  The best value for the percent conversion should be maximized, as
  defined in \texttt{conversionD\ =\ DMax(80,\ 97)}. Here, we have
  obtained a value of 95.1, which is close to the maximum value of 97.
\item
  Since we are using the desirabilty function \texttt{DTarget}, the
  values for the \texttt{thermal\ activity} should not be maximized, but
  should be close to the target. The setting
  \texttt{activityD\ =\ DTarget(55,\ 57.5,\ 60)}, as defined in
  Section~\ref{sec-defining-desirability}, states that the best value
  for the \texttt{thermal\ activity} should be close to 57.5 as
  specified by the user (and not at its maximum). Here, we have obtained
  a value of 57.5, which is exactly the target value.
\end{itemize}

\end{tcolorbox}

An alternative approach to the optimization process is to use a circular
design region instead of a cuboidal design region can be found in the
Appendix.

\section{Surrogate-Model Based Optimization Using
Desirability}\label{sec-surrogate}

In Section~\ref{sec-maximizing-desirability}, we have used a classical
optimization approach to find the best parameters for the response
surface models. This approach was used to demonstrate that the
\texttt{Python} code from \texttt{spotdesirability} leads to similar
results as Kuhn (2016)'s \texttt{R} code. Now we will use a
surrogate-model based optimization approach to find the best parameters
for the response surface models. The \texttt{spotoptim} package
implements a vectorized function \texttt{fun\_myer16a()} that computes
the two objective functions for \texttt{conversion} and
\texttt{activity}. To illustrate the vectorized evaluation, we will use
two input points: the center point of the design space and the best
point found by the optimizer from
Section~\ref{sec-maximizing-desirability}. The \texttt{fun\_myer16a()}
function takes a 2D array as input, where each row corresponds to a
different set of parameters. The function returns a 2D array with the
predicted values for \texttt{conversion} and \texttt{activity}.

\begin{Shaded}
\begin{Highlighting}[]
\NormalTok{X }\OperatorTok{=}\NormalTok{ np.array([[}\DecValTok{0}\NormalTok{, }\DecValTok{0}\NormalTok{, }\DecValTok{0}\NormalTok{], best.x])}
\NormalTok{y }\OperatorTok{=}\NormalTok{ fun\_myer16a(X)}
\BuiltInTok{print}\NormalTok{(}\SpecialStringTok{f"Objective function values:}\CharTok{\textbackslash{}n}\SpecialCharTok{\{}\NormalTok{y}\SpecialCharTok{\}}\SpecialStringTok{"}\NormalTok{)}
\end{Highlighting}
\end{Shaded}

\begin{verbatim}
Objective function values:
[[81.09       59.85      ]
 [95.10150375 57.49999992]]
\end{verbatim}

Next, we define the desirability objects. This step is identical to the
previous one, where we defined the desirability functions for
\texttt{conversion} and \texttt{activity}, see
Figure~\ref{fig-kuhn16a-desirability-plot-conversion} and
Figure~\ref{fig-kuhn16a-desirability-plot-activity}. The \texttt{DMax}
function is used for the \texttt{conversion} function, and the
\texttt{DTarget} function is used for the \texttt{activity} function.
The \texttt{DOverall} function is used to combine the two desirability
functions into an overall desirability function. The \texttt{DOverall}
function takes two arguments: the desirability object for
\texttt{conversion} and the desirability object for \texttt{activity}.

\begin{Shaded}
\begin{Highlighting}[]
\NormalTok{conversionD }\OperatorTok{=}\NormalTok{ DMax(}\DecValTok{80}\NormalTok{, }\DecValTok{97}\NormalTok{)}
\NormalTok{activityD }\OperatorTok{=}\NormalTok{ DTarget(}\DecValTok{55}\NormalTok{, }\FloatTok{57.5}\NormalTok{, }\DecValTok{60}\NormalTok{)}
\NormalTok{overallD }\OperatorTok{=}\NormalTok{ DOverall(conversionD, activityD)}
\end{Highlighting}
\end{Shaded}

Predicting the desirability for each outcome can also be vectorized. The
\texttt{predict} method of the desirability objects can take a 2D array
as input, where each row corresponds to a different set of parameters.
The method returns a 1D array with the predicted desirability values for
each set of parameters.

\begin{Shaded}
\begin{Highlighting}[]
\NormalTok{conversion\_desirability }\OperatorTok{=}\NormalTok{ conversionD.predict(y[:,}\DecValTok{0}\NormalTok{])}
\NormalTok{activity\_desirability }\OperatorTok{=}\NormalTok{ activityD.predict(y[:,}\DecValTok{1}\NormalTok{])}
\end{Highlighting}
\end{Shaded}

\begin{verbatim}
Conversion Desirability: [0.06411765 0.88832375]
Activity Desirability: [0.06       0.99999997]
\end{verbatim}

The first value in each of the two arrays represents the desirability
for the center point of the design space, while the second value
represents the desirability for the best point found by the optimizer.
The \texttt{overall\_desirability} variable contains the overall
desirability values for each set of parameters. The \texttt{all=True}
argument indicates that we want to return both the individual
desirability values and the overall desirability value. If
\texttt{all=\ False}, only the overall desirability value is returned.

\begin{Shaded}
\begin{Highlighting}[]
\NormalTok{overall\_desirability }\OperatorTok{=}\NormalTok{ overallD.predict(y, }\BuiltInTok{all}\OperatorTok{=}\VariableTok{False}\NormalTok{)}
\end{Highlighting}
\end{Shaded}

\begin{verbatim}
OverallD: [0.06202466 0.94250927]
\end{verbatim}

During the surrogate-model based optimization, the argument \texttt{all}
is set to \texttt{False}, because \texttt{spotoptim} does not need the
individual desirability values. Now we have introduced all elements
needed to perform surrogate-model based optimization using desirability
functions and the \texttt{spotoptim} package.

\begin{tcolorbox}[enhanced jigsaw, colbacktitle=quarto-callout-important-color!10!white, opacityback=0, title=\textcolor{quarto-callout-important-color}{\faExclamation}\hspace{0.5em}{Maximization and Minimization}, opacitybacktitle=0.6, arc=.35mm, left=2mm, bottomrule=.15mm, breakable, toprule=.15mm, coltitle=black, colframe=quarto-callout-important-color-frame, colback=white, rightrule=.15mm, bottomtitle=1mm, toptitle=1mm, titlerule=0mm, leftrule=.75mm]

\begin{itemize}
\tightlist
\item
  Since \texttt{spotoptim} uses minimization, but desirability should be
  maximized, \texttt{fun\_desirability} is defined to return
  \texttt{1\ -\ overall\_desirability}.
\end{itemize}

\end{tcolorbox}

\begin{Shaded}
\begin{Highlighting}[]
\KeywordTok{def}\NormalTok{ fun\_desirability(X, }\OperatorTok{**}\NormalTok{kwargs):}
\NormalTok{    y }\OperatorTok{=}\NormalTok{ fun\_myer16a(X)}
\NormalTok{    conversionD }\OperatorTok{=}\NormalTok{ DMax(}\DecValTok{80}\NormalTok{, }\DecValTok{97}\NormalTok{)}
\NormalTok{    activityD }\OperatorTok{=}\NormalTok{ DTarget(}\DecValTok{55}\NormalTok{, }\FloatTok{57.5}\NormalTok{, }\DecValTok{60}\NormalTok{)}
\NormalTok{    overallD }\OperatorTok{=}\NormalTok{ DOverall(conversionD, activityD)}
\NormalTok{    overall\_desirability }\OperatorTok{=}\NormalTok{ overallD.predict(y, }\BuiltInTok{all}\OperatorTok{=}\VariableTok{False}\NormalTok{)}
    \ControlFlowTok{return} \FloatTok{1.0} \OperatorTok{{-}}\NormalTok{ overall\_desirability}
\end{Highlighting}
\end{Shaded}

\begin{example}[Testing the Desirability
Function]\protect\hypertarget{exm-testfun}{}\label{exm-testfun}

We can test the function, which returns the two overall ``1 minus
desirability''-function values for the center point of the design space
and the best point found by the optimizer. Obviously, the best point
found by the optimizer has a lower objective function value than the
center point and therefore a higher desirability.

\begin{Shaded}
\begin{Highlighting}[]
\NormalTok{X }\OperatorTok{=}\NormalTok{ np.array([[}\DecValTok{0}\NormalTok{, }\DecValTok{0}\NormalTok{, }\DecValTok{0}\NormalTok{], best.x])}
\NormalTok{y }\OperatorTok{=}\NormalTok{ fun\_desirability(X)}
\BuiltInTok{print}\NormalTok{(}\SpecialStringTok{f"Objective function values: }\SpecialCharTok{\{}\NormalTok{y}\SpecialCharTok{\}}\SpecialStringTok{"}\NormalTok{)}
\end{Highlighting}
\end{Shaded}

\begin{verbatim}
Objective function values: [0.93797534 0.05749073]
\end{verbatim}

\end{example}

We are now ready to perform the surrogate-model based optimization using
desirability functions. The \texttt{spotoptim} package provides a class
\texttt{SpotOptim} that implements the surrogate-model based
optimization algorithm. The \texttt{SpotOptim} class takes the objective
function and the control parameters as input. It uses the same interface
as the \texttt{scipy.optimize.minimize} function, see
\url{https://docs.scipy.org/doc/scipy/reference/generated/scipy.optimize.minimize.html}.
The control parameters define the search space and other settings for
the optimization process. As the surrogate-model, we use a Gaussian
process regressor with a Matern kernel.

\begin{Shaded}
\begin{Highlighting}[]
\NormalTok{kernel }\OperatorTok{=}\NormalTok{ ConstantKernel(}\FloatTok{1.0}\NormalTok{, (}\FloatTok{1e{-}2}\NormalTok{, }\FloatTok{1e12}\NormalTok{)) }\OperatorTok{*}\NormalTok{ Matern(length\_scale}\OperatorTok{=}\FloatTok{1.0}\NormalTok{, length\_scale\_bounds}\OperatorTok{=}\NormalTok{(}\FloatTok{1e{-}4}\NormalTok{, }\FloatTok{1e2}\NormalTok{), nu}\OperatorTok{=}\FloatTok{2.5}\NormalTok{)}
\NormalTok{S\_GP }\OperatorTok{=}\NormalTok{ GaussianProcessRegressor(kernel}\OperatorTok{=}\NormalTok{kernel, n\_restarts\_optimizer}\OperatorTok{=}\DecValTok{100}\NormalTok{)}

\NormalTok{S }\OperatorTok{=}\NormalTok{ SpotOptim(fun}\OperatorTok{=}\NormalTok{fun\_desirability,         }
\NormalTok{         bounds}\OperatorTok{=}\NormalTok{[(}\OperatorTok{{-}}\FloatTok{1.7}\NormalTok{, }\FloatTok{1.7}\NormalTok{), (}\OperatorTok{{-}}\FloatTok{1.7}\NormalTok{, }\FloatTok{1.7}\NormalTok{), (}\OperatorTok{{-}}\FloatTok{1.7}\NormalTok{, }\FloatTok{1.7}\NormalTok{)],}
\NormalTok{         max\_iter}\OperatorTok{=}\DecValTok{50}\NormalTok{,}
\NormalTok{         max\_time}\OperatorTok{=}\NormalTok{inf,}
\NormalTok{         n\_initial}\OperatorTok{=}\DecValTok{15}\NormalTok{,}
\NormalTok{         var\_name}\OperatorTok{=}\NormalTok{[}\StringTok{"time"}\NormalTok{, }\StringTok{"temperature"}\NormalTok{, }\StringTok{"catalyst"}\NormalTok{],}
\NormalTok{         seed}\OperatorTok{=}\DecValTok{126}\NormalTok{,}
\NormalTok{         surrogate}\OperatorTok{=}\NormalTok{S\_GP,}
\NormalTok{         max\_surrogate\_points}\OperatorTok{=}\DecValTok{30}\NormalTok{,}
\NormalTok{         )}
\NormalTok{S.optimize()}
\end{Highlighting}
\end{Shaded}

\phantomsection\label{kuhn16a-spot}
\begin{verbatim}
 message: Optimization terminated: maximum evaluations (50) reached
                   Current function value: 0.052569
                   Iterations: 35
                   Function evaluations: 50
 success: True
     fun: 0.052568867873652114
       x: [-5.066e-01  1.700e+00 -5.873e-01]
       X: [[-3.414e-01  5.799e-01 -5.817e-01]
           [-1.469e+00  5.388e-01  9.967e-01]
           ...
           [-4.907e-01  1.700e+00 -6.224e-01]
           [-4.907e-01  1.700e+00 -6.224e-01]]
     nit: 35
    nfev: 50
       y: [ 6.923e-01  1.000e+00 ...  5.812e-02  5.812e-02]
\end{verbatim}

\begin{figure}[H]

\centering{

\pandocbounded{\includegraphics[keepaspectratio]{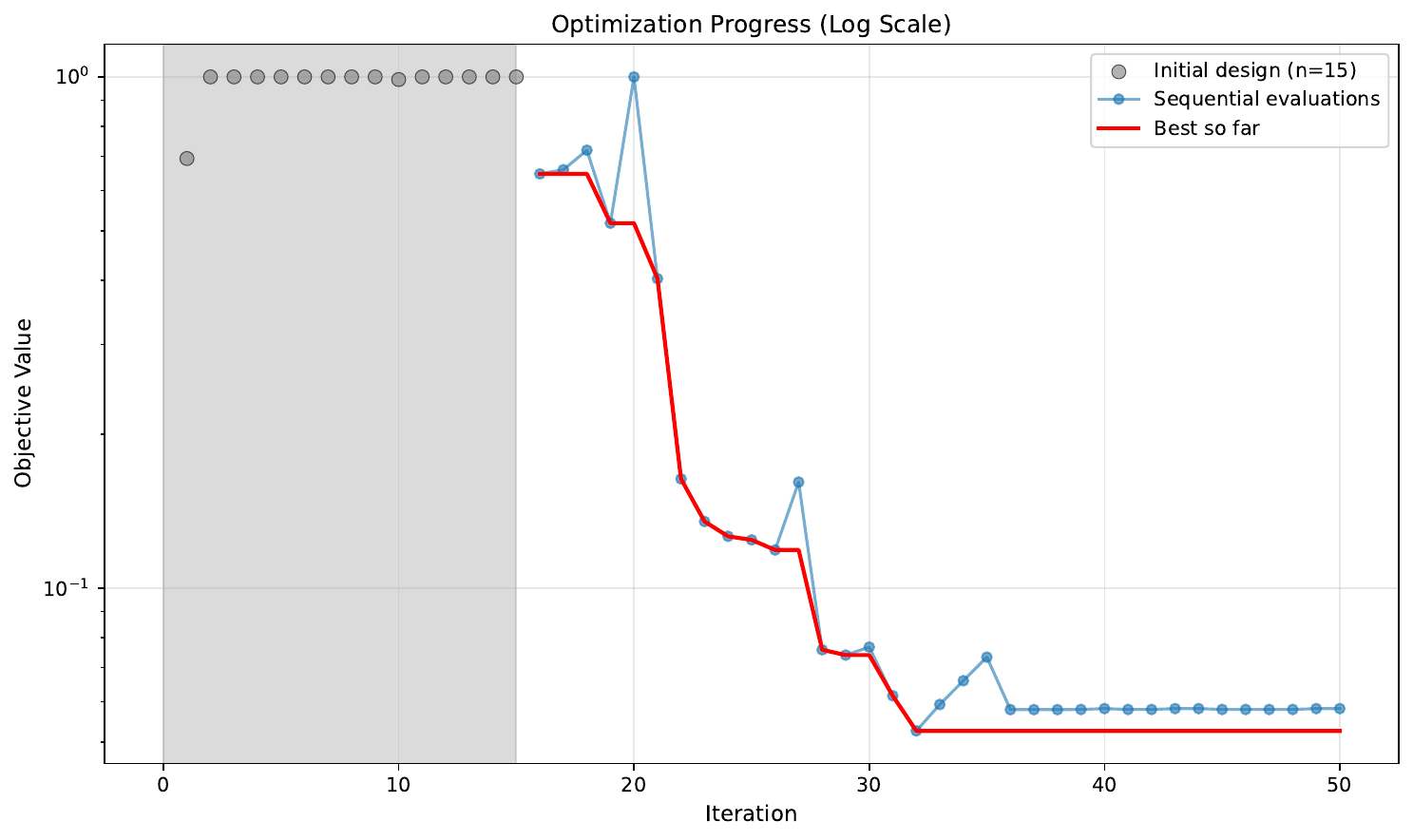}}

}

\caption{\label{fig-kuhn16a-spot-progress}The progress of the
surrogate-model based optimization using desirability functions. The
y-axis is on a logarithmic scale and represents 1 minus the overall
desirability.}

\end{figure}

The progress of the optimization process can be visualized using the
\texttt{plot\_progress} method (Figure~\ref{fig-kuhn16a-spot-progress}).

\begin{Shaded}
\begin{Highlighting}[]
\NormalTok{S.plot\_progress(log\_y}\OperatorTok{=}\VariableTok{True}\NormalTok{)}
\end{Highlighting}
\end{Shaded}

We can analyze the results in detail by accessing the attributes of the
\texttt{SpotOptim} object directly. The \texttt{min\_X} attribute
contains the best parameters found by the optimizer, and the
\texttt{min\_y} attribute contains the best desirability value. First,
we take a look at the desirability values for the best parameters found
by the optimizer. The \texttt{min\_y} attribute contains the best
desirability value. Note, we have to compute 1 minus the \texttt{min\_y}
value, because the \texttt{fun\_desirability} function returns
\texttt{1\ -\ overall\_desirability}. This results in the following best
desirability value, where values close to 1 indicate a better solution:

\begin{verbatim}
Best Desirability: 0.9474311321263479
\end{verbatim}

We can use the \texttt{min\_X} attribute to calculate the predicted
values for \texttt{conversion} and \texttt{activity} for the best
parameters found by the optimizer. Using the \texttt{fun\_myer16a}
function, we can calculate these predicted values. \texttt{SpotOptim}
does not use these ``intermediate'' values that are calculated by the
objective function and passed to the desirability function, because it
only needs the desirability values, i.e., \texttt{time},
\texttt{temperature}, and \texttt{catalyst}.

\begin{verbatim}
Best Parameters: [-0.50656334  1.7        -0.58731616]
Best Conversion: 95.3736636847183
Best Activity: 57.51854241366585
\end{verbatim}

\begin{tcolorbox}[enhanced jigsaw, colbacktitle=quarto-callout-important-color!10!white, opacityback=0, title=\textcolor{quarto-callout-important-color}{\faExclamation}\hspace{0.5em}{Analysing the Best Values from the spotoptim Optimizer}, opacitybacktitle=0.6, arc=.35mm, left=2mm, bottomrule=.15mm, breakable, toprule=.15mm, coltitle=black, colframe=quarto-callout-important-color-frame, colback=white, rightrule=.15mm, bottomtitle=1mm, toptitle=1mm, titlerule=0mm, leftrule=.75mm]

\begin{itemize}
\tightlist
\item
  Objective function values for the best parameters found by the
  optimizer are very close to the values found by the Nelder-Mead
  optimizer, see Section~\ref{sec-maximizing-desirability}.
\end{itemize}

\end{tcolorbox}

\begin{figure}

\centering{

\pandocbounded{\includegraphics[keepaspectratio]{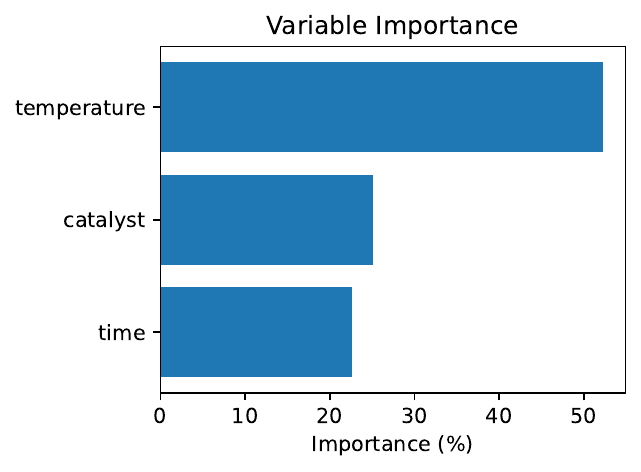}}

}

\caption{\label{fig-spot-importance}The importance of the parameters in
the optimization process based on desirability functions.}

\end{figure}

\begin{figure}

\centering{

\pandocbounded{\includegraphics[keepaspectratio]{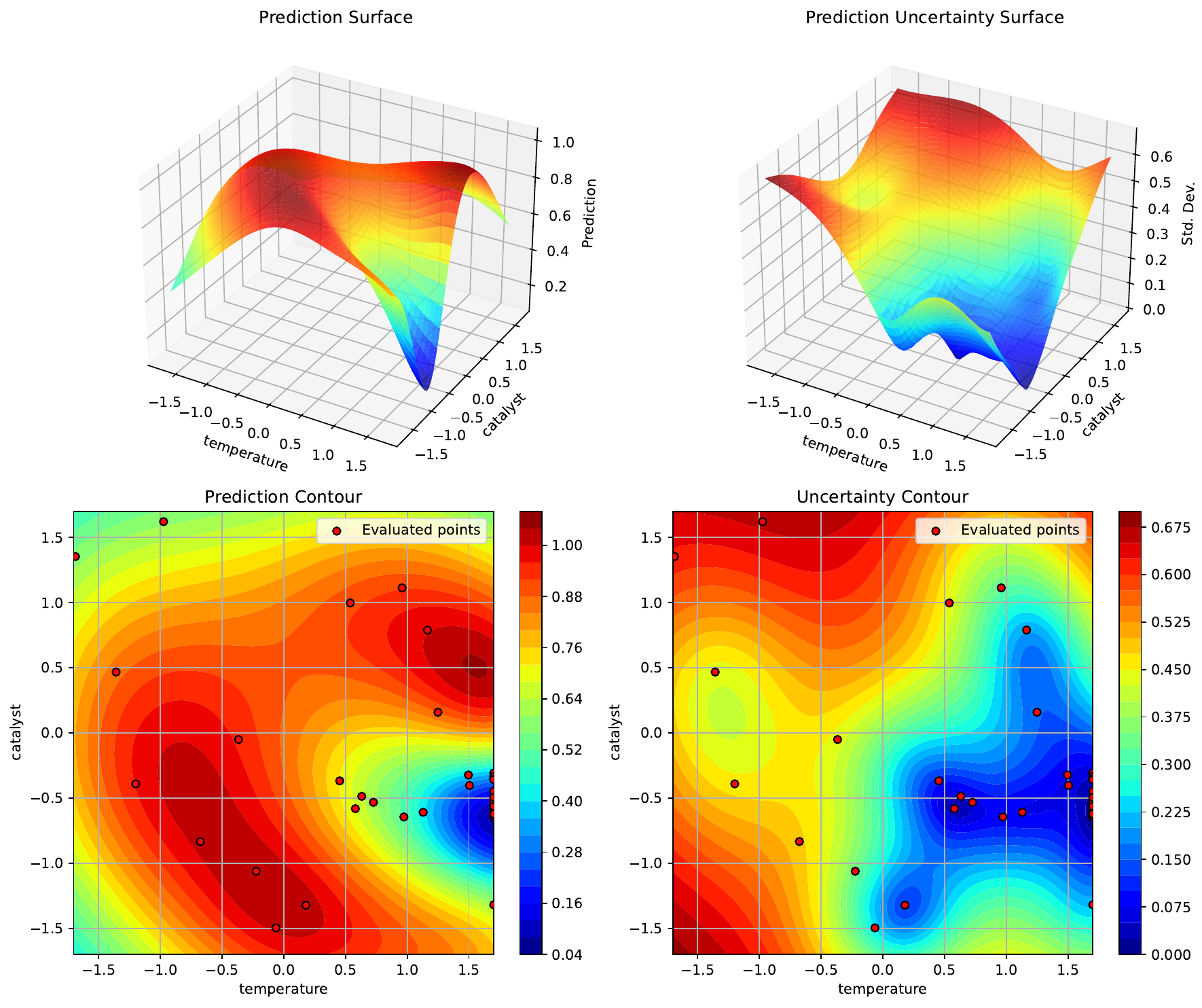}}

}

\caption{\label{fig-spot-importance-contour}Chemical-reactor
optimization: The contour plots for the two most important parameters.
Note, SpotOptim is minimizing the objective function, so we return (1 -
desirability).}

\end{figure}

Based on the information from the surrogate, which is the
\texttt{GaussianProcessRegressor} in \texttt{spotoptim}, we can analyze
the importance of the parameters in the optimization process. The
\texttt{plot\_importance} method plots the importance of each parameter,
see Figure~\ref{fig-spot-importance}. Interestingly, in contrast to the
response surface model, the \texttt{temperature} parameter is the most
important parameter in the optimization process with the
\texttt{GaussianProcessRegressor} surrogate. The
\texttt{plot\_important\_hyperparameter\_contour} method generates the
contour plots for the important parameters. The results are shown in
Figure~\ref{fig-spot-importance-contour}. The contour plots show the
importance of the parameters in the optimization process, which tries to
minimize the 1 minus desirability values. Regions with low values
present high desirability. Note: These surface plots illustrate how the
surrogate ``sees the world'' and decides where to sample next to compute
the next infill point.

\begin{verbatim}
Plotting surrogate contours for top 2 most important parameters:
  temperature: importance = 52.30% (type: float)
  catalyst: importance = 25.09% (type: float)

Generating 1 surrogate plots...
  Plotting temperature vs catalyst
\end{verbatim}

\begin{figure}

\centering{

\pandocbounded{\includegraphics[keepaspectratio]{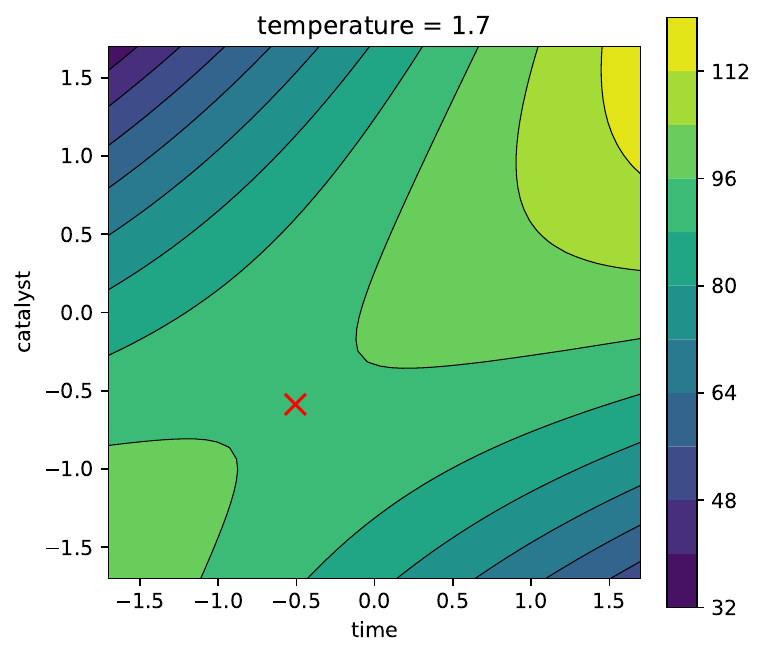}}

}

\caption{\label{fig-kuhn16a-2-surrogate}The response surface for the
percent conversion model. To plot the model contours, the temperature
variable was fixed at its best value.}

\end{figure}

\begin{figure}

\centering{

\pandocbounded{\includegraphics[keepaspectratio]{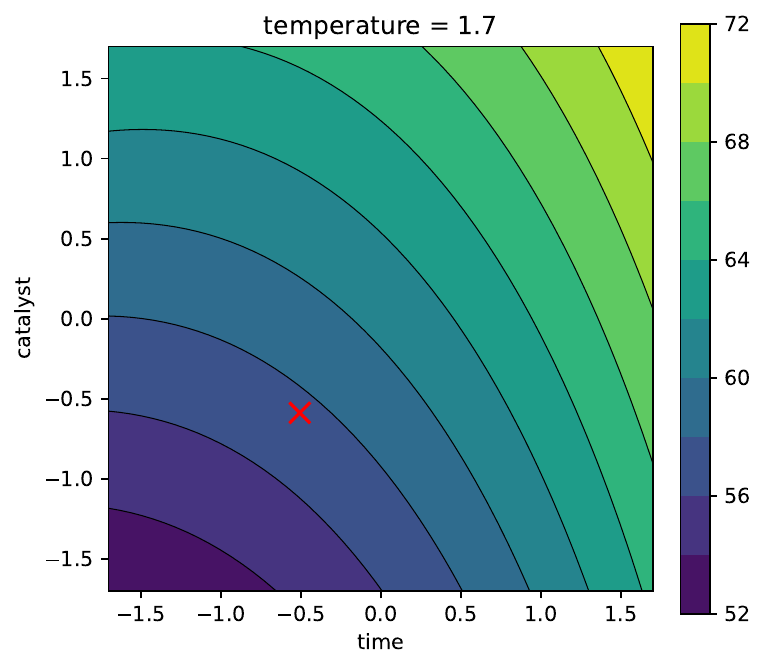}}

}

\caption{\label{fig-kuhn16a-surrogate-3}The response surface for the
thermal activity model. To plot the model contours, the temperature
variable was fixed at its best value.}

\end{figure}

Finally, we show a comparison with the response-surface model. Similar
to the procedure above that generated
Figure~\ref{fig-kuhn16a-best-conversion} and
Figure~\ref{fig-kuhn16a-best-activity}, we use the \texttt{plot\_grid}
DataFrame to generate the contour plots shown in
Figure~\ref{fig-kuhn16a-2-surrogate} and
Figure~\ref{fig-kuhn16a-surrogate-3}. To plot the model contours, the
temperature variable was fixed at its best value from the surrogate
model optimization. A comparison of the four contour plots shows that
similar results are obtained with both approaches.

\section{Hyperparameter Tuning Using Surrogate
Models}\label{sec-hyperparameter-tuning}

This section compares three different approaches to hyperparameter
tuning using the \texttt{spotoptim} package. The first approach, which
is shown in Section~\ref{sec-hpt-single}, is a single-objective
approach, where only the first objective function is used for
hyperparameter tuning. The second approach, which is shown in
Section~\ref{sec-hpt-weighted}, is a weighted multi-objective approach,
where a weighted mean of both objective functions is used for
hyperparameter tuning. The third approach, which is shown in
Section~\ref{sec-hpt-desirability}, uses a desirability function to
combine the two objective functions into a single objective function.
The desirability function is used to maximize the desirability of the
two objective functions.

The \texttt{spotoptim} package provides a method for hyperparameter
tuning using a surrogate model. We will extend the single-objective
example ``Hyperparameter Tuning with spotpython and PyTorch Lightning
for the Diabetes Data Set'' from the hyperparameter tuning cookbook
(Bartz-Beielstein 2023) in the follwoing way: Instead of using a
single-objective function, which returns the \texttt{validation\ loss}
from the neural-network training, we will use a multi-objective function
that returns two objectives: (i) \texttt{validation\ loss} and (ii)
\texttt{number\ of\ epochs}. Clearly, both objectives should be
minimized. The \texttt{validation\ loss} should be minimized to get the
best model, and the \texttt{number\ of\ epochs} should be minimized to
reduce the training time. However, if the number of training epochs is
too small, the model will not be trained properly. Therefore, we will
adopt the desirability function for the \texttt{number\ of\ epochs}
accordingly. The \texttt{Diabetes} data set from \texttt{scikit-learn}
is used for illustration in this section. It is a regression data set
that contains ten input features and a single output feature. The goal
is to predict the output feature based on the input features.

\subsection{The Single-Objective Approach}\label{sec-hpt-single}

The simplest way for handling multi-objective results is to simply
ignore all but the first objective function. This is the default
behavior of the \texttt{SpotOptim} class if the objective function,
which is specified by the \texttt{fun} argument, returns more than one
objective.

\subsubsection{Data Loading}\label{data-loading}

We use \texttt{SpotDataFromArray} to unify data handling. Here we load
the standard Diabetes dataset from \texttt{scikit-learn} and create a
\texttt{SpotDataSet}.

\begin{Shaded}
\begin{Highlighting}[]
\NormalTok{diabetes }\OperatorTok{=}\NormalTok{ load\_diabetes()}
\NormalTok{X }\OperatorTok{=}\NormalTok{ diabetes.data.astype(np.float32)}
\NormalTok{y }\OperatorTok{=}\NormalTok{ diabetes.target.reshape(}\OperatorTok{{-}}\DecValTok{1}\NormalTok{, }\DecValTok{1}\NormalTok{).astype(np.float32)}
\NormalTok{data }\OperatorTok{=}\NormalTok{ SpotDataFromArray(x\_train}\OperatorTok{=}\NormalTok{X, y\_train}\OperatorTok{=}\NormalTok{y)}
\BuiltInTok{print}\NormalTok{(}\SpecialStringTok{f"Data shape: }\SpecialCharTok{\{}\NormalTok{X}\SpecialCharTok{.}\NormalTok{shape}\SpecialCharTok{\}}\SpecialStringTok{ {-}\textgreater{} }\SpecialCharTok{\{}\NormalTok{y}\SpecialCharTok{.}\NormalTok{shape}\SpecialCharTok{\}}\SpecialStringTok{"}\NormalTok{)}
\end{Highlighting}
\end{Shaded}

\begin{verbatim}
Data shape: (442, 10) -> (442, 1)
\end{verbatim}

\subsubsection{Hyperparameter Setup}\label{hyperparameter-setup}

\texttt{spotoptim} provides a \texttt{MLP} class for training a
multi-layer perceptron (MLP) using PyTorch. From the \texttt{MLP} class
we get the default hyperparameters via the
\texttt{get\_default\_parameters} method. The \texttt{ParameterSet}
class allows for a fluid definition of the search space. We add the
\texttt{epochs} hyperparameter to the search space, because it is not
part of the default hyperparameters of the MLP class.

\begin{Shaded}
\begin{Highlighting}[]
\NormalTok{params }\OperatorTok{=}\NormalTok{ MLP.get\_default\_parameters()}
\NormalTok{params.add\_int(}\StringTok{"epochs"}\NormalTok{, low}\OperatorTok{=}\DecValTok{32}\NormalTok{, high}\OperatorTok{=}\DecValTok{1024}\NormalTok{, transform}\OperatorTok{=}\StringTok{"log(x, 2)"}\NormalTok{)}
\BuiltInTok{print}\NormalTok{(}\SpecialStringTok{f"Hyperparameters with epochs: }\SpecialCharTok{\{}\NormalTok{params}\SpecialCharTok{\}}\SpecialStringTok{"}\NormalTok{)}
\end{Highlighting}
\end{Shaded}

\begin{verbatim}
Hyperparameters with epochs: ParameterSet(
    l1=Parameter(
            name='l1',
            var_name='l1',
            bounds=Bounds(low=16, high=128),
            default=64,
            transform='log(x, 2)',
            type='int'
        ),
    num_hidden_layers=Parameter(
            name='num_hidden_layers',
            var_name='num_hidden_layers',
            bounds=Bounds(low=1, high=5),
            default=3,
            transform=None,
            type='int'
        ),
    activation=Parameter(
            name='activation',
            var_name='activation',
            bounds=['ReLU', 'Tanh', 'Sigmoid', 'LeakyReLU', 'ELU'],
            default='ReLU',
            transform=None,
            type='factor'
        ),
    lr=Parameter(
            name='lr',
            var_name='lr',
            bounds=Bounds(low=0.0001, high=100.0),
            default=10.0,
            transform='log',
            type='float'
        ),
    optimizer=Parameter(
            name='optimizer',
            var_name='optimizer',
            bounds=['Adam', 'SGD', 'RMSprop', 'AdamW'],
            default='Adam',
            transform=None,
            type='factor'
        ),
    epochs=Parameter(
            name='epochs',
            var_name='epochs',
            bounds=Bounds(low=32, high=1024),
            default=None,
            transform='log(x, 2)',
            type='int'
        ),
)
\end{verbatim}

\subsubsection{Experiment Configuration}\label{experiment-configuration}

\texttt{ExperimentControl} holds all configuration details, e.g.~the
dataset, the model class, the hyperparameters, the batch size, and the
metrics.

\begin{Shaded}
\begin{Highlighting}[]
\NormalTok{exp }\OperatorTok{=}\NormalTok{ ExperimentControl(}
\NormalTok{    dataset}\OperatorTok{=}\NormalTok{data,}
\NormalTok{    model\_class}\OperatorTok{=}\NormalTok{MLP,}
\NormalTok{    hyperparameters}\OperatorTok{=}\NormalTok{params,}
\NormalTok{    batch\_size}\OperatorTok{=}\DecValTok{32}\NormalTok{,}
\NormalTok{    metrics}\OperatorTok{=}\NormalTok{[}\StringTok{"mse"}\NormalTok{]}
\NormalTok{)}
\end{Highlighting}
\end{Shaded}

\subsubsection{Objective Function}\label{objective-function}

\texttt{TorchObjective} bridges the configuration with the optimization
loop, handling model instantiation and training. It contains, beside the
objective function, also the parameters specified in the
\texttt{ExperimentControl}. Since the \texttt{ExperimentControl} is
initialized with one metric value, the \texttt{TorchObjective} also
returns one metric value as shown in the following example.

\begin{Shaded}
\begin{Highlighting}[]
\NormalTok{objective }\OperatorTok{=}\NormalTok{ TorchObjective(experiment}\OperatorTok{=}\NormalTok{exp, seed}\OperatorTok{=}\DecValTok{42}\NormalTok{)}
\NormalTok{X\_demo }\OperatorTok{=}\NormalTok{ np.array([[}\FloatTok{0.1}\NormalTok{, }\FloatTok{0.2}\NormalTok{, }\FloatTok{0.3}\NormalTok{, }\FloatTok{0.4}\NormalTok{, }\FloatTok{0.5}\NormalTok{, }\FloatTok{9.0}\NormalTok{]])}
\NormalTok{hps }\OperatorTok{=}\NormalTok{ objective.\_get\_hyperparameters(X}\OperatorTok{=}\NormalTok{X\_demo)}
\BuiltInTok{print}\NormalTok{(}\SpecialStringTok{f"Hyperparameters: }\SpecialCharTok{\{}\NormalTok{hps}\SpecialCharTok{\}}\SpecialStringTok{"}\NormalTok{)}
\NormalTok{y\_eval }\OperatorTok{=}\NormalTok{ objective(X\_demo)}
\BuiltInTok{print}\NormalTok{(}\SpecialStringTok{f"y\_eval: }\SpecialCharTok{\{}\NormalTok{y\_eval}\SpecialCharTok{\}}\SpecialStringTok{"}\NormalTok{)}
\end{Highlighting}
\end{Shaded}

\begin{verbatim}
Hyperparameters: {'l1': 0, 'num_hidden_layers': 0, 'activation': 'ReLU', 'lr': np.float64(0.4), 'optimizer': 'Adam', 'epochs': 9}
y_eval: [[28904.8406808]]
\end{verbatim}

\subsubsection{Run Optimization}\label{run-optimization}

We initialize \texttt{SpotOptim} using the properties from our
\texttt{ParameterSet}, which are stored in the \texttt{objective}
object. Since \texttt{bounds}, \texttt{var\_type}, \texttt{var\_names},
and \texttt{var\_trans} are already part of the \texttt{objective}
object, we do not need to pass them as arguments to the
\texttt{SpotOptim} constructor. The we run the optimization via the
\texttt{optimize} method.

\begin{Shaded}
\begin{Highlighting}[]
\NormalTok{so\_optimizer }\OperatorTok{=}\NormalTok{ SpotOptim(}
\NormalTok{    fun}\OperatorTok{=}\NormalTok{objective,}
\NormalTok{    max\_iter}\OperatorTok{=}\DecValTok{50}\NormalTok{,}
\NormalTok{    max\_time}\OperatorTok{=}\NormalTok{inf,     }
\NormalTok{    n\_initial}\OperatorTok{=}\DecValTok{10}\NormalTok{,    }
\NormalTok{    seed}\OperatorTok{=}\DecValTok{42}\NormalTok{,}
\NormalTok{    surrogate}\OperatorTok{=}\NormalTok{S\_GP,}
\NormalTok{    max\_surrogate\_points}\OperatorTok{=}\DecValTok{30}\NormalTok{,}
\NormalTok{)}
\NormalTok{so\_optimizer.optimize()}
\end{Highlighting}
\end{Shaded}

\begin{verbatim}
 message: Optimization terminated: maximum evaluations (50) reached
                   Current function value: 0.209939
                   Iterations: 42
                   Function evaluations: 50
 success: True
     fun: 0.20993878639170102
       x: [128.0 5.0 'LeakyReLU' 11.214919451291014 'Adam' 1024.0]
       X: [[16.0 3.0 ... 'AdamW' 256.0]
           [16.0 4.0 ... 'SGD' 128.0]
           ...
           [128.0 5.0 ... 'Adam' 1024.0]
           [128.0 5.0 ... 'Adam' 1024.0]]
     nit: 42
    nfev: 50
       y: [ 2.787e+04  2.839e+04 ...  3.926e+00  5.453e+00]
\end{verbatim}

\subsubsection{Results}\label{results}

We can inspect the best result found. \texttt{spotoptim} provides the
method \texttt{print\_results} to print the results in a table format.

\begin{Shaded}
\begin{Highlighting}[]
\NormalTok{so\_optimizer.print\_results(precision}\OperatorTok{=}\DecValTok{2}\NormalTok{)}
\end{Highlighting}
\end{Shaded}

\begin{verbatim}
|              name |   type |   default |   lower |   upper |     tuned |   transform |
|-------------------|--------|-----------|---------|---------|-----------|-------------|
|                l1 |    int |        72 |   16.00 |  128.00 |       128 |   log(x, 2) |
| num_hidden_layers |    int |         3 |    1.00 |    5.00 |         5 |           - |
|        activation | factor |   Sigmoid |       - |       - | LeakyReLU |           - |
|                lr |  float |     50.00 |    0.00 |  100.00 |     11.21 |         log |
|         optimizer | factor |   RMSprop |       - |       - |      Adam |           - |
|            epochs |    int |       528 |   32.00 | 1024.00 |      1024 |   log(x, 2) |
\end{verbatim}

\begin{Shaded}
\begin{Highlighting}[]
\BuiltInTok{print}\NormalTok{(}\SpecialStringTok{f"Best MSE Loss: }\SpecialCharTok{\{}\NormalTok{so\_optimizer}\SpecialCharTok{.}\NormalTok{best\_y\_}\SpecialCharTok{:.4f\}}\SpecialStringTok{"}\NormalTok{)}
\BuiltInTok{print}\NormalTok{(}\StringTok{"Best Configuration:"}\NormalTok{)}
\NormalTok{best\_config }\OperatorTok{=}\NormalTok{ so\_optimizer.get\_best\_hyperparameters()}

\ControlFlowTok{if}\NormalTok{ best\_config:}
    \ControlFlowTok{for}\NormalTok{ k, v }\KeywordTok{in}\NormalTok{ best\_config.items():}
        \BuiltInTok{print}\NormalTok{(}\SpecialStringTok{f"  }\SpecialCharTok{\{}\NormalTok{k}\SpecialCharTok{\}}\SpecialStringTok{: }\SpecialCharTok{\{}\NormalTok{v}\SpecialCharTok{\}}\SpecialStringTok{"}\NormalTok{)}
\ControlFlowTok{else}\NormalTok{:}
    \BuiltInTok{print}\NormalTok{(}\StringTok{"No best configuration found."}\NormalTok{)}
\end{Highlighting}
\end{Shaded}

\begin{verbatim}
Best MSE Loss: 0.2099
Best Configuration:
  l1: 128
  num_hidden_layers: 5
  activation: LeakyReLU
  lr: 11.214919451291014
  optimizer: Adam
  epochs: 1024
\end{verbatim}

\subsubsection{Visualization}\label{visualization}

\begin{figure}[H]

\centering{

\pandocbounded{\includegraphics[keepaspectratio]{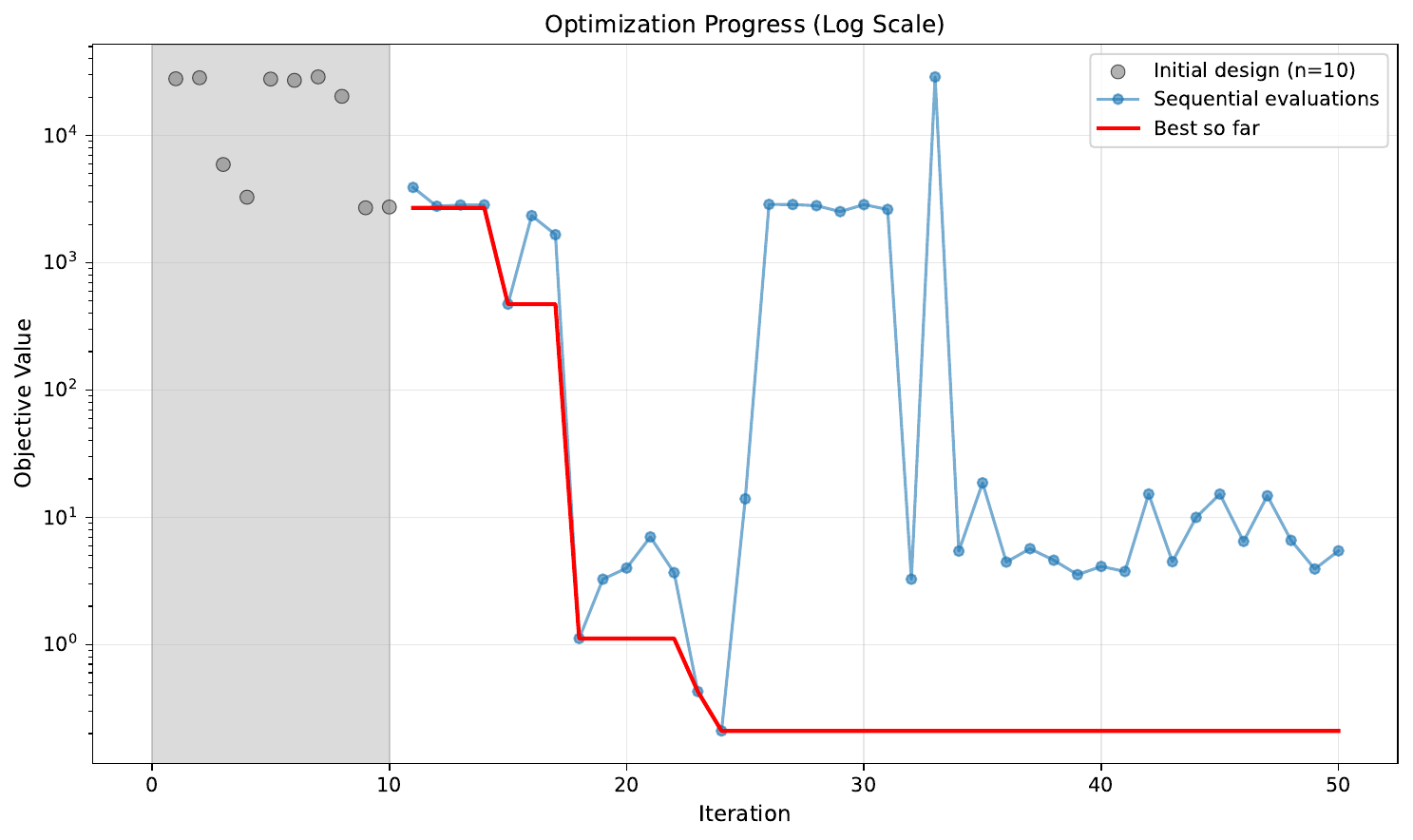}}

}

\caption{\label{fig-plot-progress-1}Progress of the HPT process using
the single-objective approach.}

\end{figure}

Plotting the optimization progress. Figure~\ref{fig-plot-progress-1}
shows the hyperparameter tuning process. The loss is plotted versus the
function evaluations. The x-axis shows the number of function
evaluations, and the y-axis shows the loss. The loss evaluated by the
objective function is plotted in blue, the best loss found is plotted in
red. The y-axis is set to a logarithmic scale for better visualization.

\begin{Shaded}
\begin{Highlighting}[]
\NormalTok{so\_optimizer.plot\_progress(log\_y}\OperatorTok{=}\VariableTok{True}\NormalTok{)}
\end{Highlighting}
\end{Shaded}

\begin{tcolorbox}[enhanced jigsaw, colbacktitle=quarto-callout-note-color!10!white, opacityback=0, title=\textcolor{quarto-callout-note-color}{\faInfo}\hspace{0.5em}{Results from the Single-Objective Approach}, opacitybacktitle=0.6, arc=.35mm, left=2mm, bottomrule=.15mm, breakable, toprule=.15mm, coltitle=black, colframe=quarto-callout-note-color-frame, colback=white, rightrule=.15mm, bottomtitle=1mm, toptitle=1mm, titlerule=0mm, leftrule=.75mm]

\begin{itemize}
\tightlist
\item
  The single-objective approach resulted in a validation loss of
  \texttt{6.694} and \texttt{1024} (\(=2^{10}\)) epochs.
\end{itemize}

\end{tcolorbox}

\subsection{Weighted Multi-Objective Function}\label{sec-hpt-weighted}

The second approach is to use a weighted mean of both objective
functions. This is done by setting the \texttt{fun\_mo2so} argument to a
custom function that computes the weighted mean of both objective
functions. The weights can be adjusted to give more importance to one
objective function over the other. Here, we define the function
\texttt{aggregate} that computes the weighted mean of both objective
functions. The first objective function is weighted with 2 and the
second objective function is weighted with 0.1. Now we consider two
objectives: \texttt{MSE} (to minimize error) and \texttt{epochs} (to
minimize computational cost), which will be aggregated using a weighted
sum.

\subsubsection{Update Experiment
Configuration}\label{update-experiment-configuration}

We update the experiment to return both metrics, i.e, \texttt{mse} and
\texttt{epochs}. This is done by setting the \texttt{metrics} argument
in the \texttt{ExperimentControl} constructor to a list of metric names.
Furthermore, the \texttt{TorchObjective} instance \texttt{objective} is
re-initialized to ensure that the changes are picked up. As we can see
from the following output, the same hyperparameters are extracted from
the input vector, but two metrics are returned. The first metric is the
\texttt{mse} and the second metric is the \texttt{epochs}

\phantomsection\label{update-experiment}
\begin{Shaded}
\begin{Highlighting}[]
\NormalTok{exp.metrics }\OperatorTok{=}\NormalTok{ [}\StringTok{"mse"}\NormalTok{, }\StringTok{"epochs"}\NormalTok{]}
\NormalTok{objective }\OperatorTok{=}\NormalTok{ TorchObjective(experiment}\OperatorTok{=}\NormalTok{exp, seed}\OperatorTok{=}\DecValTok{42}\NormalTok{)}
\NormalTok{X\_demo }\OperatorTok{=}\NormalTok{ np.array([[}\FloatTok{0.1}\NormalTok{, }\FloatTok{0.2}\NormalTok{, }\FloatTok{0.3}\NormalTok{, }\FloatTok{0.4}\NormalTok{, }\FloatTok{0.5}\NormalTok{, }\FloatTok{9.0}\NormalTok{]])}
\NormalTok{hps }\OperatorTok{=}\NormalTok{ objective.\_get\_hyperparameters(X}\OperatorTok{=}\NormalTok{X\_demo)}
\NormalTok{pprint.pp(}\SpecialStringTok{f"Hyperparameters: }\SpecialCharTok{\{}\NormalTok{hps}\SpecialCharTok{\}}\SpecialStringTok{"}\NormalTok{)}
\NormalTok{y\_eval }\OperatorTok{=}\NormalTok{ objective(X\_demo)}
\BuiltInTok{print}\NormalTok{(}\SpecialStringTok{f"y\_eval: }\SpecialCharTok{\{}\NormalTok{y\_eval}\SpecialCharTok{\}}\SpecialStringTok{"}\NormalTok{)}
\end{Highlighting}
\end{Shaded}

\begin{verbatim}
("Hyperparameters: {'l1': 0, 'num_hidden_layers': 0, 'activation': 'ReLU', "
 "'lr': np.float64(0.4), 'optimizer': 'Adam', 'epochs': 9}")
y_eval: [[2.89048407e+04 9.00000000e+00]]
\end{verbatim}

\subsubsection{Definition of the Aggregation
Function}\label{definition-of-the-aggregation-function}

We define a function to aggregate the multi-objective output \texttt{y}
with shape \texttt{(n,\ 2)} into a single scalar value for the
optimizer. \texttt{MSE} is weighted by 2.0, \texttt{epochs} by 0.1. We
can use the \texttt{fun\_mo2so} argument to pass the aggregation
function to \texttt{SpotOptim}. It is used to aggregate the
multi-objective output \texttt{y} into a single scalar value for the
optimizer.

\begin{Shaded}
\begin{Highlighting}[]
\KeywordTok{def}\NormalTok{ aggregate(y):}
    \CommentTok{\# y is (n\_samples, 2) {-}\textgreater{} [MSE, Epochs]}
\NormalTok{    weights }\OperatorTok{=}\NormalTok{ np.array([}\FloatTok{2.0}\NormalTok{, }\FloatTok{0.1}\NormalTok{])}
    \ControlFlowTok{return}\NormalTok{ np.}\BuiltInTok{sum}\NormalTok{(y }\OperatorTok{*}\NormalTok{ weights, axis}\OperatorTok{=}\DecValTok{1}\NormalTok{)}

\NormalTok{optimizer\_mo\_weighted }\OperatorTok{=}\NormalTok{ SpotOptim(}
\NormalTok{    fun}\OperatorTok{=}\NormalTok{objective, }\CommentTok{\# Now returns (n, 2)}
\NormalTok{    max\_iter}\OperatorTok{=}\DecValTok{50}\NormalTok{,}
\NormalTok{    max\_time}\OperatorTok{=}\NormalTok{inf,}
\NormalTok{    n\_initial}\OperatorTok{=}\DecValTok{15}\NormalTok{,}
\NormalTok{    seed}\OperatorTok{=}\DecValTok{42}\NormalTok{,}
\NormalTok{    fun\_mo2so}\OperatorTok{=}\NormalTok{aggregate,}
\NormalTok{    surrogate}\OperatorTok{=}\NormalTok{S\_GP,}
\NormalTok{    max\_surrogate\_points}\OperatorTok{=}\DecValTok{30}\NormalTok{,}
\NormalTok{)}
\NormalTok{optimizer\_mo\_weighted.optimize()}
\end{Highlighting}
\end{Shaded}

\phantomsection\label{optimizer-weighted}
\begin{verbatim}
 message: Optimization terminated: maximum evaluations (50) reached
                   Current function value: 151.721190
                   Iterations: 36
                   Function evaluations: 50
 success: True
     fun: 151.72118955339704
       x: [128.0 3.0 'ReLU' 3.2325995355674486 'RMSprop' 512.0]
       X: [[128.0 4.0 ... 'AdamW' 128.0]
           [32.0 2.0 ... 'SGD' 64.0]
           ...
           [64.0 3.0 ... 'RMSprop' 512.0]
           [64.0 3.0 ... 'RMSprop' 512.0]]
     nit: 36
    nfev: 50
       y: [ 5.792e+04  1.182e+04 ...  5.337e+02  1.336e+03]
\end{verbatim}

Note that \texttt{y} in the output above is the aggregated value of
\texttt{loss} and \texttt{epochs} computed by the \texttt{aggregate}
function.

\subsubsection{Results}\label{results-1}

The best objective function value of each single objective (from the
multi-objective optimization) can be retrieved as follows:

\begin{Shaded}
\begin{Highlighting}[]
\NormalTok{best\_idx }\OperatorTok{=}\NormalTok{ np.argmin(optimizer\_mo\_weighted.y\_)}
\NormalTok{best\_objectives }\OperatorTok{=}\NormalTok{ optimizer\_mo\_weighted.y\_mo[best\_idx]}
\BuiltInTok{print}\NormalTok{(}\SpecialStringTok{f"Best Objectives: MSE=}\SpecialCharTok{\{}\NormalTok{best\_objectives[}\DecValTok{0}\NormalTok{]}\SpecialCharTok{:.2f\}}\SpecialStringTok{, epochs=}\SpecialCharTok{\{}\NormalTok{best\_objectives[}\DecValTok{1}\NormalTok{]}\SpecialCharTok{:.0f\}}\SpecialStringTok{"}\NormalTok{)}
\end{Highlighting}
\end{Shaded}

\begin{verbatim}
Best Objectives: MSE=50.26, epochs=512
\end{verbatim}

\begin{Shaded}
\begin{Highlighting}[]
\BuiltInTok{print}\NormalTok{(}\SpecialStringTok{f"Best Weighted Score: }\SpecialCharTok{\{}\NormalTok{optimizer\_mo\_weighted}\SpecialCharTok{.}\NormalTok{best\_y\_}\SpecialCharTok{:.4f\}}\SpecialStringTok{"}\NormalTok{)}
\CommentTok{\# Best generic result}
\NormalTok{optimizer\_mo\_weighted.print\_results(precision}\OperatorTok{=}\DecValTok{2}\NormalTok{)}
\end{Highlighting}
\end{Shaded}

\begin{verbatim}
Best Weighted Score: 151.7212
|              name |   type |   default |   lower |   upper |   tuned |   transform |
|-------------------|--------|-----------|---------|---------|---------|-------------|
|                l1 |    int |        72 |   16.00 |  128.00 |     128 |   log(x, 2) |
| num_hidden_layers |    int |         3 |    1.00 |    5.00 |       3 |           - |
|        activation | factor |   Sigmoid |       - |       - |    ReLU |           - |
|                lr |  float |     50.00 |    0.00 |  100.00 |    3.23 |         log |
|         optimizer | factor |   RMSprop |       - |       - | RMSprop |           - |
|            epochs |    int |       528 |   32.00 | 1024.00 |     512 |   log(x, 2) |
\end{verbatim}

\subsubsection{Visualization}\label{visualization-1}

\begin{Shaded}
\begin{Highlighting}[]
\NormalTok{optimizer\_mo\_weighted.plot\_progress(log\_y}\OperatorTok{=}\VariableTok{True}\NormalTok{, mo}\OperatorTok{=}\VariableTok{True}\NormalTok{)}
\end{Highlighting}
\end{Shaded}

\begin{figure}[H]

\centering{

\pandocbounded{\includegraphics[keepaspectratio]{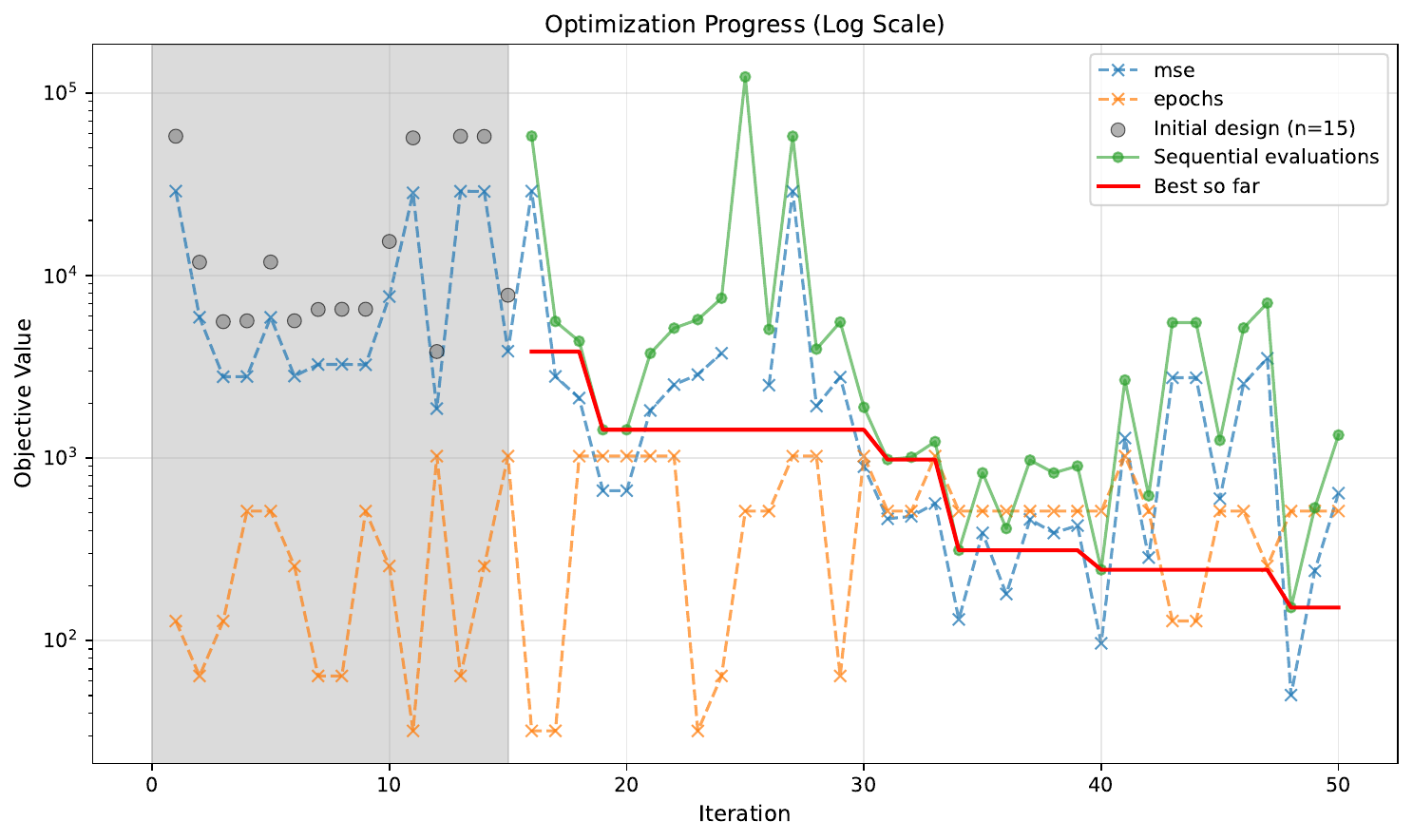}}

}

\caption{\label{fig-plot-progress-mo}Progress of the HPT process using
the weighted multi-objective approach.}

\end{figure}%

\begin{tcolorbox}[enhanced jigsaw, colbacktitle=quarto-callout-note-color!10!white, opacityback=0, title=\textcolor{quarto-callout-note-color}{\faInfo}\hspace{0.5em}{Interpretation}, opacitybacktitle=0.6, arc=.35mm, left=2mm, bottomrule=.15mm, breakable, toprule=.15mm, coltitle=black, colframe=quarto-callout-note-color-frame, colback=white, rightrule=.15mm, bottomtitle=1mm, toptitle=1mm, titlerule=0mm, leftrule=.75mm]

The plot shows how the optimizer explores the trade-off. By including
epochs in the objective, the optimizer may favor configurations that
achieve reasonable loss with fewer training epochs, or it might find
that more epochs are necessary for lower loss, depending on the weights.

\end{tcolorbox}

\subsection{Multi-Objective Hyperparameter Tuning With
Desirability}\label{sec-hpt-desirability}

Since the weights used in Section~\ref{sec-hpt-weighted} do not have a
physical meaning and are not easy to interpret, we use desirability
functions to map each objective to a {[}0, 1{]} scale and combine them.
Desirabilities are easier to interpret than weights.

\subsubsection{Setting Up the Desirability
Function}\label{setting-up-the-desirability-function}

The third approach is to use a desirability function to combine the two
objective functions into a single objective function. The desirability
function is used to maximize the desirability of the two objective
functions. The desirability function is defined by setting the
\texttt{fun\_mo2so} argument to a custom function that computes the
desirability of both objective functions. The desirability function is
defined in the following code segment. We use \texttt{DMin} to define
that we want to minimize both Loss and Epochs within specific ranges.
\texttt{DMin} transforms the incoming values to a {[}0, 1{]} scale,
where 1 is the best and 0 is the worst. The target is to minimize both
objective functions, so the low value is good (desirability 1), the high
value is bad (desirability 0). From pre-experiments we know that the
loss is around 2500 and unacceptable above 6000. \texttt{SpotOptim}
minimizes the objective, so we return \texttt{(1\ -\ desirability)}.

\begin{Shaded}
\begin{Highlighting}[]
\KeywordTok{def}\NormalTok{ desirability(y):}
    \CommentTok{\# y is (n\_samples, 2) {-}\textgreater{} [MSE, Epochs]}
    \CommentTok{\# Loss: desirable below 2000, unacceptable above 6000? }
    \CommentTok{\# Adjusted ranges based on typical Diabetes dataset values \textasciitilde{}3000{-}6000}
\NormalTok{    lossD }\OperatorTok{=}\NormalTok{ DMin(}\DecValTok{6}\NormalTok{, }\DecValTok{6000}\NormalTok{, scale}\OperatorTok{=}\DecValTok{2}\NormalTok{) }
    \CommentTok{\# Epochs: desirable 16, unacceptable 1024}
\NormalTok{    epochsD }\OperatorTok{=}\NormalTok{ DMin(}\DecValTok{32}\NormalTok{, }\DecValTok{1024}\NormalTok{, scale}\OperatorTok{=}\FloatTok{.2}\NormalTok{)}
\NormalTok{    overallD }\OperatorTok{=}\NormalTok{ DOverall(lossD, epochsD)}
    \CommentTok{\# y is passed as (n, 2)}
    \ControlFlowTok{return} \FloatTok{1.0} \OperatorTok{{-}}\NormalTok{ overallD.predict(y, }\BuiltInTok{all}\OperatorTok{=}\VariableTok{False}\NormalTok{)}
\end{Highlighting}
\end{Shaded}

\begin{figure}

\centering{

\pandocbounded{\includegraphics[keepaspectratio]{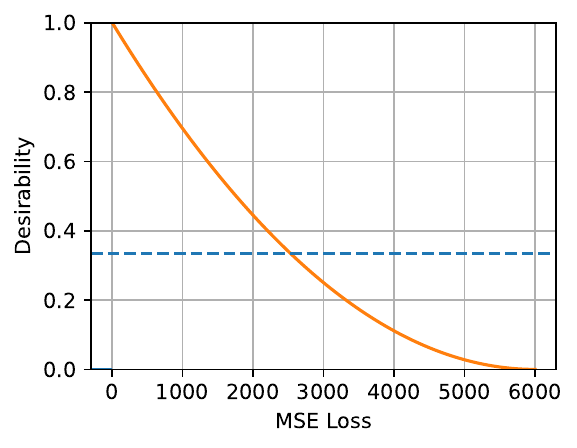}}

}

\caption{\label{fig-des-loss-1}Desirability for MSE Loss}

\end{figure}

\begin{figure}

\centering{

\pandocbounded{\includegraphics[keepaspectratio]{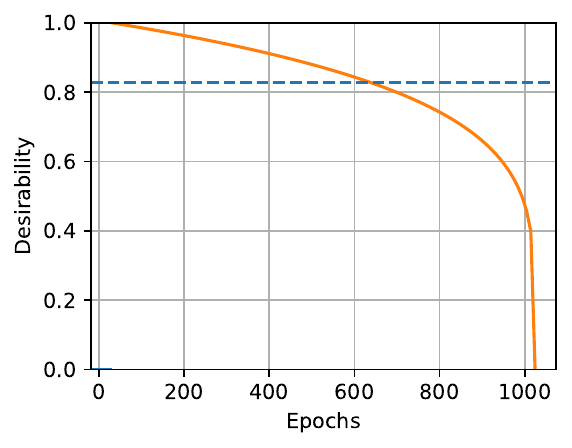}}

}

\caption{\label{fig-des-epochs-1}Desirability for Epochs}

\end{figure}

The desirability functions are visualized in Figure~\ref{fig-des-loss-1}
and Figure~\ref{fig-des-epochs-1}.

We have chosen simple desirability functions based on \texttt{DMin} for
the \texttt{validation\ loss} and \texttt{number\ of\ epochs}. The usage
of these functions might result in large plateaus where the optimizer
does not find any improvement. Therefore, we will explore more
sophisticated desirability functions in the future. These can easily be
implemented with the approach shown in
\{Section~\ref{sec-nonstandard-desirabilities}\}.

\subsubsection{Run Optimization}\label{run-optimization-1}

\begin{Shaded}
\begin{Highlighting}[]
\NormalTok{des\_optimizer }\OperatorTok{=}\NormalTok{ SpotOptim(}
\NormalTok{    fun}\OperatorTok{=}\NormalTok{objective,}
\NormalTok{    max\_iter}\OperatorTok{=}\DecValTok{50}\NormalTok{,     }\CommentTok{\# More iterations to see Pareto front better}
\NormalTok{    max\_time}\OperatorTok{=}\NormalTok{inf,      }\CommentTok{\# minutes}
\NormalTok{    n\_initial}\OperatorTok{=}\DecValTok{10}\NormalTok{,}
\NormalTok{    seed}\OperatorTok{=}\DecValTok{42}\NormalTok{,}
\NormalTok{    fun\_mo2so}\OperatorTok{=}\NormalTok{desirability,}
\NormalTok{    surrogate}\OperatorTok{=}\NormalTok{S\_GP,}
\NormalTok{    max\_surrogate\_points}\OperatorTok{=}\DecValTok{30}\NormalTok{,}
\NormalTok{)}
\NormalTok{des\_optimizer.optimize()}
\end{Highlighting}
\end{Shaded}

\phantomsection\label{des-optimizer}
\begin{verbatim}
 message: Optimization terminated: maximum evaluations (50) reached
                   Current function value: 0.041187
                   Iterations: 40
                   Function evaluations: 50
 success: True
     fun: 0.041187076557585356
       x: [128.0 5.0 'ELU' 19.14206637349006 'Adam' 256.0]
       X: [[16.0 3.0 ... 'AdamW' 256.0]
           [16.0 4.0 ... 'SGD' 128.0]
           ...
           [128.0 5.0 ... 'Adam' 256.0]
           [128.0 5.0 ... 'Adam' 256.0]]
     nit: 40
    nfev: 50
       y: [ 1.000e+00  1.000e+00 ...  5.115e-02  5.471e-02]
\end{verbatim}

\subsubsection{Results}\label{results-2}

\texttt{spotoptim}'s \texttt{print\_results} function can be used to
print the results in a table.

\phantomsection\label{results-desirability}
\begin{Shaded}
\begin{Highlighting}[]
\NormalTok{\_ }\OperatorTok{=}\NormalTok{ des\_optimizer.print\_results(precision}\OperatorTok{=}\DecValTok{2}\NormalTok{)}
\end{Highlighting}
\end{Shaded}

\begin{verbatim}
|              name |   type |   default |   lower |   upper |   tuned |   transform |
|-------------------|--------|-----------|---------|---------|---------|-------------|
|                l1 |    int |        72 |   16.00 |  128.00 |     128 |   log(x, 2) |
| num_hidden_layers |    int |         3 |    1.00 |    5.00 |       5 |           - |
|        activation | factor |   Sigmoid |       - |       - |     ELU |           - |
|                lr |  float |     50.00 |    0.00 |  100.00 |   19.14 |         log |
|         optimizer | factor |   RMSprop |       - |       - |    Adam |           - |
|            epochs |    int |       528 |   32.00 | 1024.00 |     256 |   log(x, 2) |
\end{verbatim}

\begin{Shaded}
\begin{Highlighting}[]
\NormalTok{des\_optimizer.plot\_progress(log\_y}\OperatorTok{=}\VariableTok{True}\NormalTok{, mo}\OperatorTok{=}\VariableTok{True}\NormalTok{)}
\end{Highlighting}
\end{Shaded}

\begin{figure}[H]

\centering{

\pandocbounded{\includegraphics[keepaspectratio]{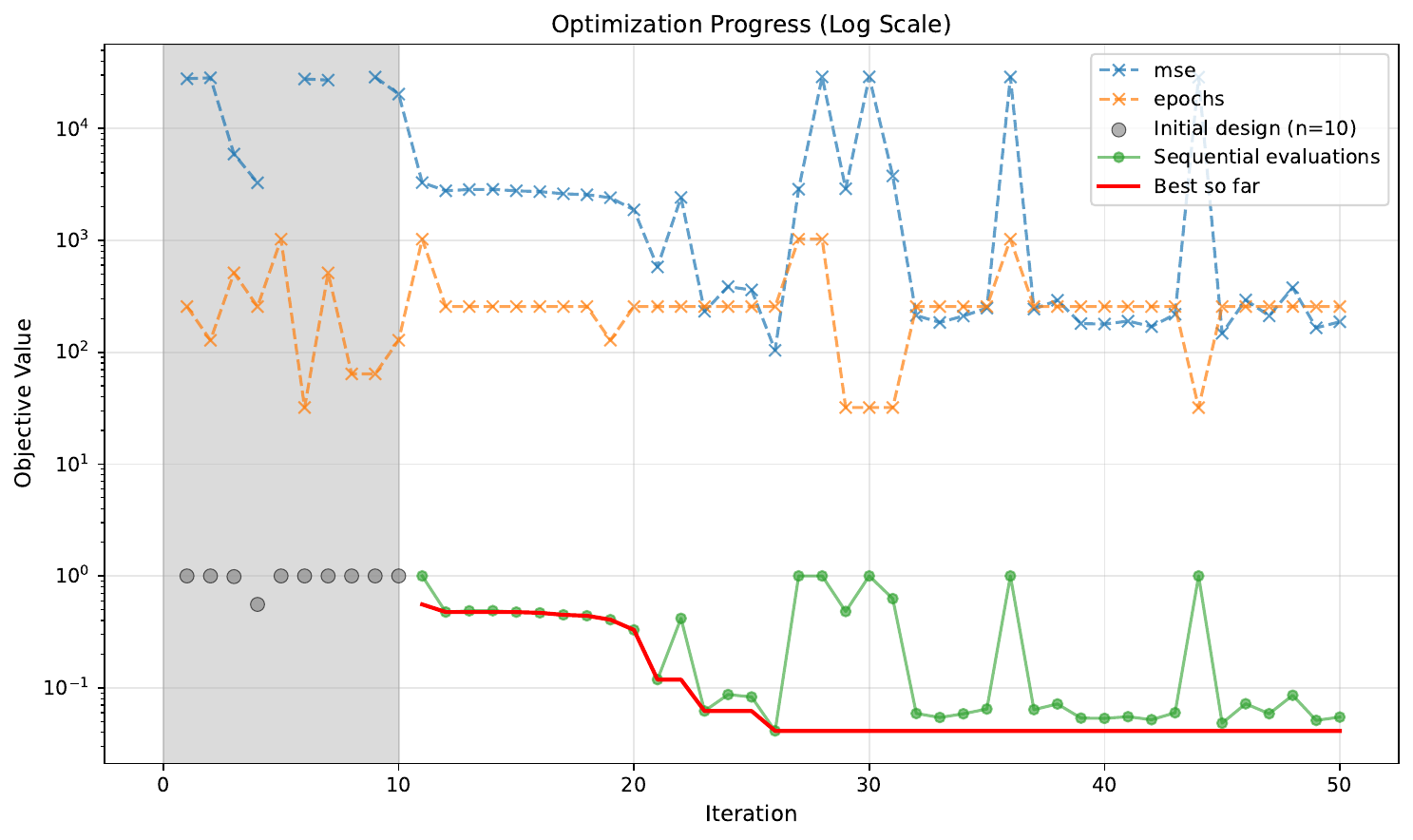}}

}

\caption{\label{fig-progress-desirability}Progress of the HPT process
using desirability functions.}

\end{figure}%

The best desirability score can be retrieved as follows:

\begin{Shaded}
\begin{Highlighting}[]
\BuiltInTok{print}\NormalTok{(}\SpecialStringTok{f"Best Desirability Score: }\SpecialCharTok{\{}\FloatTok{1.0} \OperatorTok{{-}}\NormalTok{ des\_optimizer}\SpecialCharTok{.}\NormalTok{best\_y\_}\SpecialCharTok{:.4f\}}\SpecialStringTok{"}\NormalTok{)}
\end{Highlighting}
\end{Shaded}

\begin{verbatim}
Best Desirability Score: 0.9588
\end{verbatim}

The best objective function value of each single objective (from the
multi-objective optimization) can be retrieved as follows:

\begin{Shaded}
\begin{Highlighting}[]
\NormalTok{best\_idx }\OperatorTok{=}\NormalTok{ np.argmin(des\_optimizer.y\_)}
\NormalTok{best\_objectives }\OperatorTok{=}\NormalTok{ des\_optimizer.y\_mo[best\_idx]}
\BuiltInTok{print}\NormalTok{(}\SpecialStringTok{f"Best Objectives: MSE=}\SpecialCharTok{\{}\NormalTok{best\_objectives[}\DecValTok{0}\NormalTok{]}\SpecialCharTok{:.2f\}}\SpecialStringTok{, epochs=}\SpecialCharTok{\{}\NormalTok{best\_objectives[}\DecValTok{1}\NormalTok{]}\SpecialCharTok{:.0f\}}\SpecialStringTok{"}\NormalTok{)}
\end{Highlighting}
\end{Shaded}

\begin{verbatim}
Best Objectives: MSE=103.89, epochs=256
\end{verbatim}

\begin{figure}

\centering{

\pandocbounded{\includegraphics[keepaspectratio]{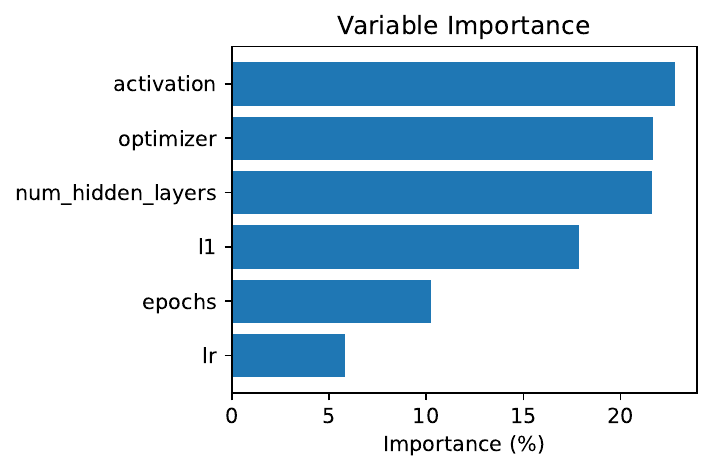}}

}

\caption{\label{fig-spot-importance-des}The importance of the parameters
in the optimization process}

\end{figure}

\begin{figure}

\centering{

\pandocbounded{\includegraphics[keepaspectratio]{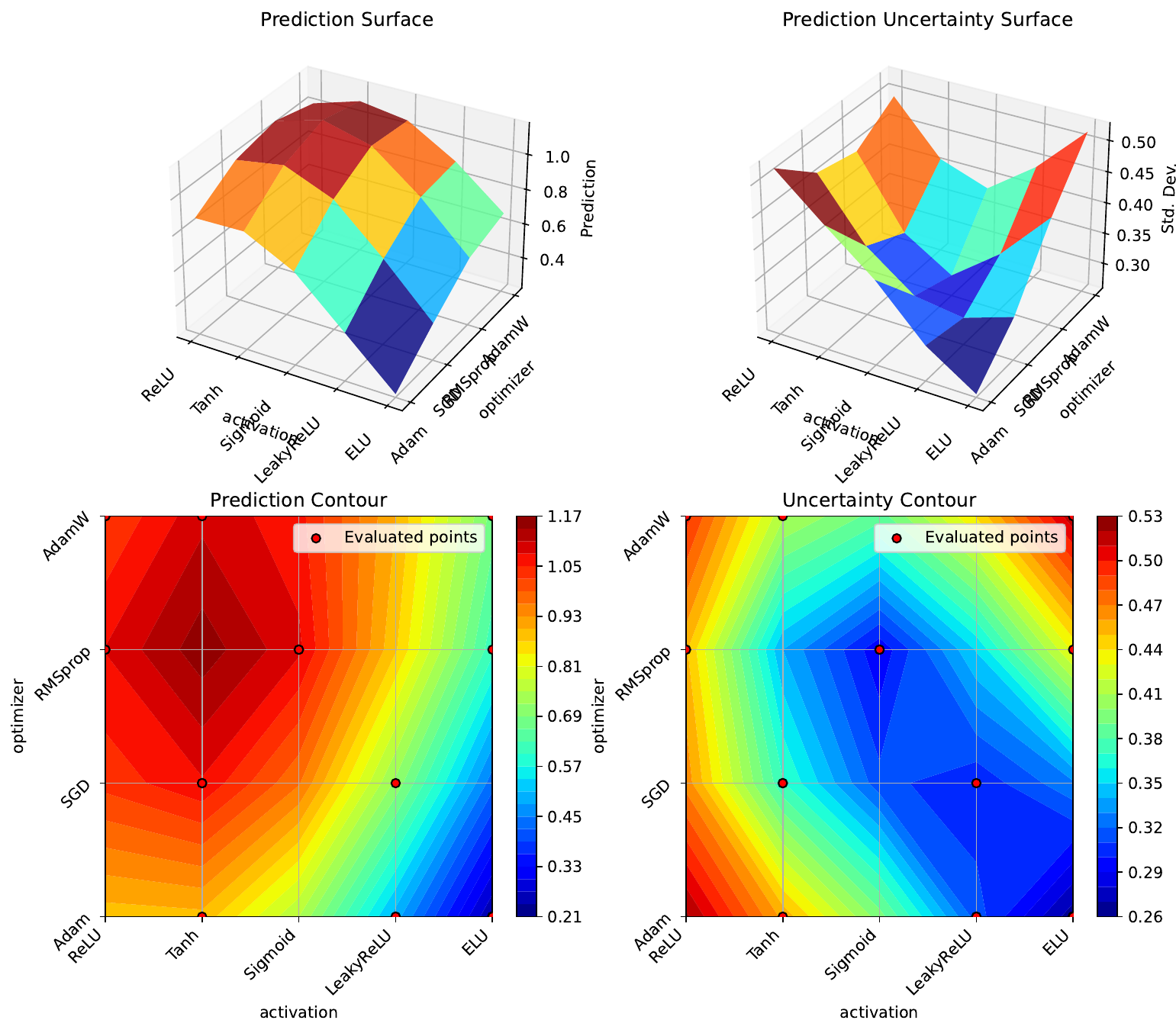}}

}

\caption{\label{fig-spot-importance-contour-des}HPT with desirability
functions: The contour plots for the two most important parameters.
Note: SpotOptim is minimizing the objective function, so we return (1 -
desirability).}

\end{figure}

The surrogate model can be used to visualize the importance of the
parameters in the optimization process as shown in
Figure~\ref{fig-spot-importance-des} and
Figure~\ref{fig-spot-importance-contour-des}.

\begin{verbatim}
Plotting surrogate contours for top 2 most important parameters:
  activation: importance = 22.78% (type: factor)
  optimizer: importance = 21.68% (type: factor)

Generating 1 surrogate plots...
  Plotting activation vs optimizer
\end{verbatim}

\subsubsection{Pareto Front}\label{pareto-front}

\begin{figure}

\centering{

\pandocbounded{\includegraphics[keepaspectratio]{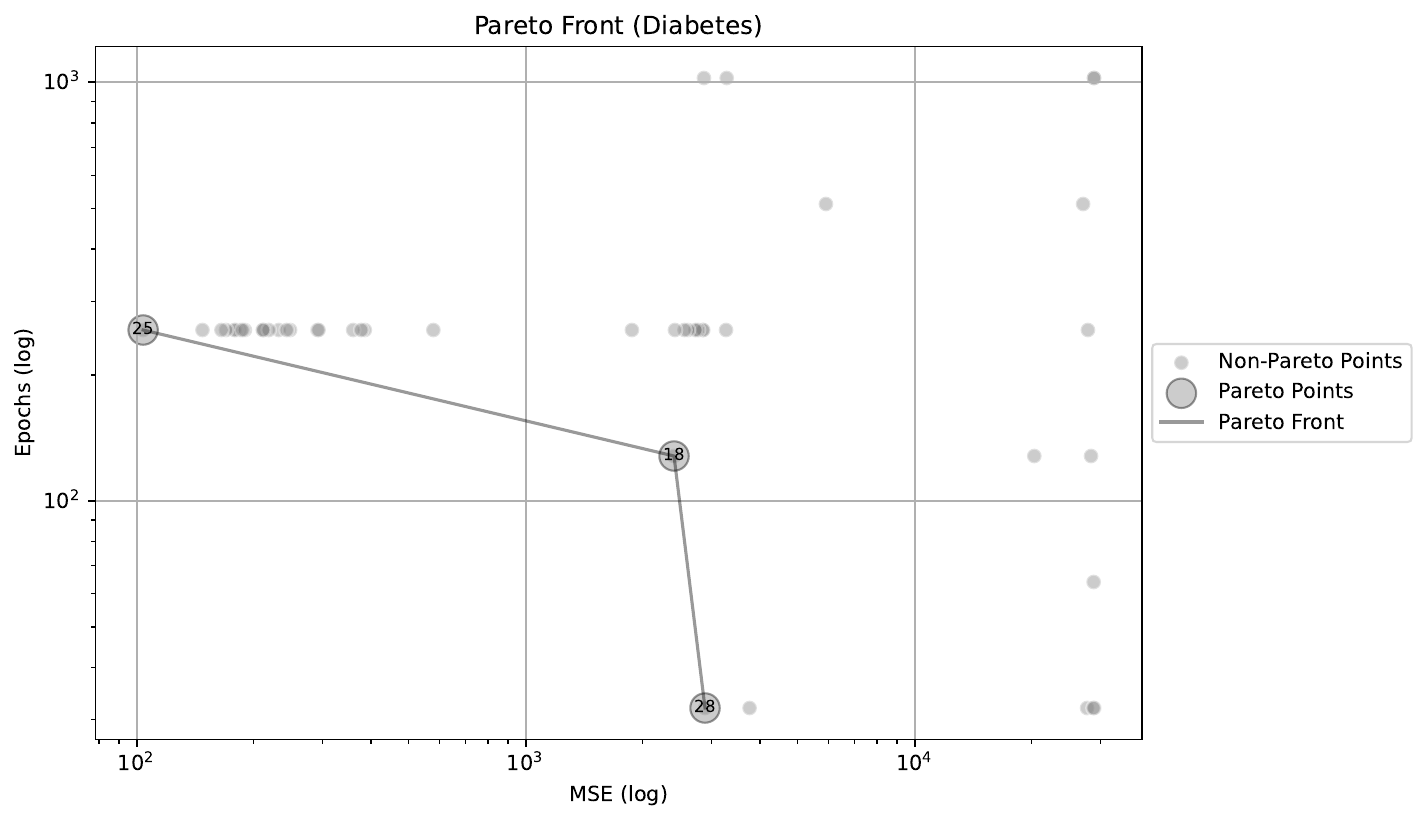}}

}

\caption{\label{fig-pareto}Pareto Front: Loss vs Epochs. Both objectives
should be minimized.}

\end{figure}

We visualize the trade-off using \texttt{spotoptim}'s \texttt{plot\_mo}
method in Figure~\ref{fig-pareto}.

\begin{tcolorbox}[enhanced jigsaw, colbacktitle=quarto-callout-note-color!10!white, opacityback=0, title=\textcolor{quarto-callout-note-color}{\faInfo}\hspace{0.5em}{Interpretation}, opacitybacktitle=0.6, arc=.35mm, left=2mm, bottomrule=.15mm, breakable, toprule=.15mm, coltitle=black, colframe=quarto-callout-note-color-frame, colback=white, rightrule=.15mm, bottomtitle=1mm, toptitle=1mm, titlerule=0mm, leftrule=.75mm]

The Pareto front highlights the non-dominated solutions. You can see how
lower MSE values generally require more epochs (or specific
hyperparameters that coincide with higher training cost), while very few
epochs yield higher MSE.

\end{tcolorbox}

\section{Using Space-Fillingness as an
Objective}\label{sec-space-fillingness}

Designed experiments are not always possible in real-world applications.
Practitioners are often forced to use existing data or to collect data
in an ad-hoc manner. In such cases, it is important to ensure that the
data is representative of the entire space of possible inputs. The
``space-fillingness'' of the experimental design is a useful metric for
this purpose. One well-established method to measure space-fillingness
is the Morris-Mitchell criterion (Morris and Mitchell 1995).

\subsection{Morris-Mitchell Criterion}\label{morris-mitchell-criterion}

The Morris-Mitchell criterion is a widely used metric for evaluating the
space-filling properties of a sampling plan (also known as a design of
experiments). It combines the concept of maximizing the minimum distance
between points with minimizing the number of pairs of points separated
by that distance. The standard Morris-Mitchell criterion, \(\Phi_q\), is
defined as:

\[
\Phi_q = \left( \sum_{i=1}^{m} J_i d_i^{-q} \right)^{1/q}
\]

where \(d_1 < d_2 < \dots < d_m\) are the distinct distances between all
pairs of points in the sampling plan, \(J_i\) is the number of pairs of
points separated by the distance \(d_i\), and \(q\) is a positive
integer (e.g., \(q=2, 5, 10, \dots\)). As \(q \to \infty\), minimizing
\(\Phi_q\) is equivalent to maximizing the minimum distance \(d_1\), and
then minimizing the number of pairs \(J_1\) at that distance, and so on.
One limitation of the standard \(\Phi_q\) is that its value depends on
the number of points in the sampling plan. This makes it difficult to
compare the quality of sampling plans with different sample sizes. To
address this, we can use a size-invariant or intensive version of the
criterion. This is achieved by normalizing the sum by the total number
of pairs, \(M = \binom{n}{2} = \frac{n(n-1)}{2}\). The intensive
criterion is defined as:

\[
\Phi_{q, \text{intensive}} = \left( \frac{1}{M} \sum_{i=1}^{m} J_i d_i^{-q} \right)^{1/q}
\]

\begin{figure}

\centering{

\pandocbounded{\includegraphics[keepaspectratio]{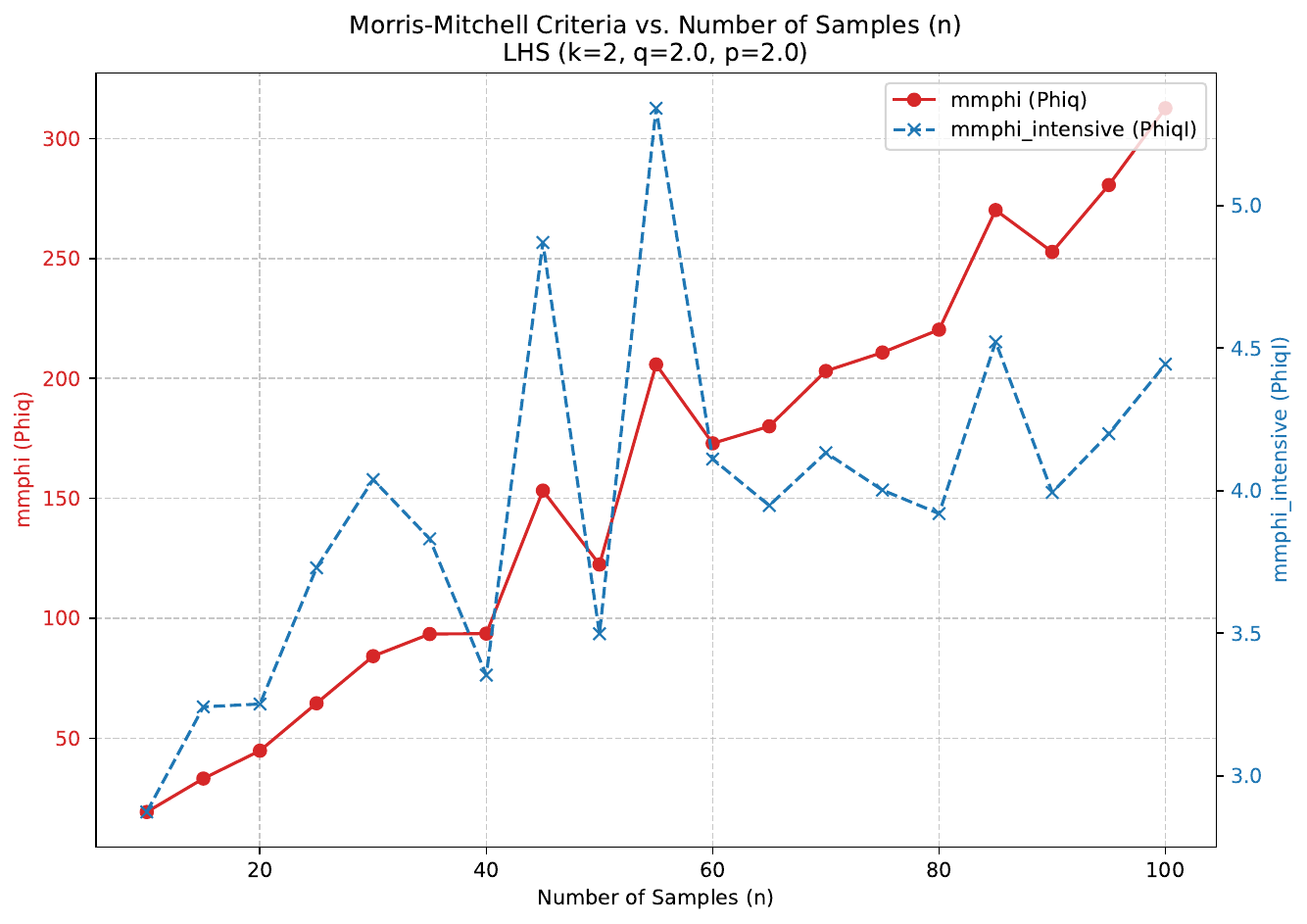}}

}

\caption{\label{fig-mm-comparison}Comparison of the Morris-Mitchell and
the intesnsified Morris-Mitchell criterion. Note: different axes scaling
was used for the two criteria.}

\end{figure}

This normalization allows for a fairer comparison between designs of
different sizes. A lower value still indicates a better space-filling
property. The \texttt{spotoptim} library provides the
\texttt{mmphi\_intensive} function to calculate this metric.
Figure~\ref{fig-mm-comparison} shows a comparison of the Morris-Mitchell
and the intesnsified Morris-Mitchell criterion.

\begin{example}[Usage of
\texttt{mmphi\_intensive}]\protect\hypertarget{exm-mmphi-intensive}{}\label{exm-mmphi-intensive}

We create a simple 3-point sampling plan in 2D and calculate the
intensive space-fillingness metric with q=2, using Euclidean distances
(p=2).

\phantomsection\label{mmphi-intensive-example}
\begin{Shaded}
\begin{Highlighting}[]
\NormalTok{X3 }\OperatorTok{=}\NormalTok{ np.array([}
\NormalTok{    [}\FloatTok{0.0}\NormalTok{, }\FloatTok{0.0}\NormalTok{],}
\NormalTok{    [}\FloatTok{0.5}\NormalTok{, }\FloatTok{0.5}\NormalTok{],}
\NormalTok{    [}\FloatTok{1.0}\NormalTok{, }\FloatTok{1.0}\NormalTok{]}
\NormalTok{])}
\NormalTok{quality, J, d }\OperatorTok{=}\NormalTok{ mmphi\_intensive(X3, q}\OperatorTok{=}\DecValTok{2}\NormalTok{, p}\OperatorTok{=}\DecValTok{2}\NormalTok{)}
\end{Highlighting}
\end{Shaded}

\phantomsection\label{mmphi-intensive-example-output}
\begin{verbatim}
Quality (Phi_q_intensive): 1.224744871391589
Multiplicities (J): [2 1]
Distinct Distances (d): [0.70710678 1.41421356]
\end{verbatim}

The returned \texttt{quality} is the \(\Phi_{q, \text{intensive}}\)
value. Lower is better. The array \texttt{J} contains the multiplicities
for each distinct distance, and \texttt{d} contains the distinct
distances found in the design. There are three points, resulting in
three distances, \(d_{12}\), \(d_{13}\), and \(d_{23}\). The smallest
distance is \(d_{12} = d_{23}\), and it occurs twice, so \(J_1 = 2\).
The largest distance is \(d_{13}\), and it occurs once, so \(J_2 = 1\).

\end{example}

\begin{example}[Computation of the
Improvement]\protect\hypertarget{exm-mm_improvement}{}\label{exm-mm_improvement}

We will simulate a poor design in 2D and calculate the Morris-Mitchell
criterion \(\Phi_{q, \text{intensive}}\). A reduction in the
\(\Phi_{q, \text{intensive}}\) value can be observed. \texttt{spotoptim}
provides the function \texttt{mm\_improvement}, which allows an
efficient computation of the improvement.
Figure~\ref{fig-mm_improvement} shows the result: first, ten clustered
points are generated, which form a poor design. Then, a new, random
point is added and the Morris-Mitchell criterion is determined. The
improvement is computed as the difference between the Morris-Mitchell
criterion of the old and the new design. The improvement is positive,
indicating that the new design is better.

\begin{figure}

\centering{

\pandocbounded{\includegraphics[keepaspectratio]{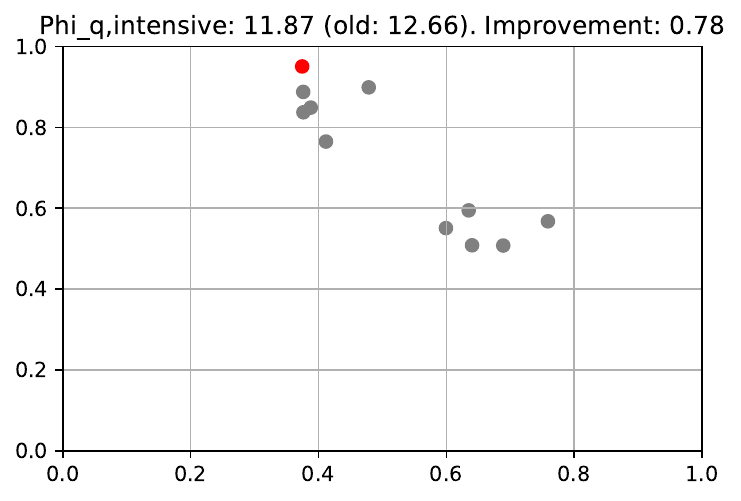}}

}

\caption{\label{fig-mm_improvement}Effect of adding a new point to a
poor design (grey). The new point is marked in red.}

\end{figure}%

\end{example}

We can also estimate the average improvement for adding a random point.

\begin{verbatim}
Average improvement: 1.0081 (stdev: 0.6462).
min: -10.7715.
max: 1.1975
lower quartile: 1.0684.
upper quartile: 1.1781
\end{verbatim}

\begin{figure}[H]

\centering{

\pandocbounded{\includegraphics[keepaspectratio]{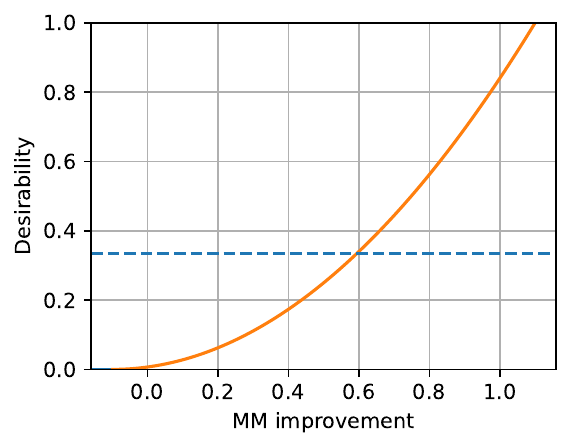}}

}

\caption{\label{fig-des-mm-imp}Desirability for Morris-Mitchell
Improvement}

\end{figure}

Since the estimated average improvement is $\approx 1$, we define the
desirability function as follows, see Figure~\ref{fig-des-mm-imp}.

\begin{Shaded}
\begin{Highlighting}[]
\NormalTok{d\_mm }\OperatorTok{=}\NormalTok{ DMax(}\OperatorTok{{-}}\FloatTok{0.1}\NormalTok{, }\FloatTok{1.1}\NormalTok{, scale}\OperatorTok{=}\DecValTok{2}\NormalTok{)}
\NormalTok{d\_mm.plot(xlabel}\OperatorTok{=}\StringTok{"MM improvement"}\NormalTok{, ylabel}\OperatorTok{=}\StringTok{"Desirability"}\NormalTok{, figsize}\OperatorTok{=}\NormalTok{(}\DecValTok{4}\NormalTok{,}\DecValTok{3}\NormalTok{))}
\end{Highlighting}
\end{Shaded}

Similar to the desirability of the space-fillingness criterion, we
estimate the expected improvement for the objective function. In this
example, we have chosen the Ackley function which is available in the
\texttt{spotoptim} package. For the definition of the desirability
function, we need the so far known best value of the objective function.
Then, we randomly generate 1\_000 points and compute the objective
function value for each point. The best estimated average improvement is
calculated as the difference between the known best value and the
minimum value of the 1\_000 samples of the objective function values.

\phantomsection\label{ackley}
\begin{verbatim}
Known best value: 4.24266557852776. Best improvement: 4.058484614461989
\end{verbatim}

\begin{figure}

\centering{

\pandocbounded{\includegraphics[keepaspectratio]{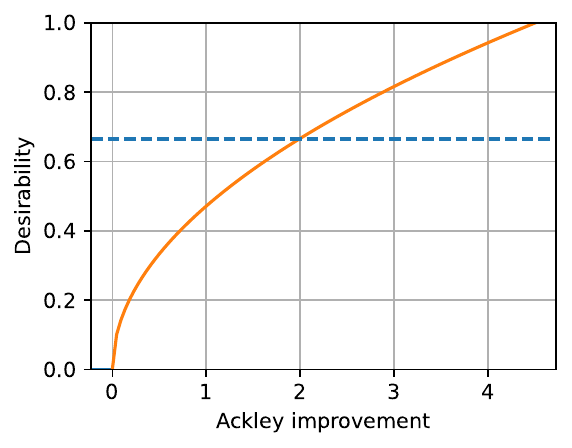}}

}

\caption{\label{fig-ackley}Desirability for Ackley Improvement}

\end{figure}

The best improvement of the Ackley function value is approximately 4.
This value is used for the defintion of the second desirability
function, which is shown in Figure~\ref{fig-ackley}.

\begin{tcolorbox}[enhanced jigsaw, colbacktitle=quarto-callout-note-color!10!white, opacityback=0, title=\textcolor{quarto-callout-note-color}{\faInfo}\hspace{0.5em}{Note:}, opacitybacktitle=0.6, arc=.35mm, left=2mm, bottomrule=.15mm, breakable, toprule=.15mm, coltitle=black, colframe=quarto-callout-note-color-frame, colback=white, rightrule=.15mm, bottomtitle=1mm, toptitle=1mm, titlerule=0mm, leftrule=.75mm]

\begin{itemize}
\tightlist
\item
  We consider in both dimensions the improvement with respect to an
  existing values, i.e., we are considering maximization problems.
\item
  Most optimizers use minimization, so we will pass
  \texttt{1\ \ -\ desirability} to the optimizer.
\end{itemize}

\end{tcolorbox}

The overall desirability \texttt{d\_overall} is set up as before. We are
using \texttt{SpotOptim} to optimize the desirability function. In
addition to the setup shown in the previous examples, we have introduced
a counter that measures how often the desirability function was called.
This counter is used to switch between the two desirability functions
after a certain number of calls:

\begin{itemize}
\tightlist
\item
  during the first phase, exploitation and exploration are combined,
  i.e., we use the overall desirability function.
\item
  during the second phase, exploitation is used, i.e., we use the
  desirability function for the Ackley improvement.
\end{itemize}

\begin{figure}

\centering{

\pandocbounded{\includegraphics[keepaspectratio]{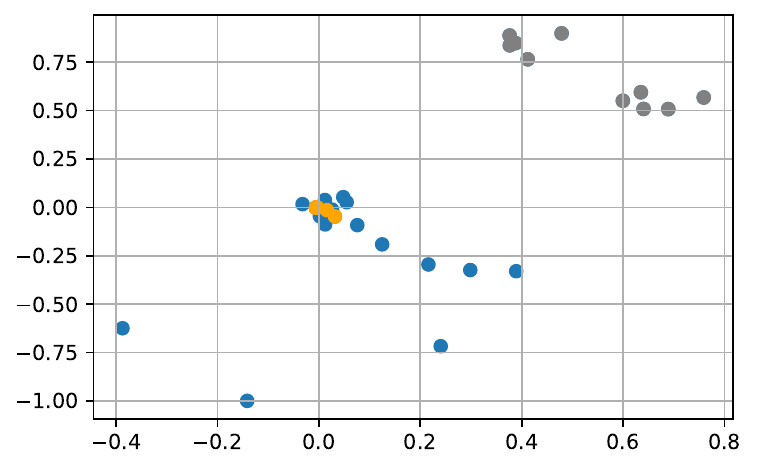}}

}

\caption{\label{fig-design-bi-objective}Design space of bi-objective
optimization. The inital points are shown in grey, the points that were
selected by the optimizer are shown in blue, when the Ackley-improvement
and space-fillingness improvement are used, and in orange, when the only
Ackley-improvement desirability is used.}

\end{figure}

The switch is done after ten calls to the desirability function. It is
clearly visible in Figure~\ref{fig-design-bi-objective}. This figure
shows the design space of the bi-objective optimization problem. The
original design points are shown in grey, the points that were selected
by the optimizer are shown in blue, when the Ackley-improvement and
space-fillingness improvement are used, and in orange, when the only
Ackley-improvement desirability is used.

\begin{verbatim}
Best x: [ 0.01872087 -0.01507067]
Best f(x): 0.022010382344720414
\end{verbatim}

\begin{figure}

\centering{

\pandocbounded{\includegraphics[keepaspectratio]{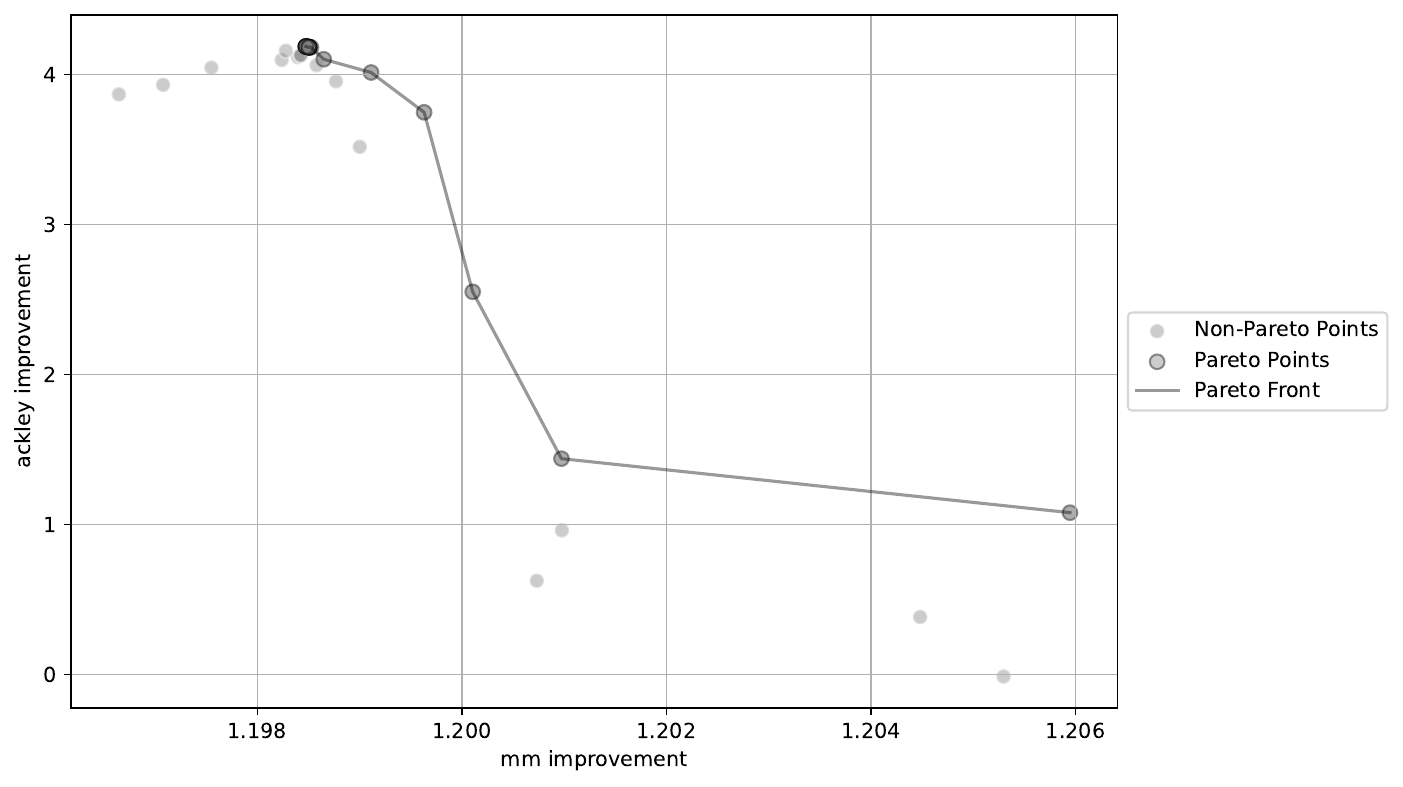}}

}

\caption{\label{fig-bi-objective}Bi-objective optimization. Improvement
of Ackley function versus Morris-Mitchell improvement.}

\end{figure}

Figure~\ref{fig-bi-objective} shows the Pareto front of the bi-objective
optimization problem.

\begin{verbatim}
Best point [ 0.01872087 -0.01507067] with value [0.08329742]
\end{verbatim}

\subsection{Infill-Point Diagnostic
Plots}\label{infill-point-diagnostic-plots}

\begin{figure}

\centering{

\pandocbounded{\includegraphics[keepaspectratio]{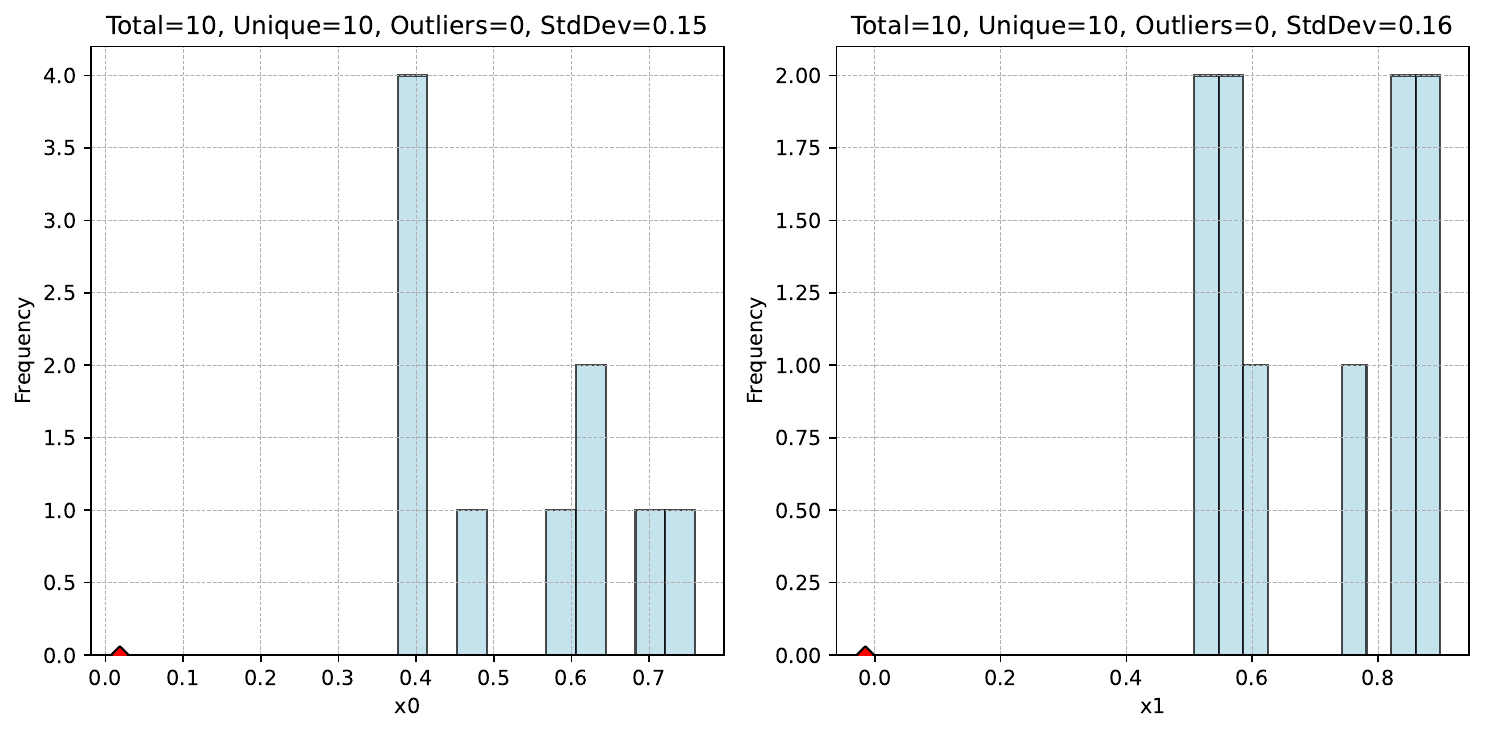}}

}

\caption{\label{fig-ip-histograms}Infill-point diagnostic plots for the
Morris-Mitchell function.}

\end{figure}

\begin{figure}

\centering{

\pandocbounded{\includegraphics[keepaspectratio]{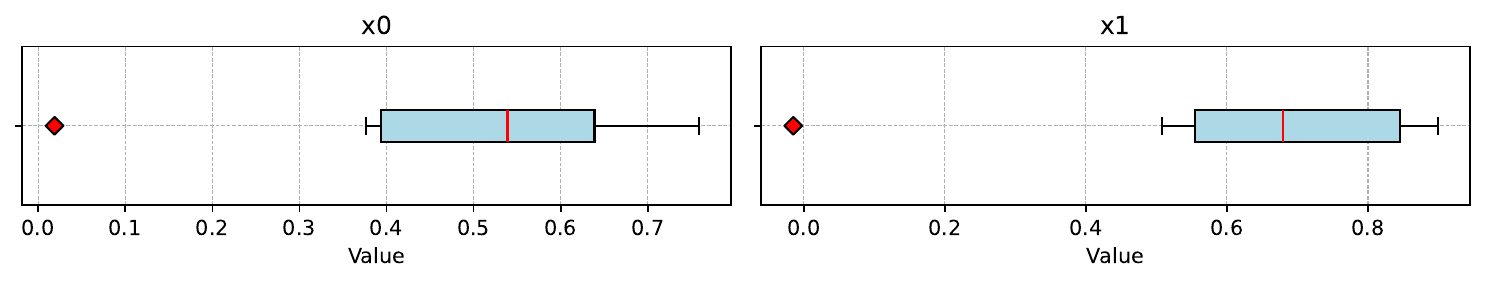}}

}

\caption{\label{fig-ip-boxplots}Infill-point diagnostic plots for the
Morris-Mitchell function.}

\end{figure}

Infill-point diagnostic plots are comprehensive tools to visualize the
location of the newly suggested best point in the context of the
existing data. We compare the distribution of the initial design points
with the chosen optimal point. The optimal point is marked in red. The
practitioner can use this plot to understand the distribution of the
initial design points and the location of the optimal point.
\texttt{spotoptim} provides two functions to create these plots:
\texttt{plot\_ip\_histograms} and \texttt{plot\_ip\_boxplots}.\\
Figure~\ref{fig-ip-histograms} shows the histograms of the initial
design points and the optimal point. Figure~\ref{fig-ip-boxplots} shows
the boxplots of the initial design points and the optimal point. In this
example, the optimal point is located in the center of the design space
(which is a known property of the Ackley function). In real-world
applications, the optimal point is unknown and the practitioner can use
these plots to understand the distribution of the initial design points
and the location of the optimal point.

\section{Conclusion}\label{sec-conclusion}

In this article, we have shown how to use the \texttt{spotdesitability}
package to perform three different multi-objective optimization tasks
using desirability functions: \gls{rsm} (Myers, Montgomery, and
Anderson-Cook 2016), surrogate model based optimization (Santner,
Williams, and Notz 2003), and hyperparameter tuning (Bartz et al. 2022).
The \texttt{spotdesirability} package is a \texttt{Python}
implementation of the \texttt{R} \texttt{desirability} package (Kuhn
2016). We have demonstrated how to define desirability functions for
different types of objectives, including maximization, minimization, and
target objectives. The intensified Morris-Mitchell criterion was
introduced and it was demonstrated how it can be used to handle the
exploration-exploitation trade-off using desirability in multi-objective
optimization.

Although the desirability function approach is one of the most widely
used methods in industry for the optimization of multiple response
processes (National Institute of Standards and Technology 2021), it is
rarely used in hyperparameter tuning. To fill this gap, we have shown
how to use the \texttt{spotdesirability} package in combination with the
\texttt{spotoptim} package to perform hyperparameter tuning using
desirability functions. The \texttt{spotoptim} package provides a
convenient way to perform surrogate model based optimization, and the
\texttt{spotdesirability} package allows us to define desirability
functions for different types of objectives. First results are
promising, but more research is needed to evaluate the performance of
the desirability function approach in hyperparameter tuning.

\section{Appendix}\label{appendix}

\subsection{\texorpdfstring{The \texttt{Python} Packages Used in This
Article}{The Python Packages Used in This Article}}\label{the-python-packages-used-in-this-article}

The following \texttt{Python} packages, classes, and functions are used
in this article:

\begin{Shaded}
\begin{Highlighting}[]
\ImportTok{import}\NormalTok{ os}
\ImportTok{from}\NormalTok{ math }\ImportTok{import}\NormalTok{ inf}
\ImportTok{import}\NormalTok{ warnings}
\ImportTok{import}\NormalTok{ numpy }\ImportTok{as}\NormalTok{ np}
\ImportTok{import}\NormalTok{ pandas }\ImportTok{as}\NormalTok{ pd}
\ImportTok{import}\NormalTok{ matplotlib.pyplot }\ImportTok{as}\NormalTok{ plt}
\ImportTok{import}\NormalTok{ pprint}

\ImportTok{from}\NormalTok{ sklearn.ensemble }\ImportTok{import}\NormalTok{ RandomForestRegressor}
\ImportTok{from}\NormalTok{ sklearn.datasets }\ImportTok{import}\NormalTok{ load\_diabetes}
\ImportTok{from}\NormalTok{ sklearn.gaussian\_process }\ImportTok{import}\NormalTok{ GaussianProcessRegressor}
\ImportTok{from}\NormalTok{ sklearn.gaussian\_process.kernels }\ImportTok{import}\NormalTok{ Matern, ConstantKernel}

\ImportTok{from}\NormalTok{ scipy.optimize }\ImportTok{import}\NormalTok{ dual\_annealing, minimize}

\ImportTok{from}\NormalTok{ spotdesirability.functions }\ImportTok{import}\NormalTok{ conversion\_pred, activity\_pred}
\ImportTok{from}\NormalTok{ spotdesirability.utils.desirability }\ImportTok{import}\NormalTok{ (DOverall, DMax, DCategorical, DMin,}
\NormalTok{                                                 DTarget, DArb, DBox)}
\ImportTok{from}\NormalTok{ spotdesirability.plot.ccd }\ImportTok{import}\NormalTok{ plotCCD}
\ImportTok{from}\NormalTok{ spotdesirability.functions.rsm }\ImportTok{import}\NormalTok{ rsm\_opt, conversion\_pred, activity\_pred}

\ImportTok{from}\NormalTok{ spotoptim.core.data }\ImportTok{import}\NormalTok{ SpotDataFromArray}
\ImportTok{from}\NormalTok{ spotoptim.core.experiment }\ImportTok{import}\NormalTok{ ExperimentControl}
\ImportTok{from}\NormalTok{ spotoptim.hyperparameters.parameters }\ImportTok{import}\NormalTok{ ParameterSet}
\ImportTok{from}\NormalTok{ spotoptim.function.torch\_objective }\ImportTok{import}\NormalTok{ TorchObjective}
\ImportTok{from}\NormalTok{ spotoptim.SpotOptim }\ImportTok{import}\NormalTok{ SpotOptim}
\ImportTok{from}\NormalTok{ spotoptim.nn.mlp }\ImportTok{import}\NormalTok{ MLP}
\ImportTok{from}\NormalTok{ spotoptim.sampling.mm }\ImportTok{import}\NormalTok{ mmphi\_intensive, mm\_improvement}
\ImportTok{from}\NormalTok{ spotoptim.eda }\ImportTok{import}\NormalTok{ plot\_ip\_boxplots, plot\_ip\_histograms}
\ImportTok{from}\NormalTok{ spotoptim.data.diabetes }\ImportTok{import}\NormalTok{ DiabetesDataset}
\ImportTok{from}\NormalTok{ spotoptim.function.mo.myer }\ImportTok{import}\NormalTok{ fun\_myer16a}
\ImportTok{from}\NormalTok{ spotoptim.plot.mo }\ImportTok{import}\NormalTok{ plot\_mo}
\ImportTok{from}\NormalTok{ spotoptim.plot.contour }\ImportTok{import}\NormalTok{ (mo\_generate\_plot\_grid,}
\NormalTok{                                     contourf\_plot)}
\ImportTok{from}\NormalTok{ spotoptim.sampling.mm }\ImportTok{import}\NormalTok{ plot\_mmphi\_vs\_n\_lhs}
\ImportTok{from}\NormalTok{ spotoptim.sampling.clustered }\ImportTok{import}\NormalTok{ Clustered}
\ImportTok{from}\NormalTok{ spotoptim.sampling.mm }\ImportTok{import}\NormalTok{ mmphi\_intensive, mm\_improvement}
\ImportTok{from}\NormalTok{ spotoptim.function }\ImportTok{import}\NormalTok{ ackley}

\NormalTok{warnings.filterwarnings(}\StringTok{"ignore"}\NormalTok{)}
\end{Highlighting}
\end{Shaded}

\subsection{Alternative Optimization Approach Using a Circular Design
Region}\label{alternative-optimization-approach-using-a-circular-design-region}

Kuhn (2016) also suggests alternatively maximizing desirability such
that experimental factors are constrained within a spherical design
region with a radius equivalent to the axial point distance:

\phantomsection\label{kuhn16a-optimization-circular-0}
\begin{Shaded}
\begin{Highlighting}[]
\CommentTok{\# Initialize the best result}
\NormalTok{best }\OperatorTok{=} \VariableTok{None}

\CommentTok{\# Perform optimization for each point in the search grid}
\ControlFlowTok{for}\NormalTok{ i, row }\KeywordTok{in}\NormalTok{ search\_grid.iterrows():}
\NormalTok{    initial\_guess }\OperatorTok{=}\NormalTok{ row.values  }\CommentTok{\# Initial guess for optimization}

    \CommentTok{\# Perform optimization using scipy\textquotesingle{}s minimize function}
\NormalTok{    result }\OperatorTok{=}\NormalTok{ minimize(}
\NormalTok{        rsm\_opt,}
\NormalTok{        initial\_guess,}
\NormalTok{        args}\OperatorTok{=}\NormalTok{(overallD, prediction\_funcs, }\StringTok{"circular"}\NormalTok{), }
\NormalTok{        method}\OperatorTok{=}\StringTok{"Nelder{-}Mead"}\NormalTok{,}
\NormalTok{        options}\OperatorTok{=}\NormalTok{\{}\StringTok{"maxiter"}\NormalTok{: }\DecValTok{1000}\NormalTok{, }\StringTok{"disp"}\NormalTok{: }\VariableTok{False}\NormalTok{\}}
\NormalTok{    )}

    \CommentTok{\# Update the best result if necessary}
    \CommentTok{\# Compare based on the negative desirability}
    \ControlFlowTok{if}\NormalTok{ best }\KeywordTok{is} \VariableTok{None} \KeywordTok{or}\NormalTok{ result.fun }\OperatorTok{\textless{}}\NormalTok{ best.fun:  }
\NormalTok{        best }\OperatorTok{=}\NormalTok{ result}
\BuiltInTok{print}\NormalTok{(}\StringTok{"Best Parameters:"}\NormalTok{, best.x)}
\BuiltInTok{print}\NormalTok{(}\StringTok{"Best Desirability:"}\NormalTok{, }\OperatorTok{{-}}\NormalTok{best.fun)}
\end{Highlighting}
\end{Shaded}

\begin{verbatim}
Best Parameters: [-0.50970524  1.50340746 -0.55595672]
Best Desirability: 0.8581520815997857
\end{verbatim}

Using these best parameters, the predicted values for conversion and
activity can be calculated as follows:

\begin{Shaded}
\begin{Highlighting}[]
\BuiltInTok{print}\NormalTok{(}\SpecialStringTok{f"Conversion pred(x): }\SpecialCharTok{\{}\NormalTok{conversion\_pred(best.x)}\SpecialCharTok{\}}\SpecialStringTok{"}\NormalTok{)}
\BuiltInTok{print}\NormalTok{(}\SpecialStringTok{f"Activity pred(x): }\SpecialCharTok{\{}\NormalTok{activity\_pred(best.x)}\SpecialCharTok{\}}\SpecialStringTok{"}\NormalTok{)}
\end{Highlighting}
\end{Shaded}

\begin{verbatim}
Conversion pred(x): 92.51922540231372
Activity pred(x): 57.499999903209876
\end{verbatim}

\begin{Shaded}
\begin{Highlighting}[]
\NormalTok{best\_temperature }\OperatorTok{=}\NormalTok{ best.x[}\DecValTok{1}\NormalTok{]}
\CommentTok{\# remove the temperature variable from the best parameters}
\NormalTok{best\_point }\OperatorTok{=}\NormalTok{ \{}\StringTok{"time"}\NormalTok{: best.x[}\DecValTok{0}\NormalTok{], }\StringTok{"temperature"}\NormalTok{: best.x[}\DecValTok{1}\NormalTok{], }\StringTok{"catalyst"}\NormalTok{: best.x[}\DecValTok{2}\NormalTok{]\}}
\CommentTok{\# set the values of temperature to the best temperature in the df}
\CommentTok{\# and recalculate the predicted values}
\NormalTok{plot\_grid\_best }\OperatorTok{=}\NormalTok{ plot\_grid\_best.copy()}
\NormalTok{plot\_grid\_best[}\StringTok{"temperature"}\NormalTok{] }\OperatorTok{=}\NormalTok{ best\_temperature}
\CommentTok{\# Recalculate the predicted values for conversion and activity}
\NormalTok{plot\_grid\_best[}\StringTok{"conversionPred"}\NormalTok{] }\OperatorTok{=}\NormalTok{ plot\_grid\_best[[}\StringTok{"time"}\NormalTok{, }\StringTok{"temperature"}\NormalTok{, }\StringTok{"catalyst"}\NormalTok{]].}\BuiltInTok{apply}\NormalTok{(}\KeywordTok{lambda}\NormalTok{ row: conversion\_pred(row.values), axis}\OperatorTok{=}\DecValTok{1}\NormalTok{)}
\NormalTok{plot\_grid\_best[}\StringTok{"activityPred"}\NormalTok{] }\OperatorTok{=}\NormalTok{ plot\_grid\_best[[}\StringTok{"time"}\NormalTok{, }\StringTok{"temperature"}\NormalTok{, }\StringTok{"catalyst"}\NormalTok{]].}\BuiltInTok{apply}\NormalTok{(}\KeywordTok{lambda}\NormalTok{ row: activity\_pred(row.values), axis}\OperatorTok{=}\DecValTok{1}\NormalTok{)}
\end{Highlighting}
\end{Shaded}

\begin{Shaded}
\begin{Highlighting}[]
\NormalTok{contourf\_plot(}
\NormalTok{    plot\_grid\_best,}
\NormalTok{    x\_col}\OperatorTok{=}\StringTok{"time"}\NormalTok{,}
\NormalTok{    y\_col}\OperatorTok{=}\StringTok{"catalyst"}\NormalTok{,}
\NormalTok{    z\_col}\OperatorTok{=}\StringTok{"conversionPred"}\NormalTok{,}
\NormalTok{    facet\_col}\OperatorTok{=}\StringTok{"temperature"}\NormalTok{,}
\NormalTok{    highlight\_point}\OperatorTok{=}\NormalTok{best\_point,}
\NormalTok{)}
\end{Highlighting}
\end{Shaded}

\begin{figure}[H]

\centering{

\pandocbounded{\includegraphics[keepaspectratio]{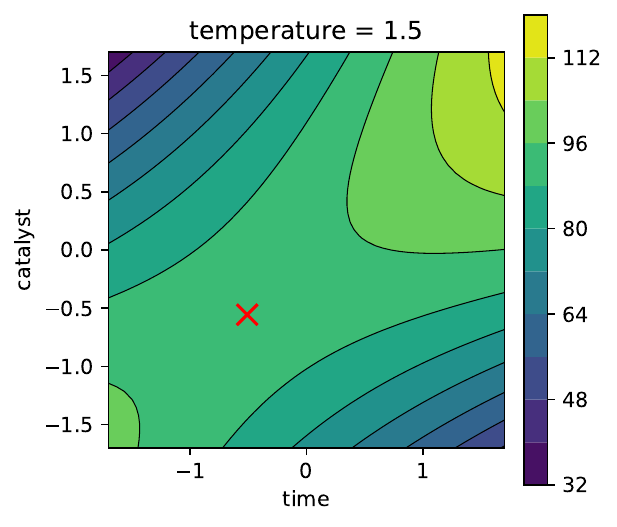}}

}

\caption{\label{fig-kuhn16a-best-conversion-circular}The response
surface for the percent conversion model. To plot the model contours,
the temperature variable was fixed at the best value found by the
optimizer.}

\end{figure}%

\begin{Shaded}
\begin{Highlighting}[]
\NormalTok{contourf\_plot(}
\NormalTok{    plot\_grid\_best,}
\NormalTok{    x\_col}\OperatorTok{=}\StringTok{"time"}\NormalTok{,}
\NormalTok{    y\_col}\OperatorTok{=}\StringTok{"catalyst"}\NormalTok{,}
\NormalTok{    z\_col}\OperatorTok{=}\StringTok{"activityPred"}\NormalTok{,}
\NormalTok{    facet\_col}\OperatorTok{=}\StringTok{"temperature"}\NormalTok{,}
\NormalTok{    highlight\_point}\OperatorTok{=}\NormalTok{best\_point,}
\NormalTok{)}
\end{Highlighting}
\end{Shaded}

\begin{figure}[H]

\centering{

\pandocbounded{\includegraphics[keepaspectratio]{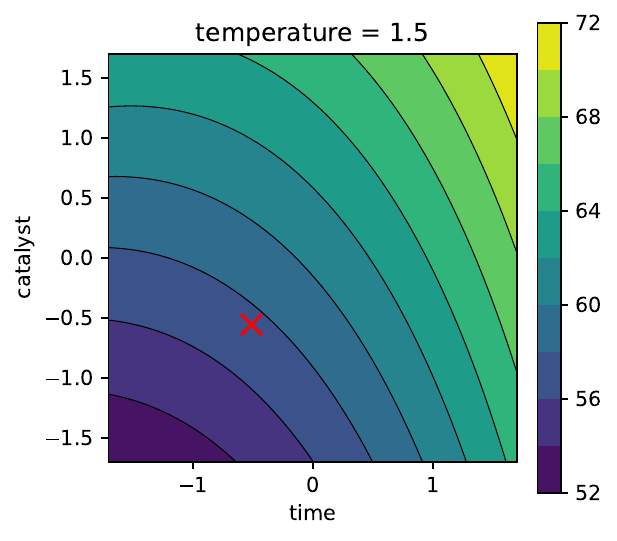}}

}

\caption{\label{fig-kuhn16a-best-activity-circular}The response surface
for the thermal activity model. To plot the model contours, the
temperature variable was fixed at the best value found by the
optimizer.}

\end{figure}%

Kuhn (2016) comments that the process converges to relative sub-optimual
values. He suggests that using a radius of 2 achieves an overall
desirability equal to one, even if the solution slightly extrapolates
beyond the design region.

\newpage{}

\newpage{}

\section*{References}\label{references}
\addcontentsline{toc}{section}{References}

\phantomsection\label{refs}

\end{document}